\documentclass[lettersize,journal]{IEEEtran}

\usepackage{amsmath,amsfonts}
\usepackage{array}
\usepackage[caption=false,font=normalsize,labelfont=sf,textfont=sf]{subfig}
\usepackage{textcomp}
\usepackage{stfloats}
\usepackage{url}
\usepackage{verbatim}
\usepackage{graphicx}
\usepackage{xcolor,float}
\usepackage{bm}
\usepackage{algorithmic,algorithm}
\usepackage{balance}

\newtheorem{definition}{Definition}
\newtheorem{assumption}{Assumption}
\newtheorem{theorem}{Theorem}
\newtheorem{lemma}{Lemma}
\newtheorem{remark}{Remark}

\def\BibTeX{{\rm B\kern-.05em{\sc i\kern-.025em b}\kern-.08em
    T\kern-.1667em\lower.7ex\hbox{E}\kern-.125emX}}
\usepackage{balance}
\begin{document}
\title{Dual-domain Defenses for Byzantine-resilient \\ Decentralized Resource Allocation}
\author{Runhua Wang, Qing Ling and Zhi Tian
\thanks{Runhua Wang and Qing Ling are with the School of Computer Science and Engineering and the Guangdong Provincial Key Laboratory of Computational Science, Sun Yat-Sen University, Guangzhou, Guangdong 510006, China. Zhi Tian is with the Department of Electrical and Computer Engineering, George Mason University, Fairfax, VA 22030, USA. Corresponding author: Qing Ling (lingqing556@mail.sysu.edu.cn).}
\thanks{Qing Ling is supported in part by NSF China grant 62373388, Guangdong Basic and Applied Basic Research Foundation grants 2021B1515020094 and 2023B1515040025, as well as R\&D project of Pazhou Lab (Huangpu) grant 2023K0606.
A short and preliminary version of this paper has been presented in ICASSP 2024 \cite{R-Wang-2024-icassp}.}
}


\maketitle

\begin{abstract}
This paper investigates the problem of decentralized resource allocation in the presence of Byzantine attacks. Such attacks occur when an unknown number of malicious agents send random or carefully crafted messages to their neighbors, aiming to prevent the honest agents from reaching the optimal resource allocation strategy. We characterize these malicious behaviors with the classical Byzantine attacks model, and propose a class of Byzantine-resilient decentralized resource allocation algorithms augmented with dual-domain defenses. The honest agents receive messages containing the (possibly malicious) dual variables from their neighbors at each iteration, and filter these messages with robust aggregation rules. Theoretically, we prove that the proposed algorithms can converge to neighborhoods of the optimal resource allocation strategy, given that the robust aggregation rules are properly designed. Numerical experiments are conducted to corroborate the theoretical results.
\end{abstract}


\section{INTRODUCTION}\label{sec 1}
\IEEEPARstart{R}ESOURCE allocation, which aims at assigning limited resources to a group of agents to minimize their costs, is a fundamental problem in network optimization.
Existing resource allocation algorithms can be categorized into distributed and decentralized approaches. Distributed resource allocation algorithms rely on a central agent to coordinate all the agents, which often leads to the communication bottleneck at the central agent, and thus results in limited scalability \cite{b-Jayash-Koshal-2011, b-Berkay-Turan-2021}. Consequently, decentralized resource allocation algorithms, which rely on coordination among neighboring agents, have become attractive alternatives. They have been widely applied in various fields, such as smart grids, transportation systems, wireless sensor networks, etc \cite{b-Zhong-Fan-2013, b-Md-Noor-A-Rahim-2022, b-Ayaz-Ahmad-2015}. Solving the decentralized resource allocation problem requires collaboration between neighboring agents. However, such collaboration is not always reliable since some of the agents could be malicious. The aim of this paper is to develop effective decentralized resource allocation algorithms that are resilient to the attacks from the malicious agents.


\noindent \textbf{Decentralized Resource Allocation Algorithms.}
Existing decentralized resource allocation algorithms can be categorized as continuous-time \cite{b-Shu-Liang-2020,b-Wenwen-Jia-2022,b-Yanan-Zhu-2021,b-Kaihong-Lu-2022} and discrete-time \cite{b-Hariharan-Lakshmanan-2008,b-LXiao-2006,b-Euhanna-Ghadimi-2013,b-Thinh-T.Doan-2021,b-Yun-Xu-2017,b-Jiaqi-Zhang-2020,b-AngeliaNedic-2018,b-SulaimanA-Alghunaim-2020}. In this paper, we focus on discrete-time algorithms. The primary challenge in algorithm design is to satisfy the global resource constraint. Weighted gradient methods have been proposed to guarantee global constraint satisfaction with the aid of feasible initialization \cite{b-Hariharan-Lakshmanan-2008,b-LXiao-2006,b-Euhanna-Ghadimi-2013}, but they turn out to be sensitive to perturbations.
The work of \cite{b-Hariharan-Lakshmanan-2008} is based on time-varying networks, while \cite{b-LXiao-2006} considers fixed networks. The work of \cite{b-Euhanna-Ghadimi-2013} utilizes historical information to accelerate the algorithm. On the other hand, primal-dual algorithms handle the global resource constraint via introducing a dual variable
\cite{b-Thinh-T.Doan-2021,b-Yun-Xu-2017,b-Jiaqi-Zhang-2020,b-AngeliaNedic-2018,b-SulaimanA-Alghunaim-2020}. The works of \cite{b-Thinh-T.Doan-2021,b-Yun-Xu-2017} develop decentralized Lagrangian methods, which precisely solve the primal sub-problems while perform a dual gradient step at each iteration. The work of \cite{b-Jiaqi-Zhang-2020} employs the push-pull gradient method to solve the dual problem and proposes a dual gradient tracking algorithm for unbalanced networks. For non-smooth resource allocation problems, decentralized proximal primal-dual algorithms are developed in \cite{b-AngeliaNedic-2018,b-SulaimanA-Alghunaim-2020}.

The decentralized resource allocation algorithms discussed above perform well when all the agents are honest. However, malicious agents, either spontaneously or by manipulation, are always threats to decentralized networks. These agents do not follow the given algorithmic protocol, but send random or crafted messages to their honest neighbors for the sake of misleading the optimization process. To characterize such behaviors, we use the classical Byzantine attacks model and term the malicious agents as Byzantine agents \cite{b-LeslieLamport-1982, b-ZhixiongYang-2020}. We briefly review some general Byzantine-resilient \textit{optimization} algorithms and few Byzantine-resilient \textit{resource allocation} algorithms, as follows.

\noindent \textbf{Byzantine-resilient Algorithms.} In the context of distributed optimization, Byzantine-resilient algorithms have been extensively studied. The main idea behind the algorithm design is to use various robust aggregation rules, such as coordinate-wise median \cite{D-Yin-2018-icml}, Krum \cite{Y-Chen-2017-MACS, C-Xie-2018-arXiv} and geometric median \cite{P-Blanchard-2017-NIPS}, to filter out malicious messages. However, directly extending this idea into decentralized optimization often cannot guarantee consensus, and thus yields large optimization errors \cite{b-ZhaoxianWu-2023}.

Given a general Byzantine-resilient decentralized optimization problem, honest agents cooperate to reach a consensual optimal solution that minimizes their average cost function. This is different to the resource allocation problem, where the honest agents are expected to obtain different optimal solutions (namely, allocated resources). Some works focus on deterministic problems \cite{b-Lili-Su-2021,b-Zhixiong-Yang-2019, b-Shreyas-Sundaram-2019,b-Lili-Su-2020,b-Cheng-Fang-2022,b-Waseem-Abbas-2022,b-Kananart-Kuwaranancharoen-2020} and some others consider stochastic problems \cite{b-ZhaoxianWu-2023,b-LieHe-2022}. Their common feature is to let each honest agent aggregate possibly malicious messages (namely, optimization variables) received from its neighbors in a robust manner.

For Byzantine-resilient decentralized optimization problems with deterministic cost functions, when the optimization variable is a scalar, \cite{b-Lili-Su-2021,b-Zhixiong-Yang-2019} proposes the trimmed mean (TM) robust aggregation rule, with which each honest agent discards the smallest $b$ and the largest $b$ messages received from its neighbors, followed by averaging the remaining messages and its own. Here $b$ is an estimated upper bound of the number of Byzantine neighbors. A similar approach in \cite{b-Shreyas-Sundaram-2019} lets each honest agent filter $b$ received messages larger and $b$ received messages smaller than its own message, also followed by averaging. For high-dimensional problems, \cite{b-Lili-Su-2020,b-Cheng-Fang-2022} extends TM to coordinate-wise TM (CTM), such that each honest agent performs the TM operation at each dimension. The work of \cite{b-Waseem-Abbas-2022} introduces the notion of centerpoint, which is an extension of the robust median aggregation rule to the high-dimensional scenario. In \cite{b-Kananart-Kuwaranancharoen-2020}, each iteration involves two filtering steps: distance-based and dimension-wise removals. Distance-based removal calculates the Euclidean distances between the received messages and the agent's own message, sorts the distances, and removes $b$ messages with the largest distances. Additionally, messages with extreme values in any dimension are removed.

When the cost functions are stochastic, TM and CTM are also applicable. Besides, the work of \cite{b-ZhaoxianWu-2023} proposes iterative outlier scissor (IOS), in which each honest agent iteratively discards $b$ messages that are the farthest from the average of the remaining received messages. The work of \cite{b-LieHe-2022} proposes self-centered clipping (SCC), in which each honest agent uses its own optimization variable as the center, clips the received messages, and then runs weighted average.

Although the aforementioned Byzantine-resilient decentralized optimization algorithms are proved to be effective, they cannot be directly applied to solve the resource allocation problem. The local optimization variables of the honest agents are coupled with a consensus constraint in the former but with a global resource constraint in the latter. Therefore, in a decentralized resource allocation algorithm, filtering ``outliers'' from the neighboring optimization variables becomes meaningless. To fill this gap, \cite{b-Berkay-Turan-2021} proposes a primal-dual Byzantine-resilient resource allocation algorithm from a robust optimization perspective, but the proposed algorithm is only applicable in a distributed network with a central server. A Byzantine-resilient decentralized resource allocation (BREDA) algorithm is developed in \cite{b-Runhua-Wang-2022}. In addition to the updates of primal and dual variables, each honest agent maintains an auxiliary variable that dynamically tracks the average of all honest agents' primal variables. Then, CTM is applied to aggregate the neighboring auxiliary variables.

\noindent \textbf{Our Contributions.} This paper focuses on the challenging and less-studied Byzantine-resilient decentralized resource allocation problem, and makes the following contributions:

\noindent \textbf{C1)} We propose a class of primal-dual Byzantine-resilient decentralized resource allocation algorithms with dual-domain defenses. The key intuition is that the honest agents should reach a consensual dual variable. Therefore, we can let each honest agent filter the received neighboring dual variables with properly designed robust aggregation rules, including but not limited to CTM, IOS and SCC.

\noindent \textbf{C2)} Compared with BREDA that defends against Byzantine attacks in the primal domain \cite{b-Runhua-Wang-2022}, the proposed algorithms utilize dual-domain defenses, and have the following advantages: (i) maintaining less variables and simpler updates; (ii) allowing more general robust aggregation rules than CTM; (iii) being able to reach dual consensus.

\noindent \textbf{C3)} Theoretically, we prove that if the robust aggregation rules are properly designed, the proposed algorithms converge to neighborhoods of the optimal primal-dual pair, and the honest agents are guaranteed to reach consensus in the dual domain even at presence of Byzantine attacks. With numerical experiments, we verify Byzantine-resilience of the proposed algorithms and its advantages over BREDA.

Compared to the short, preliminary conference version \cite{R-Wang-2024-icassp}, this journal version has been significantly extended. We have included comprehensive derivations for the algorithm design, detailed theoretical analysis, as well as additional numerical experiments that deepen the insights presented in \cite{R-Wang-2024-icassp}. These extensions not only reinforce the theoretical foundation but also enhance the practical relevance of our proposed algorithms.


\noindent\textbf{Paper Organization:} This paper is organized as follows. In Section \ref{sec 2}, we formulate the decentralized resource allocation problem under Byzantine attacks. Section \ref{sec 3} proposes an attack-free decentralized resource allocation algorithm that operates in the dual domain, and shows its failure under Byzantine attacks. Section \ref{sec 4} further proposes a class of Byzantine-resilient decentralized resource allocation algorithms. Section \ref{sec 5} establishes convergence of the proposed Byzantine-resilient decentralized resource allocation algorithms. Numerical experiments are given in Section \ref{sec 6}. Section \ref{sec 7} summarizes this paper and discusses future research directions.

\noindent\textbf{Notation:} Throughout this paper, $(\cdot)^{ \top }$ stands for the transposition of a vector or a matrix, $\|\cdot\|$ stands for the $\ell_{2}$-norm of a vector or a matrix, $\|\cdot\|_{F}$ denotes the Frobenius norm of a matrix, and $\left \langle \cdot , \cdot \right \rangle $ represents the inner product of vectors. We define $\widetilde{\bm{1}}\in \mathbb{R}^{J} $ and $\bm{1}\in \mathbb{R}^{H} $ as all-one column vectors while $I\in \mathbb{R}^{H \times H}$ as an identity matrix, where $J$ is the number of all agents and $H$ is the number of honest agents.

\section{PROBLEM STATEMENT}\label{sec 2}
We consider a decentralized resource allocation problem that involves a network of autonomous agents. The network is modeled as an undirected, connected graph $\widetilde{\mathcal{G}} ( \mathcal{J},\widetilde{\mathcal{E}} )$ with the set of vertices $\mathcal{J}:=\left \{ 1, \cdots,J \right \} $ and the set of edges $\widetilde{\mathcal{E}} $. If $( i,j) \in \widetilde{\mathcal{E}}$, then the two agents $i$ and $j$ are neighbors and can communicate with each other. For agent $i$, define the set of its neighbors as $\mathcal{N}_{i}= \{ j\mid ( i,j ) \in \widetilde{\mathcal{E}} \} $. Each agent $i$ possesses a strongly convex local cost function $f_{i}\left ( \bm{\theta}_{i} \right )$, where $\bm{\theta}_{i}  \in \mathbb{R}^{D}$ stands for the amount of local resources and belongs to a compact, convex set $C_{i}$. The average amount of local resources, denoted as $\frac{1}{J}\sum_{i \in \mathcal{J}}\bm{\theta}_{i}$,
equals to a constant vector $\bm{s} \in \mathbb{R}^{D}$. When all the agents are honest, the decentralized resource allocation problem is formulated as
\begin{align}
\label{eq_primal-problem}
\begin{split}
\underset{\widetilde{\bm{\Theta}}}{\min} \quad & \widetilde{f}(\widetilde{\bm{\Theta}})=\frac{1}{J}\sum_{i \in \mathcal{J}}f_{i}\left ( \bm{\theta}_{i} \right ),\\
s. t. \quad &   \frac{1}{J}\sum_{i \in \mathcal{J}}\bm{\theta}_{i} = \bm{s}, \quad \bm{\theta}_{i} \in C_{i}, ~ \forall i \in \mathcal{J}, \\
\end{split}
\end{align}
where $\widetilde{\bm{\Theta}}=[\bm{\theta}_{1},\cdots, \bm{\theta}_{J}] \in \mathbb{R}^{JD}$ concatenates all the local variables and $\widetilde{C}$ is the Cartesian product of $C_{i}$ for all $i\in \mathcal{J}$. 

The decentralized resource allocation problem in the form of \eqref{eq_primal-problem} arises in, for example, economic dispatch in smart grids \cite{b-Qiao-Li-2019, b-Huaqing-Li-2020}. The goal is to obtain an optimal generation strategy that minimizes the total generation cost, while satisfying a global power demand constraint and local generator constraints, through cooperation among a network of generators. We will introduce the economic dispatch problem in detail in Section \ref{sec 6}, and focus on the case that some of the generators are malicious.

\begin{figure}[htbp]
\centerline{\includegraphics[width=7.2cm]{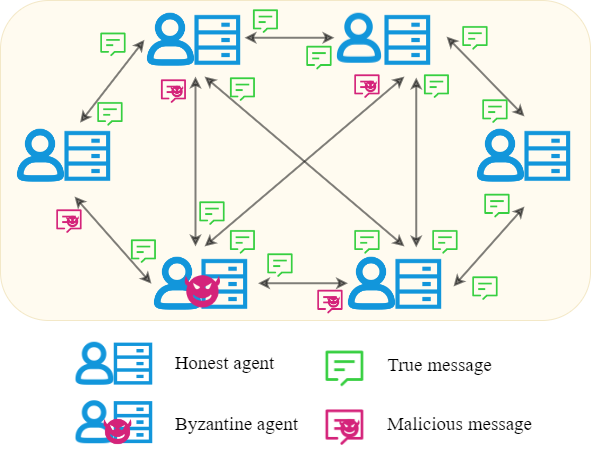}}
\caption{Decentralized resource allocation under Byzantine attacks.}
\label{fig_system}
\end{figure}

When some of the agents are Byzantine, as shown in Fig. \ref{fig_system}, solving \eqref{eq_primal-problem} is an impossible task, because they will not collaborate with the honest agents during the optimization process. Denote the set of Byzantine agents as $\mathcal{B}$ and the set of honest agents as $\mathcal{H}: =\mathcal{J}\setminus \mathcal{B}$. The numbers of Byzantine agents and honest agents are denoted as $B$ and $H$, respectively. Note that the number and identities of Byzantine agents are not known in advance, but we can roughly estimate an upper bound of the number. For notational convenience, we number the honest agents from $1$ to $H$, and the Byzantine agents from $H+1$ to $H+B$. Consider a subgraph $\mathcal{G}(\mathcal{H}, \mathcal{ \mathcal{E} })$ of $\widetilde{\mathcal{G}} ( \mathcal{J},\widetilde{\mathcal{ \mathcal{E} }})$, where $\mathcal{E} = \{(i,j) \in \widetilde{\mathcal{E}}; i, j \in \mathcal{H} \}$ is the set of edges between the honest agents. We assume $\mathcal{G}(\mathcal{H}, \mathcal{ \mathcal{E} })$ to be connected too so that the honest agents can cooperate. The goal of the honest agents is to solve
\begin{align}
\label{eq_oracle-problem}
\begin{split}
\underset{ \bm{\Theta}}{\min} \quad & f\left ( \bm{\Theta} \right ):=\frac{1}{ H }\sum_{i \in \mathcal{H}}f_{i}\left ( \bm{\theta}_{i} \right ),\\
s.t. \quad &  \frac{1}{ H }\sum_{i\in \mathcal{H}}\bm{\theta}_{i} = \bm{s}, \quad \bm{\theta}_{i} \in C_{i}, ~ \forall i \in \mathcal{H}, \\
\end{split}
\end{align}
where $\bm{\Theta}=[\bm{\theta}_{1},\cdots, \bm{\theta}_{H}] \in \mathbb{R}^{HD}$ concatenates all the local variables of the honest agents and $C$ is the Cartesian product of $C_{i}$ for all $i\in \mathcal{H}$.

In \eqref{eq_oracle-problem}, we modify the optimization objective to the average cost of the honest agents and consider the average resource constraint of the honest agents.
We do not modify the average resource supply $s$, as the Byzantine agents may still occupy some resources. Adjusting $s$ will inevitably affect the resources allocated to the honest agents.

However, solving \eqref{eq_oracle-problem} is still challenging since the honest agents cannot distinguish their Byzantine neighbors, while the latter can send arbitrarily malicious messages during the optimization process. Therefore, in this paper, we focus on developing Byzantine-resilient decentralized resource allocation algorithms to approximately solve \eqref{eq_oracle-problem}.

\section{ATTACK-FREE DECENTRALIZED RESOURCE ALLOCATION}\label{sec 3}
This section begins with reviewing an attack-free decentralized resource allocation algorithm, which operates in the dual domain, to solve \eqref{eq_primal-problem}.

\vspace{-2mm}
\subsection{Algorithm Development}
The Lagrangian function of (\ref{eq_primal-problem}) is
\begin{align}
\label{eq_Lagrangian-primal-problem}
\begin{split}
\widetilde{\mathcal{L}}( \widetilde{\bm{\Theta}} ;  \widetilde{\bm{\lambda }}) :=
\frac{1}{J}\sum_{i \in \mathcal{J}}f_{i}\left ( \bm{\theta}_{i} \right )+ \widetilde{\bm{\lambda }}^{\top} (\frac{1}{J}\sum_{i \in \mathcal{J}}\bm{\theta}_{i}-\bm{s}),
\end{split}
\end{align}
where $\widetilde{\bm{\lambda }}\in \mathbb{R}^D$ is the dual variable. Hence, the dual function $\widetilde{d}(\widetilde{\bm{\lambda }}) := \min_{\widetilde{\bm{\Theta}}\in \widetilde{C}} \widetilde{\mathcal{L}}( \widetilde{\bm{\Theta}} ;  \widetilde{\bm{\lambda }})$ is given by
\begin{align}
\label{eq_dual-function_primal}
\widetilde{d}(\widetilde{\bm{\lambda }}):=&\underset{\widetilde{\bm{\Theta}}\in \widetilde{C}}{\min}\{\frac{1}{J}\sum_{i \in \mathcal{J}}f_{i}\left ( \bm{\theta}_{i} \right )+\widetilde{\bm{\lambda }}^{\top}( \frac{1}{J}\sum_{i \in \mathcal{J}}\bm{\theta}_{i}-\bm{s})\} \\
=& \frac{1}{J}\sum_{i \in \mathcal{J}} \underset{\bm{\theta}_{i}\in C_{i}}{\min}\{f_{i}(\bm{\theta}_{i})+\widetilde{\bm{\lambda }}^{\top}\bm{\theta}_{i}\}-\widetilde{\bm{\lambda }}^{\top}\bm{s}\notag \\
=& \frac{1}{J}\sum_{i \in \mathcal{J}}(-\underset{\bm{\theta}_{i}\in C_{i}}{\max}\{-f_{i}(\bm{\theta}_{i})-\widetilde{\bm{\lambda }}^{\top}\bm{\theta}_{i}\})-\widetilde{\bm{\lambda }}^{\top}\bm{s}\notag \\
=& \frac{1}{J} \sum_{i \in \mathcal{J}} -\widetilde{F}_{i}^{*}(-\widetilde{\bm{\lambda }})-\widetilde{\bm{\lambda }}^{\top}\bm{s},\notag
\end{align}
where $\widetilde{F}_{i}^{*}(\widetilde{\bm{\lambda }}):=\max_{\bm{\theta}_{i}\in C_{i}}\{\widetilde{\bm{\lambda }}^{\top}\bm{\theta}_{i}-f_{i}(\bm{\theta}_{i})\}$. With it, we write the dual problem of \eqref{eq_primal-problem} as a minimization problem in the form of
\begin{align}
\label{eq_dual-problem_primal}
\underset{\widetilde{\bm{\lambda }}\in \mathbb{R}^D}{\min} ~ \widetilde{g}(\widetilde{\bm{\lambda }})= -d(\widetilde{\bm{\lambda}})=\sum_{i \in \mathcal{J}}\widetilde{g}_{i}(\widetilde{\bm{\lambda }}),
\end{align}
where $\widetilde{g}_{i}(\widetilde{\bm{\lambda }}):=\frac{1}{J}\widetilde{F}_{i}^{*}(-\widetilde{\bm{\lambda }})+\frac{1}{J}\widetilde{\bm{\lambda }}^{\top}\bm{s}$.

Because $f_{i}(\cdot)$ is strongly convex, according to the conjugate correspondence theorem in \cite{b-Amir-Beck-2017}, its conjugate function $\widetilde{F}_{i}^{*}(\cdot)$ is smooth. By Danskin's theorem \cite{b-Dimitri-P.-Bertsekas-1999}, the gradient $\nabla \widetilde{F}^{*}_{i}(\bm{\lambda}_{i})=\arg\max_{\bm{\theta}_{i}\in C_{i}}\{\bm{\lambda }_{i}^{\top}\bm{\theta}_{i}-f_{i}(\bm{\theta}_{i})\}$. Hence, we have
\begin{align}
\label{eq_gradient_dual_function_primal}
\nabla \widetilde{g}_{i}(\bm{\lambda}_{i})=\frac{1}{J}\bm{s}-\frac{1}{J}\arg\underset{\bm{\theta}_{i}\in C_{i}}{\min}\{\bm{\lambda}_{i}^{\top}\bm{\theta}_{i}+f_{i}(\bm{\theta}_{i})\}.
\end{align}

According to the above discussions, the optimization problem \eqref{eq_dual-problem_primal} can be solved through decentralized gradient methods \cite{b-Thinh-T.Doan-2021,b-AngeliaNedic-2009,b-B-Johansson-2008}. To do so, we let each agent holds a local dual variable $\bm{\lambda}_{i} \in \mathbb{R}^D$. The updates of primal and dual variables for all agents $i\in \mathcal{J}$ in the attack-free decentralized resource allocation algorithm at iteration $k+1$ are given by
\begin{align}
\bm{\theta}_{i}^{k}&=\arg\underset{\bm{\theta}_{i}\in C_{i}}{\min}\{\bm{\theta}_{i}^{\top}\bm{\lambda}_{i}^{k}+f_{i}(\bm{\theta}_{i})\}, \label{eq_DEGRA-1} \\
\bm{\lambda}_{i}^{k+\frac{1}{2}}&=\bm{\lambda}_{i}^{k}-\gamma^{k}\nabla \widetilde{g}_{i}(\bm{\lambda}_{i}^{k})=\bm{\lambda}_{i}^{k}-\gamma^{k}(\frac{1}{J}\bm{s}-\frac{1}{J}\bm{\theta}_{i}^{k}), \label{eq_DEGRA-2} \\
\bm{\lambda}_{i}^{k+1}&=\sum_{j\in \mathcal{N}_{i}\cup \{i\}}\widetilde{e}_{ij}\bm{\lambda}_{j}^{k+\frac{1}{2}}. \label{eq_DEGRA-3}
\end{align}
Therein, $\gamma^k > 0$ is the step size and $\widetilde{e}_{ij}\geq 0$ is the weight assigned by agent $i$ to agent $j$. Note that $\widetilde{e}_{ij}>0$ if and only if $(i, j)\in \widetilde{\mathcal{E}}$ or $i=j$. We collect these weights in $\widetilde{E}=[\widetilde{e}_{ij}]\in \mathbb{R}^{J\times J}$, which is assumed to be doubly stochastic. Such an attack-free decentralized resource allocation algorithm is summarized in Algorithm \ref{alg1}.

\begin{algorithm}[H]
\caption{Attack-free decentralized resource allocation algorithm}\label{alg1}
\begin{algorithmic}
\STATE
\STATE Initialization: All agents $i\in \mathcal{J}$ initialize $\bm{\lambda}_{i}^{0}=\bm{\lambda}^{0}$.
\FOR{$k=0,1,2,...$}
  \FOR{all agents $i\in \mathcal{J}$}
  \STATE Compute $\bm{\theta}_{i}^{k}=\arg\min_{\bm{\theta}_{i}\in C_{i}} \{\bm{\theta}_{i}^{\top}\bm{\lambda}_{i}^{k}+f_{i}(\bm{\theta}_{i})\}$.
  \STATE Compute $\bm{\lambda}_{i}^{k+\frac{1}{2}}=\bm{\lambda}_{i}^{k}-\gamma^{k}(\frac{1}{J}\bm{s}-\frac{1}{J}\bm{\theta}_{i}^{k})$.
  \STATE Broadcast $\bm{\lambda}_{i}^{k+\frac{1}{2}}$ to its neighbors.
  \STATE Receive $\bm{\lambda}_{j}^{k+\frac{1}{2}}$ from its neighbors.
  \STATE Aggregate $\bm{\lambda}_{i}^{k+1}=\sum_{j\in \mathcal{N}_{i}\cup \{i\}}\widetilde{e}_{ij}\bm{\lambda}_{j}^{k+\frac{1}{2}}$.
  \ENDFOR
\ENDFOR
\end{algorithmic}
\end{algorithm}

\subsection{Failure of Attack-free Decentralized Resource Allocation Algorithm under Byzantine Attacks}
When all the agents are honest, the decentralized resource allocation algorithm outlined in \eqref{eq_DEGRA-1}--\eqref{eq_DEGRA-3} can effectively solve \eqref{eq_primal-problem}; readers are referred to \cite{b-Thinh-T.Doan-2021,b-AngeliaNedic-2009,b-B-Johansson-2008}. However, it fails in the presence of Byzantine attacks. At iteration $k+1$, each honest agent $i\in \mathcal{H}$ updates $\bm{\lambda}_{i}^{k+1}$ based on $\bm{\lambda}_{i}^{k+\frac{1}{2}}$ from its own and $\bm{\lambda}_{j}^{k+\frac{1}{2}}$ from its neighbors $j \in \mathcal{N}_i$. An honest neighbor $j \in \mathcal{N}_i \cap \mathcal{H}$ faithfully sends the message $\bm{\lambda}_{j}^{k+\frac{1}{2}}$, but a Byzantine neighbor $j \in \mathcal{N}_i \cap \mathcal{B}$ may send an arbitrarily malicious message $*$ instead of the true message $\bm{\lambda}_{j}^{k+\frac{1}{2}}$. We define the message sent by agent $j$ as
\begin{align}
\label{broadcasting_message}
\check{\bm{\lambda}}_{j}^{k+\frac{1}{2}} =\left\{\begin{matrix}\bm{\lambda}_{j}^{k+\frac{1}{2}}, \quad & j\in \mathcal{H},
 \\ *, \quad & j \in \mathcal{B}.
\end{matrix}\right.
\end{align}

The malicious messages sent by the Byzantine agents prevent the honest agents from obtaining the optimal dual variable and corresponding resource allocation strategy. We provide a simple example to illustrate their impact. Assume that the local cost function of agent $i$ is $f_{i}(\bm{\theta}_{i})=\bm{\theta}_{i}^{2}$, the local resource constraint set is $C_{i}=[0,100]$, and the average resource is $\bm{s}=50$. The optimal dual variable and resource allocation of agent $i$ are $\bm{\lambda}_{i}^{*}=-100$ and $\bm{\theta}_{i}^{*}=50$, respectively. According to \eqref{eq_DEGRA-1}, the update of $\bm{\theta}_{i}^{k+1}$ is $\bm{\theta}_{i}^{k+1}=\Pi_{[0,100]}(-\frac{\bm{\lambda}_{i}^{k+1}}{2})$, the projection of $-\frac{\bm{\lambda}_{i}^{k+1}}{2}$ onto $[0,100]$. A Byzantine agent $j$ can manipulate $\bm{\lambda}_{i}^{k+1}$ by \eqref{eq_DEGRA-3} to be either $0$ or $-200$ through sending a proper $\check{\bm{\lambda}}_{j}^{k+\frac{1}{2}}$. In consequence, honest agent $i$ will obtain resource allocation of either $\bm{\theta}_{i}^{k+1}=0$ or $\bm{\theta}_{i}^{k+1}=100$, which are faraway from the optimal solution.

\section{BYZANTINE-RESILIENT DECENTRALIZED RESOURCE ALLOCATION}\label{sec 4}
In light of the influence of Byzantine attacks to decentralized resource allocation, we propose a class of Byzantine-resilient decentralized resource allocation algorithms to approximately solve \eqref{eq_oracle-problem} in this section.
\subsection{Algorithm Development}
As we have shown in Section \ref{sec 3}, the decentralized resource allocation algorithm outlined in \eqref{eq_DEGRA-1}--\eqref{eq_DEGRA-3} fails in the presence of Byzantine attacks. This is due to the vulnerability of the weighted average aggregation in \eqref{eq_DEGRA-3} to Byzantine attacks. To address this issue, we replace the weighted average with proper robust aggregation rules, and propose a class of Byzantine-resilient decentralized resource allocation algorithms. The updates of each honest agent $i \in \mathcal{H}$ are given by
\begin{align}
\bm{\theta}_{i}^{k}&=\arg\underset{\bm{\theta}_{i}\in C_{i}}{\min}\{\bm{\theta}_{i}^{\top}\bm{\lambda}_{i}^{k}+f_{i}(\bm{\theta}_{i})\},\label{eq_BRE-DEGRA-1}\\
\bm{\lambda}_{i}^{k+\frac{1}{2}}&=\bm{\lambda}_{i}^{k}-\gamma^{k}(\frac{1}{J}\bm{s}-\frac{1}{J}\bm{\theta}_{i}^{k}),\label{eq_BRE-DEGRA-2}\\
\bm{\lambda}_{i}^{k+1}&=AGG_{i}(\bm{\lambda}_{i}^{k+\frac{1}{2}},\{\check{\bm{\lambda}}_{j}^{k+\frac{1}{2}}\}_{j\in \mathcal{N}_{i}}),\label{eq_BRE-DEGRA-3}
\end{align}
where $AGG_{i}(\cdot)$ denotes a certain robust aggregation rule of honest agent $i$. The proposed Byzantine-resilient decentralized resource allocation algorithm is summarized in Algorithm \ref{alg2}.

In this paper, we mainly consider the applications of three well-appreciated robust aggregation rules: CTM, IOS and SCC. Further, we will show that a wide class of robust aggregation rules enable the updates of \eqref{eq_BRE-DEGRA-1}--\eqref{eq_BRE-DEGRA-3} to converge to neighborhoods of the optimal resource allocation strategy of \eqref{eq_oracle-problem}. The remaining design is to delineate the conditions for ``proper'' robust aggregation rules.

\noindent \textbf{Robust Aggregation Rules.} Intuitively, for an honest agent $i$, we expect that the output of $AGG_{i}(\bm{\lambda}_{i}^{k+\frac{1}{2}}, \{\check{\bm{\lambda}}_{j}^{k+\frac{1}{2}} \}_{j\in \mathcal{N}_{i}})$ is close to a proper weighted average of the messages from its honest neighbors and its own local dual variable, denoted as $\bar{\bm{\lambda}}_{i}^{k+\frac{1}{2}}$ $:=\sum_{j\in (\mathcal{N}_{i}\cap\mathcal{H})\cup \{i\}}e_{ij}\bm{\lambda}_{j}^{k+\frac{1}{2}}$ with the weights $\{e_{ij}\}_{j \in \mathcal{H}}$ satisfying $\sum_{j\in (\mathcal{N}_{i}\cap\mathcal{H})\cup \{i\}}e_{ij}=1$. We use the maximal value of $\{\|\bm{\lambda}_{j}^{k+\frac{1}{2}}-\bar{\bm{\lambda}}_{i}^{k+\frac{1}{2}}\| \}_{j\in (\mathcal{N}_{i}\cap \mathcal{H})\cup \{i\}}$ as the metric to quantify the proximity. Therefore, we follow \cite{b-ZhaoxianWu-2022, b-Haoxiang-Ye-2023} to characterize a set of robust aggregation rules with a virtual weight matrix and a contraction constant.

\begin{definition}
\label{d1}
Consider a set of robust aggregation rules $\{AGG_i\}_{i \in \mathcal{H}}$. If there exist a constant $\rho \geq 0$ and a matrix $E\in \mathbb{R}^{H \times H}$ whose elements satisfy $e_{ij} \in (0,1]$ when $j\in (\mathcal{N}_{i}\cap\mathcal{H})\cup \{i\}$, $e_{ij}=0$ when $j\notin (\mathcal{N}_{i}\cap\mathcal{H})\cup \{i\}$, and $\sum_{j\in (\mathcal{N}_{i}\cap\mathcal{H})\cup \{i\}}e_{ij}=1$ for any $i \in \mathcal{H}$, such that it holds
\begin{align}
\label{eq_definition1}
    & \|AGG_{i}(\bm{\lambda}_{i},\{\check{\bm{\lambda}}_{j}\}_{j\in \mathcal{N}_{i}})-\bar{\bm{\lambda}}_{i}\| \\
\le & \rho \max_{j\in (\mathcal{N}_{i}\cap \mathcal{H})\cup \{i\}}\|\bm{\lambda}_{j}-\bar{\bm{\lambda}}_{i}\| \nonumber
\end{align}
for any $i\in \mathcal{H}$, then $\rho$ is the contraction constant and $E$ is the virtual weight matrix associated with the set of robust aggregation rules $\{AGG_i\}_{i \in \mathcal{H}}$. Here $\bar{\bm{\lambda}}_{i}:=\sum_{j\in (\mathcal{N}_{i}\cap\mathcal{H})\cup \{i\}}e_{ij}\bm{\lambda}_{j}$.
\end{definition}

In the next section, we will prove that if a robust aggregation rule satisfies Definition \ref{d1}, it is ``proper'' if the associated $\rho$ is small and $E$ is doubly stochastic.

\begin{remark}
The work of \cite{b-Haoxiang-Ye-2023} has demonstrated that CTM, IOS and SCC all satisfy Definition \ref{d1} under network conditions stricter than connectedness of the honest agents, and specified their corresponding $\rho$ and $E$. For example, for each honest agent, CTM requires the number of its honest neighbors to exceed $2b$. Note that the pair of $(\rho,E)$ is not unique. Finding the best pair is beyond the scope of this paper.

There also exist other robust aggregation rules, such as the total-variation-based \cite{b-JiePeng-2021} and attack-detection-based \cite{b-Kailkhura-2016} ones, which do not satisfy Definition 1. We will investigate these approaches in our future work.
\end{remark}

\begin{algorithm}[H]
\caption{Byzantine-resilient decentralized resource allocation algorithm}\label{alg2}
\begin{algorithmic}
\STATE
\STATE Initialization: All agents $i$ initialize $\bm{\lambda}_{i}^{0}=\bm{\lambda}^{0}$.
\FOR{$k=0,1,2,...$}
  \FOR{all honest agents $i\in \mathcal{H}$}
  \STATE Compute $\bm{\theta}_{i}^{k}=\arg\min_{\bm{\theta}_{i}\in C_{i}}\{\bm{\theta}_{i}^{\top}\bm{\lambda}_{i}^{k}+f_{i}(\bm{\theta}_{i})\}$.
  \STATE Compute $\bm{\lambda}_{i}^{k+\frac{1}{2}}=\bm{\lambda}_{i}^{k}-\gamma^{k}(\frac{1}{J}\bm{s}-\frac{1}{J}\bm{\theta}_{i}^{k})$.
  \STATE Broadcast $\bm{\lambda}_{i}^{k+\frac{1}{2}}$ to its neighbors.
  \STATE Receive $\check{\bm{\lambda}}_{j}^{k+\frac{1}{2}}$ from its neighbors.
  \STATE Aggregate $\bm{\lambda}_{i}^{k+1}=AGG_{i}(\bm{\lambda}_{i}^{k+\frac{1}{2}}, \{\check{\bm{\lambda}}_{j}^{k+\frac{1}{2}}\}_{j\in \mathcal{N}_{i}})$.
  \ENDFOR
  \FOR{all Byzantine agents $i\in \mathcal{B}$}
  \STATE Broadcast $\check{\bm{\lambda}}_{i}^{k+\frac{1}{2}}=*$ to its neighbors.
  \ENDFOR
\ENDFOR
\end{algorithmic}
\end{algorithm}

\subsection{Advantages over BREDA}
Our proposed algorithms have several advantages over BREDA \cite{b-Runhua-Wang-2022}: simplicity, generality and dual consensus. First, at each iteration of BREDA, each honest agent needs to update a primal variable, a dual variable, and an auxiliary variable that tracks the average of the honest primal variables. By contrast, at each iteration of our proposed algorithms, each honest agent only updates two local variables, one is primal and the other is dual. Second, the robust aggregation rule of BREDA is confined to CTM; using other robust aggregation rules lacks convergence guarantee. However, CTM does not fit for the scenario that an honest agent has a large number of Byzantine neighbors, because the number of discarded messages has to be at least twice. This is unfavorable especially when the underlying network is sparse. Instead, our proposed algorithms allow a wide class of robust aggregation rules that satisfy Definition \ref{d1}. Third, BREDA guarantees the local auxiliary variables to be nearly consensual, but the local dual variables are not necessarily so. We will validate this fact in the numerical experiments. Since the optimal dual variable stands for the shadow price of the resources \cite{b-FPKelly-1998}, reaching consensus of the local dual variables is important in various applications. Our proposed algorithms have such a guarantee, as shown in the next section.

\section{CONVERGENCE ANALYSIS}\label{sec 5}
This section analyzes convergence of the attack-free and Byzantine-resilient decentralized resource allocation algori- thms, outlined in \eqref{eq_DEGRA-1}--\eqref{eq_DEGRA-3} and \eqref{eq_BRE-DEGRA-1}--\eqref{eq_BRE-DEGRA-3}, respectively.

We begin with several assumptions.

\begin{assumption}
\label{a1}
For any $i\in \mathcal{J}$, the local cost function $f_{i}(\cdot)$ is $u_{f}$-strongly convex and $L_{f}$-smooth, and the local constraint set $C_{i}$ is compact and convex.
\end{assumption}

\begin{assumption}
\label{a2}
There exist $\widetilde{\bm{\Theta}}$ and $\bm{\Theta}$ in the relative interiors of $\widetilde{C}$ and $C$, such that the constraints $\frac{1}{J}\sum_{i \in \mathcal{J}}\bm{\theta}_{i}=\bm{s}$ and $\frac{1}{H}\sum_{i\in \mathcal{H}}\bm{\theta}_{i}=\bm{s}$ satisfy, respectively.
\end{assumption}

Assumptions \ref{a1} and \ref{a2} are common in investigating resource allocation problems, and are satisfied by many applications \cite{b-Berkay-Turan-2021, b-Qiao-Li-2019, b-Huaqing-Li-2020}. With Assumptions \ref{a1} and \ref{a2}, the duality gaps of \eqref{eq_primal-problem} and \eqref{eq_oracle-problem} are both 0. In addition, the negative dual functions to minimize are also strongly convex and smooth.

\begin{assumption}
\label{a3}
The graphs $\widetilde{\mathcal{G}} ( \mathcal{J},\widetilde{\mathcal{ \mathcal{E} }} )$ and ${\mathcal{G}} ( \mathcal{J},{\mathcal{ \mathcal{E} }} )$ are both undirected and connected. The weight matrices $\widetilde{E}$ and $E$ are doubly stochastic and row stochastic, respectively, and satisfy
\begin{align}
\widetilde{\kappa} & := \|\widetilde{E}-\frac{1}{J}\widetilde{\bm{1}}\widetilde{\bm{1}}^{\top}\|^{2} < 1, \label{eq_a6_tttt} \\
\kappa & := \|E-\frac{1}{H}\bm{1}\bm{1}^{\top} E\|^{2} < 1. \label{eq_a6}
\end{align}
\end{assumption}

We have emphasized that the connectedness of $\widetilde{\mathcal{G}}$ and ${\mathcal{G}}$ is necessary. The requirement \eqref{eq_a6_tttt} is common in decentralized optimization. It holds when $\widetilde{e}_{ij}>0$ if and only if $(i, j)\in \widetilde{\mathcal{E}}$ or $i=j$. The requirement \eqref{eq_a6} on the associated virtual weight matrix $E$ is in the same form of \eqref{eq_a6_tttt} if $E$ is doubly stochastic, but we allow $E$ to be only row stochastic.

\subsection{Convergence of Attack-free Decentralized Resource Allocation Algorithm}
Denote $(\widetilde{\bm{\Theta}}^{*},\widetilde{\bm{\lambda}}^{*})$ as the optimal primal-dual pair of \eqref{eq_primal-problem}, in which $\widetilde{\bm{\Theta}}^{*}\in \mathbb{R}^{JD}$ and $\widetilde{\bm{\lambda}}^{*}\in \mathbb{R}^{D}$. The following theorem shows the convergence of the attack-free decentralized allocation algorithm \eqref{eq_DEGRA-1}--\eqref{eq_DEGRA-3}.

\begin{theorem}
\label{t1} Consider $\widetilde{\bm{\Theta}}^{k+1}$ and $\{\bm{\lambda}_{i}^{k+1}\}_{i\in \mathcal{J}}$ generated by the attack-free decentralized resource allocation algorithm \eqref{eq_DEGRA-1}--\eqref{eq_DEGRA-3} and suppose that no Byzantine agents are present. If Assumptions \ref{a1}--\ref{a3} hold, then with a proper decreasing step size $\gamma^{k}=O(\frac{1}{k})$, we have \\
$a)$ $\lim_{k \to +\infty} \sum_{i\in \mathcal{J}}\|\bm{\lambda}_{i}^{k+1}-\widetilde{\bm{\lambda}}^{*}\|=0$, \\
$b)$ $\lim_{k \to +\infty} \|\widetilde{\bm{\Theta}}^{k+1}-\widetilde{\bm{\Theta}}^{*}\|=0$.
\end{theorem}

Theorem \ref{t1} shows that the local primal and dual variables generated by \eqref{eq_DEGRA-1}--\eqref{eq_DEGRA-3} converge to their optima. This matches the classical conclusion for the decentralized gradient method \cite{b-Thinh-T.Doan-2021,b-AngeliaNedic-2009,b-B-Johansson-2008}. Those works assume convex and possibly non-smooth cost functions, while we assume strongly convex and smooth cost functions, with which we have performance guarantee for the ensuing Byzantine-resilient algorithms. The proof of Theorem \ref{t1} and the conditions on the step size $\gamma^{k}$ are presented in Appendix \ref{appendix1}.

\subsection{Convergence of Byzantine-resilient Decentralized Resource Allocation Algorithm}

Similarly, denote $({\bm{\Theta}}^{*},{\bm{\lambda}}^{*})$ as the optimal primal-dual pair of \eqref{eq_oracle-problem}, in which ${\bm{\Theta}}^{*}\in \mathbb{R}^{HD}$ and ${\bm{\lambda}}^{*}\in \mathbb{R}^{D}$. The following theorem shows the convergence of the Byzantine-resilient decentralized allocation algorithm \eqref{eq_BRE-DEGRA-1}--\eqref{eq_BRE-DEGRA-3}.
\begin{theorem}
\label{t2} Consider ${\bm{\Theta}}^{k+1}$ and $\{\bm{\lambda}_{i}^{k+1}\}_{i\in \mathcal{H}}$ generated by the Byzantine-resilient decentralized resource allocation algorithm \eqref{eq_BRE-DEGRA-1}--\eqref{eq_BRE-DEGRA-3}. Suppose that Byzantine agents are present but the used robust aggregation rule satisfies \eqref{eq_definition1} in Definition \ref{d1}. If Assumptions \ref{a1}--\ref{a3} hold and the contraction constant $\rho$ satisfies
$$\rho<\frac{1-\kappa}{8\sqrt{H}},$$
then with a proper decreasing step size $\gamma^{k}=O(\frac{1}{k})$, we have \\
$a)$ $\limsup_{k\to +\infty}\sum_{i\in \mathcal{H}}\|\bm{\lambda}_{i}^{k+1}-\bm{\lambda}^{*}\|\le \sqrt{ \frac{192\delta^{2}H^{2}}{\beta^{2}}}\cdot \sqrt{1+\frac{9}{\epsilon^{3}}}\cdot\sqrt{4\rho^{2}H+\chi^{2}},$\\
$b)$ $\lim_{k\to +\infty}\sum_{i\in \mathcal{H}}\|\bm{\lambda}_{i}^{k+1}-\bar{\bm{\lambda}}^{k+1}\|= 0,$\\
$c)$ $\limsup_{k\to +\infty}\|\bm{\Theta}^{k+1}-\bm{\Theta}^{*}\|\le \frac{1}{u_{f}}\cdot \sqrt{ \frac{192\delta^{2}}{\beta^{2}}}\cdot \sqrt{1+\frac{9}{\epsilon^{3}}}\cdot\sqrt{4\rho^{2}H+\chi^{2}},$\\
where $\bar{\bm{\lambda}}^{k+1}:=\frac{1}{H}\sum_{i\in \mathcal{H}}\bm{\lambda}_{i}^{k+1}$, $\beta=\frac{1}{H(u_{f}+L_{f})}$, $\epsilon=\kappa-8\rho\sqrt{H}$, and $\chi^{2}:=\frac{1}{H}\|E^{\top}\bm{1}-\bm{1}\|^{2}$ quantifies the non-doubly stochasticity of $E$.
\end{theorem}


The proof of Theorem \ref{t2} and the conditions on the step size $\gamma^{k}$ are presented in Appendix \ref{appendix2}. Theorem \ref{t2} demonstrates that if the robust aggregation rule is properly designed such that the associated contraction constant $\rho$ is sufficiently small, then the local primal and dual variables generated by \eqref{eq_BRE-DEGRA-1}--\eqref{eq_BRE-DEGRA-3} converge to neighborhoods of their optima. Sizes of the neighborhoods are determined by the associated contraction constant $\rho$ and virtual weight matrix $E$ (more precisely, $\chi^2$). Notably, the local dual variables are guaranteed to reach consensus even under Byzantine attacks.

Compared to the proof of Theorem \ref{t1}, that of Theorem \ref{t2} is more challenging. First, under the Byzantine attacks and with the robust aggregation rule, dual-domain consensus is no longer merited. We discover that $\rho$ must be sufficiently small for reaching consensus. Second, due to the imperfectness during the aggregation, each iteration incurs an error determined by $\rho$ and $\chi^2$. We have to handle such an error within the analysis. Note that when $\rho=0$ and $E$ is doubly stochastic, Theorem \ref{t2} reduces to Theorem \ref{t1}.

Our analysis is related to but significantly different from that in \cite{b-ZhaoxianWu-2023}. The work of \cite{b-ZhaoxianWu-2023} considers a general Byzantine-resilient decentralized stochastic non-convex \textit{optimization} problem, and analyzes robust aggregation rules that satisfy Definition \ref{d1} in the primal domain. By contrast, we consider a strongly convex \textit{resource allocation} problem, and analyze in the dual domain. The different assumptions lead to different convergence metrics, and the corresponding technical tools are different, too.

\begin{remark}
Although Algorithm \ref{alg2} and its convergence analysis in Theorem \ref{t2} are only applicable to \eqref{eq_primal-problem} with an equality constraint, they can be extended to handle the problem with an inequality constraint as well. The extensions can be achieved by incorporating a non-negative projection operation in updating the dual variables of Algorithm \ref{alg2} and utilizing the non-expansive property of projection in the convergence analysis of Theorem \ref{t2}.
\end{remark}

\section{NUMERICAL EXPERIMENTS}\label{sec 6}
In this section, we conduct numerical experiments to show the performance of the proposed Byzantine-resilient decen- tralized resource allocation algorithms.

\begin{figure*}[htbp]
\centerline{\includegraphics[width=16.6cm]{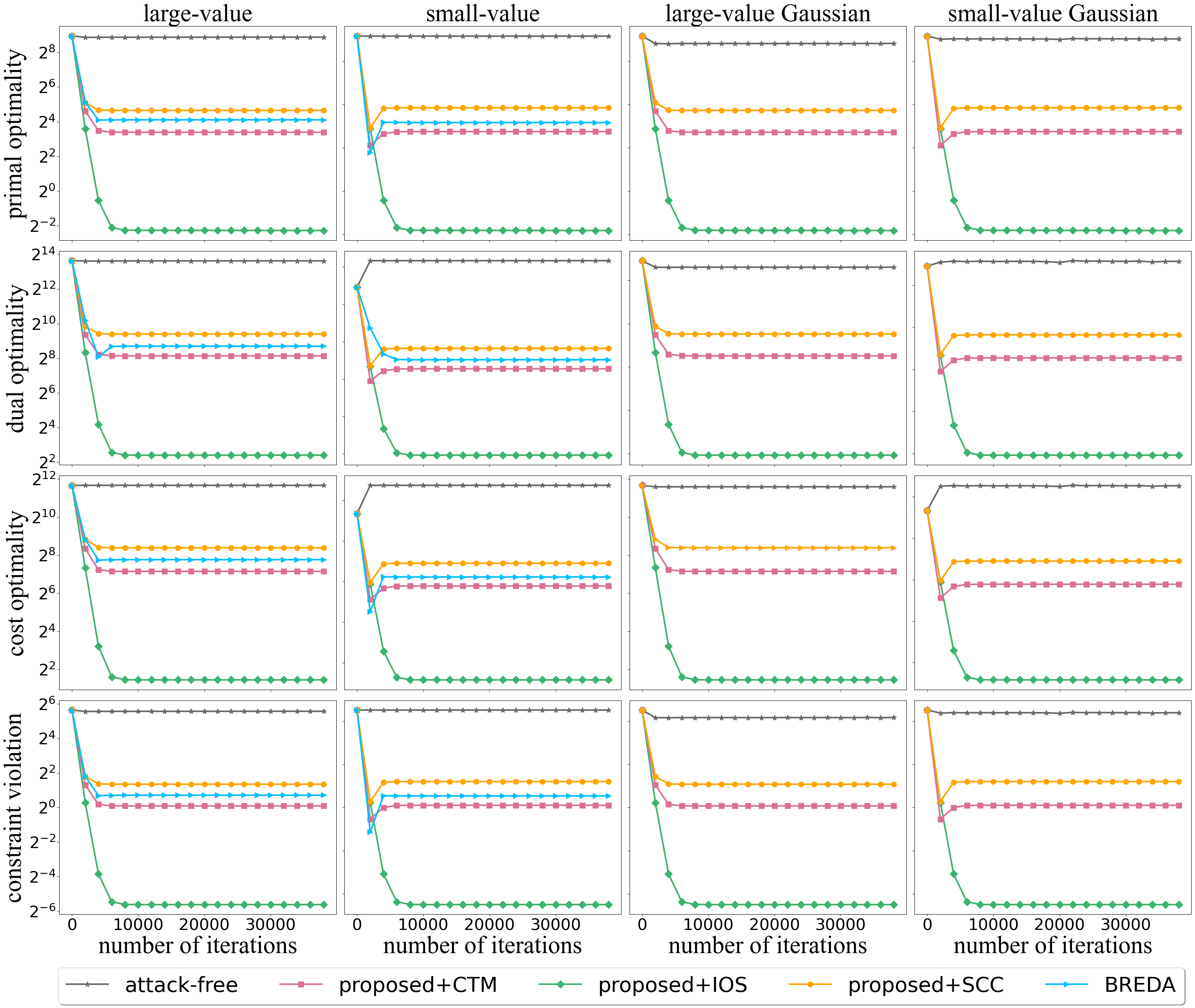}}
\caption{Primal optimality, dual optimality, cost optimality and constraint violation of the compared algorithms with optimal parameters in Case 1. The number of Byzantine agents is set as $6$.}
\label{fig_1}
\end{figure*}

\begin{figure*}[htbp]
\centerline{\includegraphics[width=16.8cm]{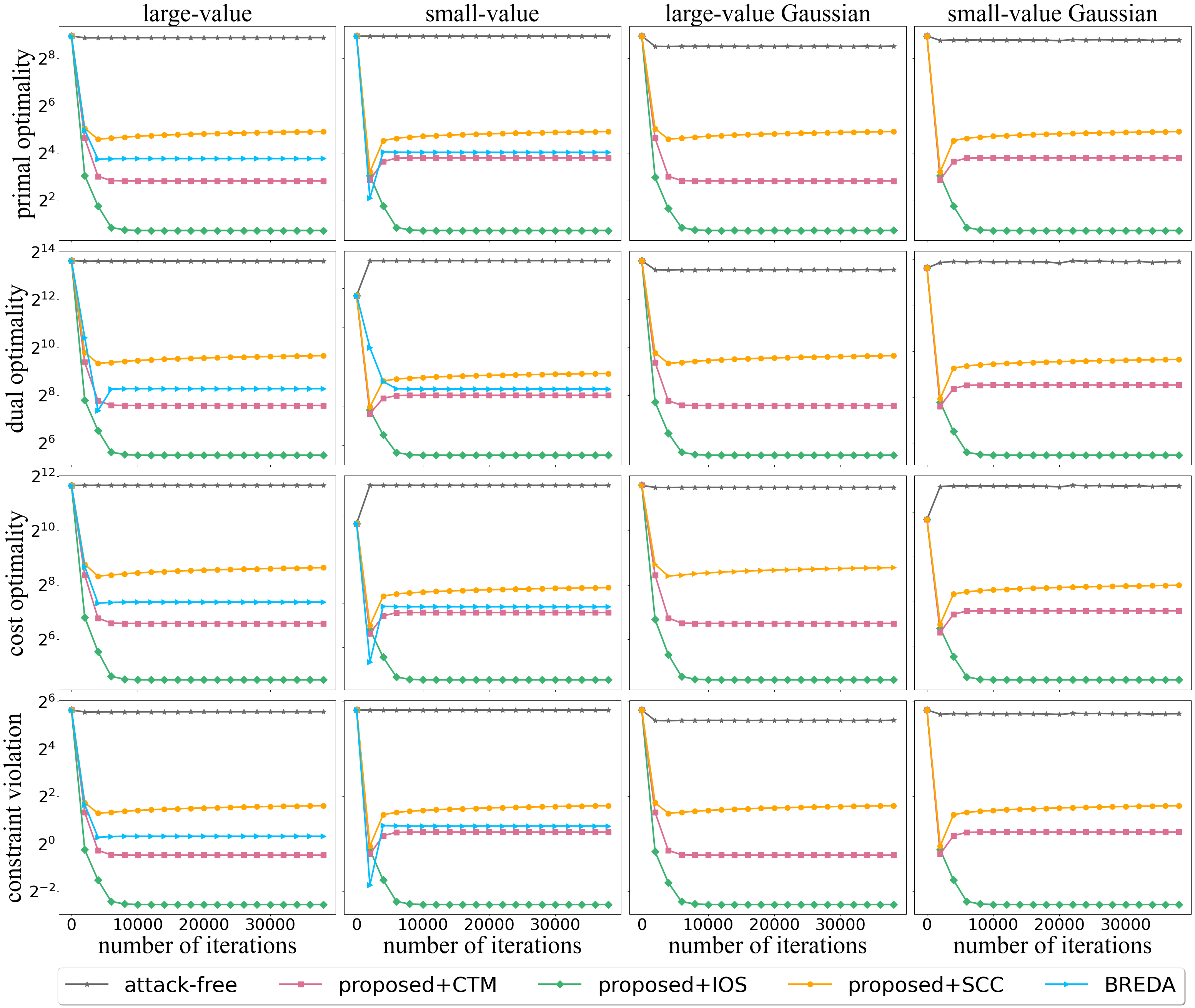}}
\caption{Primal optimality, dual optimality, cost optimality and constraint violation of the compared algorithms with non-optimal parameters in Case 1. The number of Byzantine agents is set as $6$.}
\label{fig_2}
\end{figure*}

\begin{figure*}[htbp]
\centerline{\includegraphics[width=16.8cm]{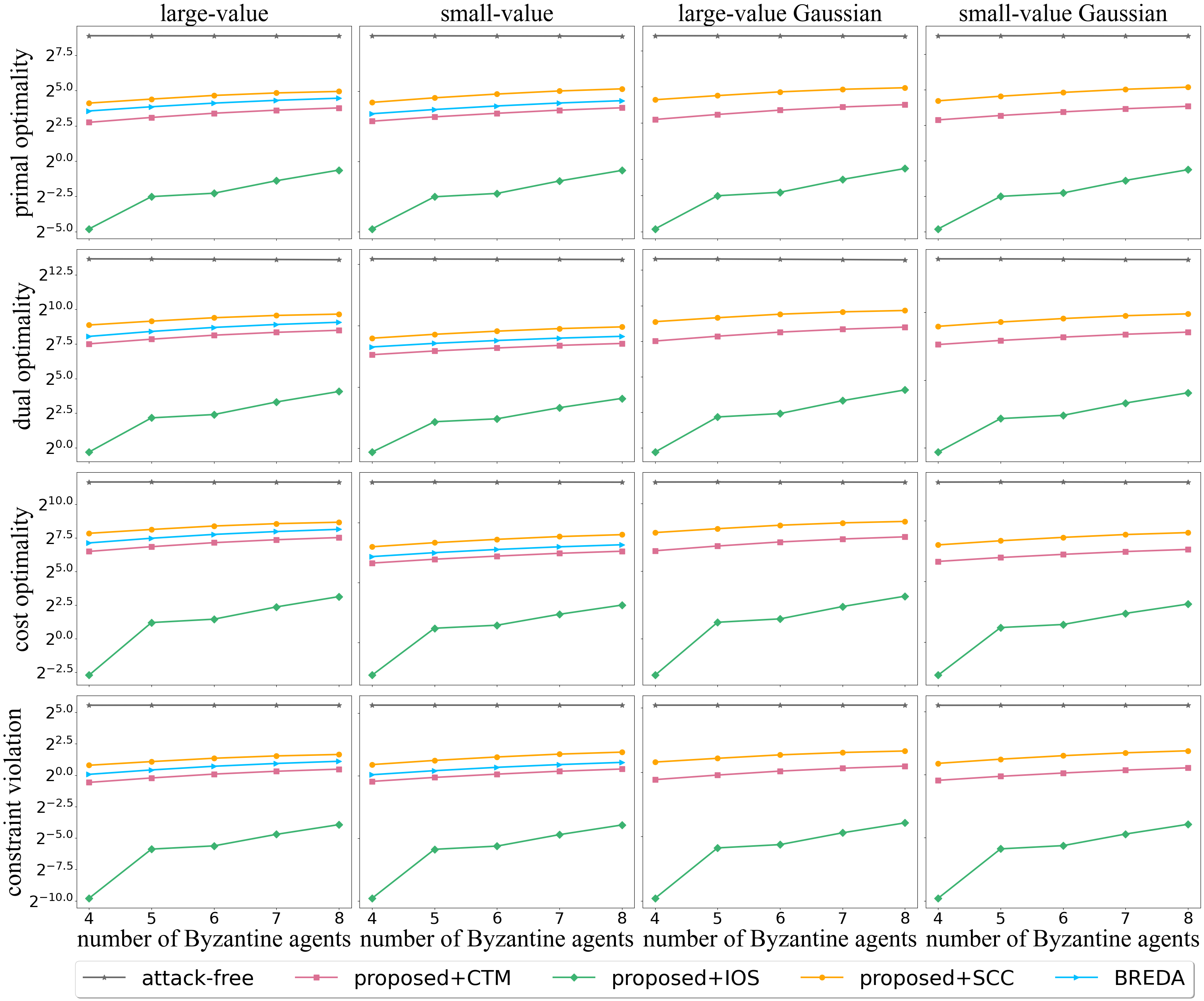}}
\caption{Primal optimality, dual optimality, cost optimality and constraint violation of the compared algorithms with optimal parameters in Case 1. The number of Byzantine agents is set as $4$, $5$, $6$, $7$, and $8$.}
\label{fig_3}
\end{figure*}

\begin{figure*}[htbp]
\centerline{\includegraphics[width=16.8cm]{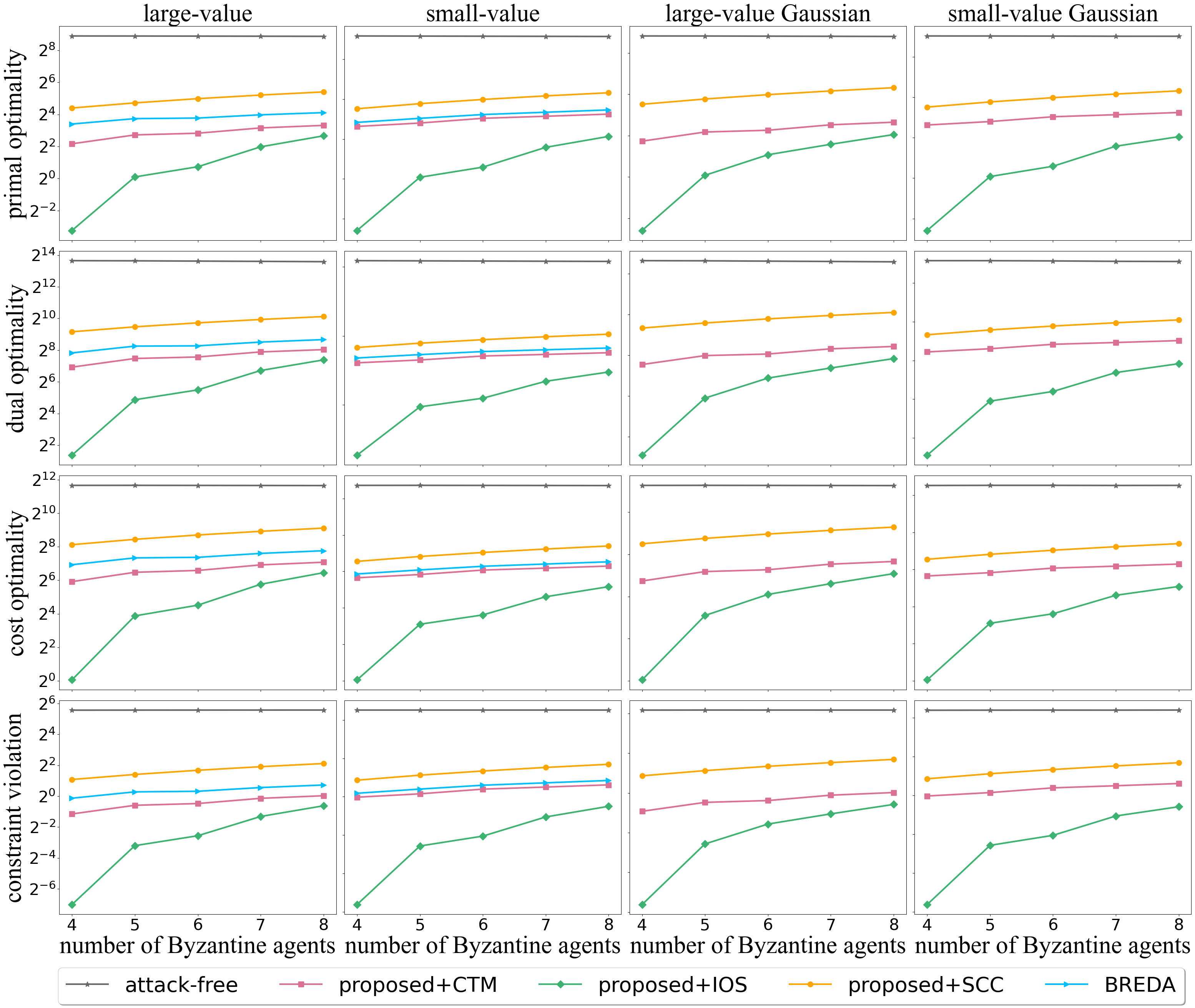}}
\caption{Primal optimality, dual optimality, cost optimality and constraint violation of the proposed algorithms with non-optimal parameters in Case 1. The number of Byzantine agents is set as $4$, $5$, $6$, $7$, and $8$.}
\label{fig_4}
\end{figure*}

\subsection{Case 1: Synthetic Problem}
We first test on a synthetic and scalar case with $D=1$. Consider a randomly generated network consisting of $J=100$ agents, where each agent has $15$ neighbors. The weight $\widetilde{e}_{ij}$ is set to $\frac{1}{16}$ if and only if $(i,j)\in \widetilde{\mathcal{E}}$ or $i=j$. The total amount of resources is $5000$ such that $\bm{s}=50$. The local constraint of each agent $i$ is $\bm{\theta}_{i} \in C_{i} = [0, 100]$. Each agent $i$ has a local cost function $f_{i}\left ( \bm{\theta}_{i} \right )=a_{i}(\bm{\theta}_{i}-b_{i})^{2}$, in which $a_i\sim \mathcal{U} (1,2)$ and $b_i\sim \mathcal{N} (2,0.6^{2})$ with $\mathcal{U}(\cdot, \cdot)$ standing for uniform distribution and $\mathcal{N} (\cdot, \cdot)$ for Gaussian distribution. Such quadratic cost functions is also used in \cite{b-LXiao-2006,b-Yun-Xu-2017,b-Jiaqi-Zhang-2020}.

We randomly select $B=6$ Byzantine agents by default, but allow each agent to have at most $4$ Byzantine neighbors. For the proposed algorithms, we test four types of Byzantine attacks: large-value, small-value, large-value Gaussian, and small-value Gaussian. With large-value attacks, a Byzantine agent sets its message as $-0.01$. With small-value attacks, a Byzantine agent sets its message as $-600$. With large-value Gaussian attacks, a Byzantine agent sets its message following a Gaussian distribution with mean $-30$ and variance $5^2$. With small-value Gaussian attacks, a Byzantine agent sets its message following a Gaussian distribution with mean $-300$ and variance $40^2$. We consider three popular robust aggregation rules: CTM, IOS and SCC. The step size is $\gamma^{k}=(k+1)^{-0.1}$, which is faster than the conservative theoretical step size in the order of $O(\frac{1}{k})$.

We use the attack-free decentralized resource allocation algorithm \eqref{eq_DEGRA-1}--\eqref{eq_DEGRA-3} as a baseline. Another baseline is BREDA. Note that BREDA defends against Byzantine attacks in the primal domain, whereas our proposed algorithms defend in the dual domain. To enable fair comparisons, for the dual-domain large-value attacks, we generate the corresponding primal-domain attacks such that their effects on the primal variables are almost the same, for our proposed algorithms and BREDA, respectively. Similarly, we also generate the corresponding primal-domain small-value attacks. Thus, with large-value and small-value attacks in BREDA, a Byzantine agent sets its message as $100$ and $0$, respectively. Note that it is difficult to generate the corresponding primal-domain large-value and small-value Gaussian attacks, and we do not compare with BREDA under these attacks.

\begin{table*}[!ht]
\hspace{-1.7cm}
\begin{minipage}{0.47\linewidth}
\renewcommand\arraystretch{1.52}
\centering
\caption{BOUNDS OF $\rho^{2}$ AND $\chi^{2}$ FOR CASE 1}
\begin{tabular}{|c|c|c|c|} \hline
 & $\rho^{2}$ & $\chi^{2}$ & $\rho^{2}+\chi^{2}$ \\ \hline
CTM & $0.44$ & $0.0031$ & $0.44$\\ \hline
IOS & $0.11$ & $0$ & $0.11$\\ \hline
SCC & $2.75$ & $0$ & $2.75$\\ \hline
\end{tabular}
\label{table-1}
\end{minipage}
\hspace{-1cm} 
\begin{minipage}{0.6\linewidth}
\centering
\renewcommand{\arraystretch}{1.18}
\caption{DUAL CONSENSUS ERRORS WITH OPTIMAL PARAMETERS FOR CASE 1}
\begin{tabular}{|c|c|c|c|c|} \hline
 & large-value & small-value & large-value Gaussian & small-value Gaussian \\ \hline
BREDA & 105.70 & 121.09 & / & /\\ \hline
proposed+CTM & 1.20e-02 & 1.07e-02 & 1.20e-02 & 1.07e-02\\ \hline
proposed+IOS & 1.09e-02 & 1.09e-02 & 1.09e-02 & 1.09e-02\\ \hline
proposed+SCC & 3.36e-02 & 3.16e-02 & 3.36e-02 & 3.16e-02\\ \hline
\end{tabular}
\label{table-2}
\end{minipage}
\end{table*}

\begin{table*}[!ht] 
\centering 
\renewcommand{\arraystretch}{1.18}
\caption{DUAL CONSENSUS ERRORS WITH NON-OPTIMAL PARAMETERS FOR CASE 1}
\begin{tabular}{|c|c|c|c|c|} \hline 
       & large-value & small-value & large-value Gaussian & small-value Gaussian\\ \hline
BREDA  &	130.48 &	144.95 & / & /\\ \hline
proposed+CTM  & 1.95e-02 & 1.58e-02 & 1.95e-02 & 1.58e-02\\ \hline
proposed+IOS  & 9.08e-02 & 9.08e-02 & 9.08e-02 & 9.08e-02\\ \hline
proposed+SCC  & 3.65e-02 & 3.70e-02 & 3.65e-02 & 3.70e-02\\ \hline
\end{tabular}\\
\label{table-3}
\end{table*}

To observe the sensitivity of CTM, IOS and SCC to their parameters, we consider two scenarios: optimal parameters and non-optimal parameters. For the scenario of optimal parameters, each honest agent sets the parameters $b$ of CTM and IOS as the number of its Byzantine neighbors. In SCC, the clipping threshold $\tau$ is determines according to Theorem 3 in \cite{b-LieHe-2022}. The results are shown in Figs. \ref{fig_1} and \ref{fig_3}. On the other hand, for the scenario of non-optimal parameters, all honest agents set $b=4$, which corresponds to the upper bound of the number of Byzantine neighbors, in CTM and IOS. In SCC, the clipping threshold is set to $\tau=0.2$. The results are shown in Figs. \ref{fig_2} and \ref{fig_4}. Performance metrics are primal optimality $\|\bm{\Theta}^{k}-\bm{\Theta}^{*}\|$, dual optimality $\sum_{i\in \mathcal{H}}\|\bm{\lambda}_{i}^{k}-\bm{\lambda}^{*}\|$, cost optimality $\|f(\bm{\Theta}^{k})-f(\bm{\Theta}^{*})\|$, constraint violation $\|\frac{1}{\mathcal{H}}\sum_{i\in \mathcal{H}}\bm{\theta}_{i}^{k}-\bm{s}\|$, and dual consensus error $\sum_{i\in \mathcal{H}}\|\bm{\lambda}_{i}^{k}-\bar{\bm{\lambda}}^{k}\|^{2}$.

Fig. \ref{fig_1} illustrates that the attack-free decentralized resource allocation algorithm \eqref{eq_DEGRA-1}--\eqref{eq_DEGRA-3} fails under all Byzantine attacks. By contrast, the proposed algorithms and BREDA demonstrate satisfactory Byzantine-resilience. Among the robust aggregation rules used in our proposed algorithms, IOS performs the best and CTM is better than SCC in terms of primal optimality, dual optimality, cost optimality, and constraint violation.
To see the reason, recall that Theorem \ref{t2} shows the primal optimality and dual optimality are both in the order of $O(\rho^{2}+\chi^{2})$. We calculate the corresponding bounds of $\rho^{2}+\chi^{2}$ in Table \ref{table-1} according to Lemmas 3--5 in \cite{b-Haoxiang-Ye-2023}. From the smallest to the largest are IOS, CTM and SCC, which validates our theoretical findings. In the numerical experiments, we observe that even though the contraction factor $\rho$ of SCC is greater than 1, the primal-dual optimality error of SCC converges to a fixed value other than explodes. However, in the theoretical analysis, we require $\rho<\frac{1-\kappa}{8\sqrt{H}}< 1$ to ensure convergence. This stricter requirement arises because we must guarantee consensus of the dual variables under Byzantine attacks. The discrepancy between our experimental and theoretical results is understandable, since the theoretical analysis must account for the worst-case scenarios and thus require more conservative conditions.

Fig. \ref{fig_1} also reveals that BREDA is worse than the proposed algorithms with proper robust aggregation rules. To further highlight the advantages of our proposed algorithms, we list the dual consensus errors in Table \ref{table-2}. No matter the types of Byzantine attacks and robust aggregation rules, the proposed algorithms are all able to achieve nearly perfect dual consensus. By contrast, BREDA cannot guarantee dual consensus. This phenomenon reveals the benefits of the dual-domain defenses.

When the parameters are non-optimal, the above conclusions still hold. Fig. \ref{fig_2} and Table \ref{table-3} demonstrate the advantages of our proposed dual-domain defense algorithms over BREDA in the scenario with non-optimal parameters.

In Figs. \ref{fig_3} and \ref{fig_4}, we check the sensitivity of the compared algorithms to the number of Byzantine agents $B$ by setting $B$ as $4$, $5$, $6$, $7$ and $8$. The attack-free decentralized resource allocation algorithm \eqref{eq_DEGRA-1}--\eqref{eq_DEGRA-3} fails for any $B$. By contrast, both BREDA and our proposed algorithms demonstrate satisfactory resilience, and their performance is steady when $B$ varies.

\begin{figure*}[htbp]
\centerline{\includegraphics[width=16.8cm]{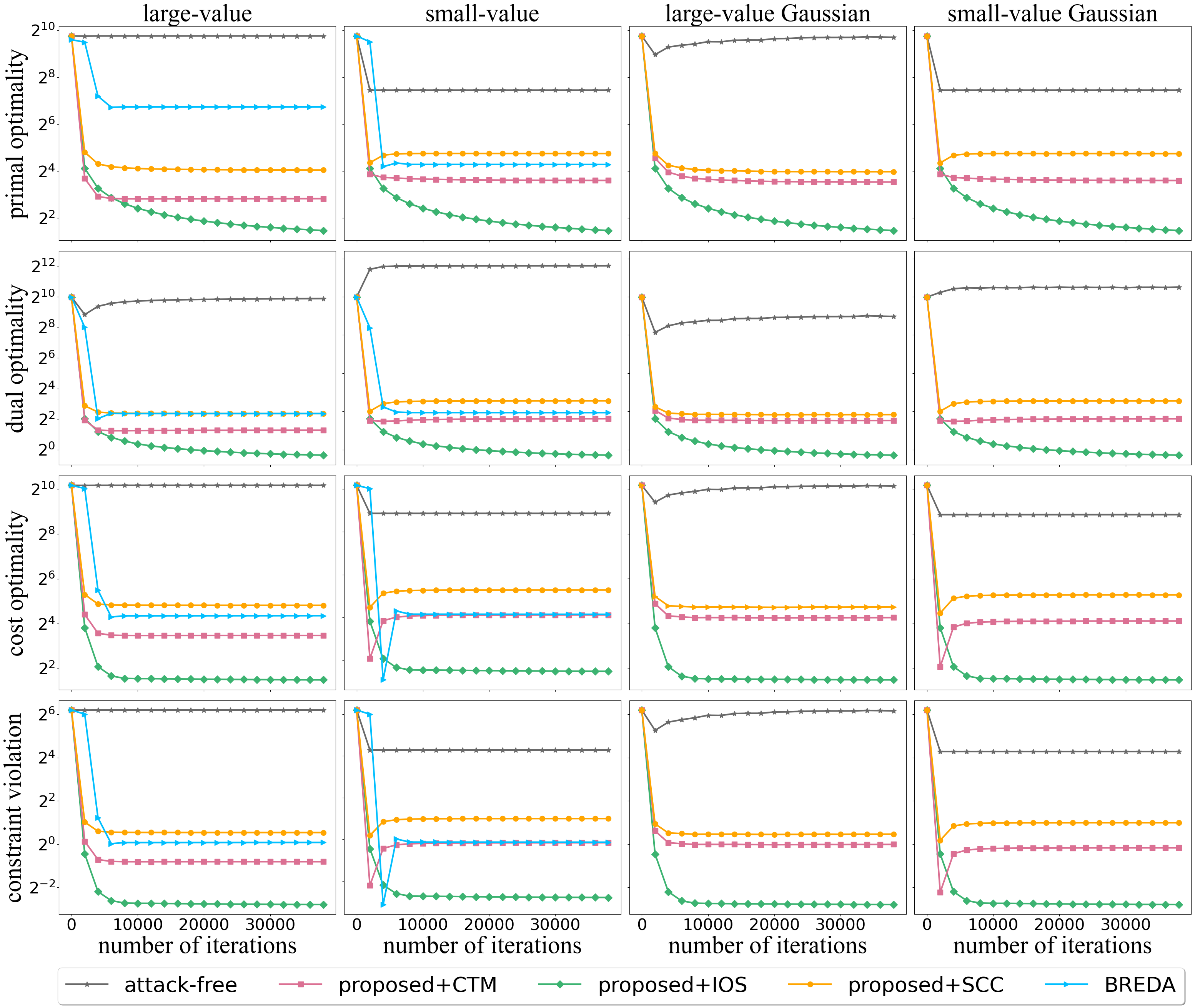}}
\caption{Primal optimality, dual optimality, cost optimality and constraint violation of the compared algorithms with optimal parameters in Case 2.}
\label{fig_5}
\end{figure*}

\begin{table*}[!ht]
\hspace{-1.7cm}
\begin{minipage}{0.47\linewidth}
\renewcommand\arraystretch{1.52}
\centering
\caption{BOUNDS OF $\rho^{2}$ AND $\chi^{2}$ FOR CASE 2 }
\begin{tabular}{|c|c|c|c|r|l|} \hline 
     & $\rho^{2}$ & $\chi^{2}$ & $\rho^{2}+\chi^{2}$ \\ \hline
CTM  & $0.024$ & $0.11$ & $0.134$\\ \hline
IOS  & $ 0.006$ & $0$ & $0.006$\\ \hline
SCC  & $0.965$ & $0$ & $0.965$\\ \hline
\end{tabular}
\label{table-4}
\end{minipage}
\hspace{-1cm} 
\begin{minipage}{0.6\linewidth}
\centering
\renewcommand{\arraystretch}{1.18}
\caption{DUAL CONSENSUS ERRORS WITH OPTIMAL PARAMETERS FOR CASE 2}
\begin{tabular}{|c|c|c|c|c|} \hline 
       & large-value & small-value & large-value Gaussian & small-value Gaussian\\ \hline
BREDA  & 	0.51 &	0.49 & / & /\\ \hline
proposed+CTM  & 2.16e-04 & 3.28e-03 & 3.48e-03 & 3.28e-03\\ \hline
proposed+IOS  & 3.37e-03 & 3.37e-03 & 3.37e-03 & 3.37e-03\\ \hline
proposed+SCC  & 3.55e-03 & 3.23e-03 & 3.54e-03 & 3.23e-03\\ \hline
\end{tabular}
\label{table-5}
\end{minipage}
\end{table*}

\subsection{Case 2: Economic Dispatch for IEEE 118-Bus Test System}
We next consider a power dispatch problem for the IEEE 118-bus test system, which contains 54 generators \cite{b-IEEE-118}. Each generator $i$ has a local power $\bm{\theta}_{i}$ and a corresponding cost function $f_{i}(\bm{\theta}_{i})=\eta_{i}\bm{\theta}_{i}^{2}+\zeta_{i}\bm{\theta}_{i}+\xi_{i}$, where $\eta_{i}\in [0.0024, 0.0697]$, $\zeta_{i}\in [8.3391, 37.6968]$, and $\xi_{i}\in [6.78, 74.33]$. The local constraint of each agent $i$ is $\bm{\theta}_{i}\in [\bm{\theta}_{i}^{\min},\bm{\theta}_{i}^{\max}]$, where $\bm{\theta}_{i}^{\min}\in [5,150]$ and $\bm{\theta}_{i}^{\max}\in [30,420]$. The total amount of resources is set as 6000, such that $\bm{s}=\frac{6000}{54}$ \cite{b-Thinh-T.Doan-2021}. To test the performance of the proposed algorithms, we randomly select one Byzantine agent out of the 54 generators and apply different types of Byzantine attacks, including large-value, small-value, large-value Gaussian, and small-value Gaussian. For large-value attacks, the Byzantine generator sets its message as $-0.01$, whereas for small-value attacks, the Byzantine generator sets its message as $-100$. For large-value Gaussian attacks, the Byzantine generator sets its message following a Gaussian distribution with mean $-10$ and variance $5^{2}$. For small-value Gaussian attacks, the Byzantine generator sets its message following a Gaussian distribution with mean $-50$ and variance $10^{2}$. We also design the corresponding larger-value and smaller-value attacks for BREDA, where the Byzantine generator sets its message as $420$ and $5$, respectively. The weight matrix $\widetilde{E}$ is constructed according to the Metropolis constant weight rule \cite{b-Wei-Shi-2015}. The parameters $b$ and $\tau$ are optimal. The step size for the proposed algorithms is determined as $\gamma^{k}=(k+1)^{-0.7}$.

Fig. \ref{fig_5} demonstrates the failure of the attack-free decentralized resource allocation algorithm, as well as the resilience of the proposed algorithms and BREDA against various Byzantine attacks. We also calculate the corresponding bounds of $\rho^{2}+\chi^{2}$ of the robust aggregation rules IOS, CTM, and SCC, as presented in Table \ref{table-4}.
Observe that a smaller bound of $\rho^{2}+\chi^{2}$ leads to better performance, which has been predicted by our theoretical findings.

According to Fig. \ref{fig_5}, BREDA performs worse than the proposed algorithms with proper robust aggregation rules. We calculate the dual consensus errors of the proposed algorithms with different robust aggregation rules and BREDA, as presented in Table \ref{table-5}. The proposed algorithms achieve nearly consensual dual variables and BREDA does not.

\vspace{-0.1cm}
\section{CONCLUSIONS AND FUTURE WORK}\label{sec 7}
In this paper, we address the challenging Byzantine-resili- ence issue in decentralized resource allocation. We propose a class of Byzantine-resilient algorithms that leverage robust aggregation rules within a dual-domain defense framework. Given that the robust aggregation rules are properly designed, we prove that the primal and dual variables of the honest agents converge to the neighborhoods of their optima, while the dual variables are able to reach consensus. This dual-domain defense approach not only simplifies the algorithmic updates but also enhances the overall Byzantine-resilience. Our numerical experiments further demonstrate the resilience of the proposed algorithms against various Byzantine attacks, confirming their practical utility.

In the future, we plan to extend our algorithm development and theoretical analysis to stochastic and online decentralized resource allocation problems under Byzantine attacks, which are of particular importance for time-sensitive applications.

\begin{appendices}
\section{Proof of Theorem \ref{t2}}\label{appendix2}
\subsection{Part $a$ of Theorem \ref{t2}}
According to the update of $\bm{\lambda}_{i}^{k+1}$ in Algorithm \ref{alg2}, we have
\begin{align}
\label{eq_proof_theorem2-b-1}
&\|\bar{\bm{\lambda}}^{k+1}-\bm{\lambda}^{*}\|^{2}\\
=&\|\frac{1}{H}\sum_{i\in \mathcal{H}}
AGG_{i}(\bm{\lambda}_{i}^{k+\frac{1}{2}},\{\check{\bm{\lambda}}_{j}^{k+\frac{1}{2}}\}_{j\in \mathcal{N}_{i}})-\bm{\lambda}^{*}\|^{2}\notag\\
=&\|\frac{1}{H}\sum_{i\in \mathcal{H}}
AGG_{i}(\bm{\lambda}_{i}^{k+\frac{1}{2}},\{\check{\bm{\lambda}}_{j}^{k+\frac{1}{2}}\}_{j\in \mathcal{N}_{i}})-\bar{\bm{\lambda}}^{k+\frac{1}{2}}+\bar{\bm{\lambda}}^{k+\frac{1}{2}}-\bm{\lambda}^{*}\|^{2}\notag \\
\le& \frac{1}{v_{1}}\|\frac{1}{H}\sum_{i\in \mathcal{H}}AGG_{i}(\bm{\lambda}_{i}^{k+\frac{1}{2}},\{\check{\bm{\lambda}}_{j}^{k+\frac{1}{2}}\}_{j\in \mathcal{N}_{i}})-\bar{\bm{\lambda}}^{k+\frac{1}{2}}\|^{2}\notag\\
&+\frac{1}{1-v_{1}}\|\bar{\bm{\lambda}}^{k+\frac{1}{2}}-\bm{\lambda}^{*}\|^{2} \notag\\
=&\frac{1}{v_{1}}\|\frac{1}{H}\sum_{i\in \mathcal{H}}AGG_{i}(\bm{\lambda}_{i}^{k+\frac{1}{2}},\{\check{\bm{\lambda}}_{j}^{k+\frac{1}{2}}\}_{j\in \mathcal{N}_{i}})-\frac{1}{H}\sum_{i\in\mathcal{H}}\bar{\bm{\lambda}}_{i}^{k+\frac{1}{2}}\notag\\
&+\frac{1}{H}\sum_{i\in\mathcal{H}}\bar{\bm{\lambda}}_{i}^{k+\frac{1}{2}}-\bar{\bm{\lambda}}^{k+\frac{1}{2}}\|^{2}+\frac{1}{1-v_{1}}\|\bar{\bm{\lambda}}^{k+\frac{1}{2}}-\bm{\lambda}^{*}\|^{2}\notag\\
\le&\frac{2}{v_{1}}\|\frac{1}{H}\sum_{i\in \mathcal{H}}AGG_{i}(\bm{\lambda}_{i}^{k+\frac{1}{2}},\check{\bm{\lambda}}_{j}^{k+\frac{1}{2}}(j\in \mathcal{N}_{i}))-\frac{1}{H}\sum_{i\in\mathcal{H}}\bar{\bm{\lambda}}_{i}^{k+\frac{1}{2}}\|^{2}\notag
\end{align}
\begin{align}
&+\frac{2}{v_{1}}\|\frac{1}{H}\sum_{i\in\mathcal{H}}\bar{\bm{\lambda}}_{i}^{k+\frac{1}{2}}-\bar{\bm{\lambda}}^{k+\frac{1}{2}}\|^{2}+\frac{1}{1-v_{1}}\|\bar{\bm{\lambda}}^{k+\frac{1}{2}}-\bm{\lambda}^{*}\|^{2}\notag\\
\le & \underbrace{\frac{2}{v_{1}H}\sum_{i\in \mathcal{H}}\|AGG_{i}(\bm{\lambda}_{i}^{k+\frac{1}{2}},\{\check{\bm{\lambda}}_{j}^{k+\frac{1}{2}}\}_{j\in \mathcal{N}_{i}})-\bar{\bm{\lambda}}_{i}^{k+\frac{1}{2}}\|^{2}}_{T_{1}}\notag\\
&+\underbrace{\frac{2}{v_{1}}\|\frac{1}{H}\sum_{i\in\mathcal{H}}\bar{\bm{\lambda}}_{i}^{k+\frac{1}{2}}-\bar{\bm{\lambda}}^{k+\frac{1}{2}}\|^{2}}_{T_{2}}+\underbrace{\frac{1}{1-v_{1}}\|\bar{\bm{\lambda}}^{k+\frac{1}{2}}-\bm{\lambda}^{*}\|^{2}}_{T_{3}},\notag
\end{align}
where $v_{1}$ is any positive constant in $(0,1)$. To derive the first inequality, we use $\|\bm{a}+\bm{b}\|^{2}\le \frac{1}{v}\|\bm{a}\|^{2}+\frac{1}{1-v}\|\bm{b}\|^{2}$ for any positive constant $v\in (0,1)$. The last inequality holds because $(a_{1}+\cdots+a_{H})^{2}\le H(a_{1}^{2}+\cdots+a_{H}^{2})$. Next, we analyze $T_{1}$, $T_{2}$ and $T_{3}$ in turn.

\noindent\textbf{Bounding $T_{1}$:}
According to \eqref{eq_definition1} in Definition \ref{d1}, $T_{1}$ can be bounded by
\begin{align}
\label{eq_proof_theorem2-b-T1-1}
T_{1}\le& \frac{2}{v_{1}H}\sum_{i\in \mathcal{H}}\rho^{2}\max_{j\in \mathcal{N}_{i}\cap\mathcal{H}\cup \{i\}}\|\bm{\lambda}_{j}^{k+\frac{1}{2}}-\bar{\bm{\lambda}}_{i}^{k+\frac{1}{2}}\|^{2}
\end{align}
\begin{align}
=&\frac{2\rho^{2}}{v_{1}H}\sum_{i\in \mathcal{H}}\max_{j\in \mathcal{N}_{i}\cap\mathcal{H}\cup \{i\}}\|\bm{\lambda}_{j}^{k+\frac{1}{2}}-\bar{\bm{\lambda}}^{k+\frac{1}{2}}+\bar{\bm{\lambda}}^{k+\frac{1}{2}}-\bar{\bm{\lambda}}_{i}^{k+\frac{1}{2}}\|^{2} \notag \\
\le & \frac{4\rho^{2}}{v_{1}H}\sum_{i\in \mathcal{H}}\max_{j\in \mathcal{N}_{i}\cap\mathcal{H}\cup \{i\}}\|\bm{\lambda}_{j}^{k+\frac{1}{2}}-\bar{\bm{\lambda}}^{k+\frac{1}{2}}\|^{2}\notag\\
&+ \frac{4\rho^{2}}{v_{1}H}\sum_{i\in \mathcal{H}}\|\bar{\bm{\lambda}}_{i}^{k+\frac{1}{2}}-\bar{\bm{\lambda}}^{k+\frac{1}{2}}\|^{2}\notag \\
\le & \frac{4\rho^{2}}{v_{1}H}\sum_{i\in \mathcal{H}}\max_{i\in \mathcal{H}}\|\bm{\lambda}_{i}^{k+\frac{1}{2}}-\bar{\bm{\lambda}}^{k+\frac{1}{2}}\|^{2}\notag\\
&+ \frac{4\rho^{2}}{v_{1}H}\sum_{i\in \mathcal{H}}\max_{i\in \mathcal{H}}\|\bm{\lambda}_{i}^{k+\frac{1}{2}}-\bar{\bm{\lambda}}^{k+\frac{1}{2}}\|^{2}\notag\\
=& \frac{8\rho^{2}}{v_{1}}\max_{i\in \mathcal{H}}\|\bm{\lambda}_{i}^{k+\frac{1}{2}}-\bar{\bm{\lambda}}^{k+\frac{1}{2}}\|^{2}.\notag
\end{align}

Define $\Lambda=[\cdots,\bm{\lambda}_{i},\cdots]\in \mathbb{R}^{H\times D} $ that collects $\bm{\lambda}_{i}$ of all honest agents $i\in \mathcal{H}$. Combining the fact $\max_{i\in \mathcal{H}}\|\bm{\lambda}_{i}^{k+\frac{1}{2}}-\bar{\bm{\lambda}}^{k+\frac{1}{2}}\|^{2}\le \|\Lambda^{k+\frac{1}{2}}-\frac{1}{H}\bm{1}\bm{1}^{\top}\Lambda^{k+\frac{1}{2}}
\|^{2}_{F}$ and \eqref{eq_proof_theorem2-b-T1-1}, we obatin
\begin{align}
\label{eq_proof_theorem2-b-T1}
T_{1}\le \frac{8\rho^{2}}{v_{1}}\|\Lambda^{k+\frac{1}{2}}-\frac{1}{H}\bm{1}\bm{1}^{\top}\Lambda^{k+\frac{1}{2}}
\|^{2}_{F}.
\end{align}

\noindent\textbf{Bounding $T_{2}$:} By the definition of $\bar{\bm{\lambda}}_{i}=\sum_{j\in \mathcal{N}_{i}\cap\mathcal{H}\cup \{i\}}e_{i,j}\bm{\lambda}_{j}$, we have
\begin{align}
\label{eq_proof_theorem2-b-T2-1}
T_{2}=& \frac{2}{v_{1}}\|\frac{1}{H}\sum_{i\in\mathcal{H}}\sum_{j\in \mathcal{N}_{i}\cap\mathcal{H}\cup \{i\}}e_{i,j}\bm{\lambda}_{j}^{k+\frac{1}{2}}-\bar{\bm{\lambda}}^{k+\frac{1}{2}}\|^{2}\\
=& \frac{2}{v_{1}}\|\frac{1}{H}\bm{1}^{\top}E\Lambda^{k+\frac{1}{2}}-\frac{1}{H}\bm{1}^{\top}\Lambda^{k+\frac{1}{2}}
\|^{2}\notag\\
=& \frac{2}{v_{1}}\|\frac{1}{H}\bm{1}^{\top}(E\Lambda^{k+\frac{1}{2}}-\frac{1}{H}\bm{1}\bm{1}^{\top}\Lambda^{k+\frac{1}{2}})
\|^{2} \notag\\
=& \frac{2}{v_{1}H^{2}} \|\bm{1}^{\top}(E\Lambda^{k+\frac{1}{2}}-\frac{1}{H}\bm{1}\bm{1}^{\top}\Lambda^{k+\frac{1}{2}})
\|^{2}  \notag\\
=& \frac{2}{v_{1}H^{2}} \|(\bm{1}^{\top}E-\bm{1}^{\top})(\Lambda^{k+\frac{1}{2}}-\frac{1}{H}\bm{1}\bm{1}^{\top}\Lambda^{k+\frac{1}{2}})
\|^{2} \notag\\
\le&\frac{2}{v_{1}H^{2}} \|E^{\top}\bm{1}-\bm{1}\|^{2}\|\Lambda^{k+\frac{1}{2}}-\frac{1}{H}\bm{1}\bm{1}^{\top}\Lambda^{k+\frac{1}{2}}
\|^{2}. \notag
\end{align}
To drive the last equality, we use Definition \ref{d1} that the virtual weight matrix $E$ is row stochastic.

Define $\chi^{2}=\frac{1}{H}\|E^{\top}\bm{1}-\bm{1}\|^{2}$ to quantify how non-column stochastic the virtual weight matrix $E$ is. Applying the fact $\|\cdot\|^{2}\le \|\cdot\|_{F}^{2}$ to the right-hand side of \eqref{eq_proof_theorem2-b-T2-1}, we have
\begin{align}
\label{eq_proof_theorem2-b-T2}
T_{2}\le \frac{2\chi^{2}}{v_{1}H}\|\Lambda^{k+\frac{1}{2}}-\frac{1}{H}\bm{1}\bm{1}^{\top}\Lambda^{k+\frac{1}{2}}
\|_{F}^{2}.
\end{align}

\noindent\textbf{Bounding $T_{3}$:}
Averaging both sides of \eqref{eq_BRE-DEGRA-2} over $i\in \mathcal{H}$, we have
\begin{align}
\label{eq_proof_theorem2-b-T3-1}
\bar{\bm{\lambda}}^{k+\frac{1}{2}}=\bar{\bm{\lambda}}^{k}-\frac{\gamma^{k}}{H}\sum_{i\in \mathcal{H}}(\frac{1}{J}\bm{s}-\frac{1}{J}\bm{\theta}_{i}^{k}).
\end{align}
The dual problem of \eqref{eq_oracle-problem} can be written as a minimization problem in the form of
\begin{align}
\label{eq_dual-problem_oracle}
\underset{\bm{\lambda }\in \mathbb{R}^D}{\min} ~ g(\bm{\lambda })=\sum_{i\in\mathcal{H}}g_{i}(\bm{\lambda }),
\end{align}
where $g_{i}(\bm{\lambda }):=\frac{1}{H}F_{i}^{*}(-\bm{\lambda })+\frac{1}{H}\bm{\lambda }^{\top}\bm{s}$ and $F_{i}^{*}(\bm{\lambda}):=\max_{\bm{\theta}_{i}\in C_{i}} \{\bm{\lambda}^{\top}\bm{\theta}_{i}-f_{i}(\bm{\theta}_{i})\}$.
Based on the definition of $g_{i}(\bm{\lambda }):=\frac{1}{H}F_{i}^{*}(-\bm{\lambda })+\frac{1}{H}\bm{\lambda }^{\top}\bm{s}$ and Danskin's theorem \cite{b-Dimitri-P.-Bertsekas-1999}, we have
\begin{align}
\label{eq_gradient_dual_function_oracle}
\nabla g_{i}(\bm{\lambda}_{i})=\frac{1}{H}\bm{s}-\frac{1}{H}\arg\underset{\bm{\theta}_{i}\in C_{i}}{\min}\{\bm{\lambda}_{i}^{\top}\bm{\theta}_{i}+f_{i}(\bm{\theta}_{i})\}.
\end{align}

Combining \eqref{eq_BRE-DEGRA-1}, \eqref{eq_proof_theorem2-b-T3-1} and \eqref{eq_gradient_dual_function_oracle}, we can obtain
\begin{align}
\label{eq_proof_theorem2-b-T3-2}
\bar{\bm{\lambda}}^{k+\frac{1}{2}}=\bar{\bm{\lambda}}^{k}-\frac{\gamma^{k}}{J}\sum_{i\in \mathcal{H}}\nabla g_{i}(\bm{\lambda}_{i}^{k}).
\end{align}

Substituting \eqref{eq_proof_theorem2-b-T3-2} into $T_{3}$, we have
\begin{align}
\label{eq_proof_theorem2-b-T3-3}
T_{3}=&\frac{1}{1-v_{1}}\|\bar{\bm{\lambda}}^{k}-\bm{\lambda}^{*}-\frac{\gamma^{k}}{J}\sum_{i\in \mathcal{H}}\nabla g_{i}(\bar{\bm{\lambda}}^{k})\\
&+\frac{\gamma^{k}}{J}\sum_{i\in \mathcal{H}}\nabla g_{i}(\bar{\bm{\lambda}}^{k})-\frac{\gamma^{k}}{J}\sum_{i\in \mathcal{H}}\nabla g_{i}(\bm{\lambda}_{i}^{k})\|^{2}\notag\\
=&\frac{1}{1-v_{1}}\|\bar{\bm{\lambda}}^{k}-\bm{\lambda}^{*}-\frac{\gamma^{k}}{J}\sum_{i\in \mathcal{H}}\nabla g_{i}(\bar{\bm{\lambda}}^{k})\|^{2}\notag\\
&+\frac{(\gamma^{k})^{2}}{1-v_{1}}\|\frac{1}{J}\sum_{i\in \mathcal{H}}(\nabla g_{i}(\bar{\bm{\lambda}}^{k})-\nabla g_{i}(\bm{\lambda}_{i}^{k}))\|^{2}+\frac{2\gamma^{k}}{1-v_{1}}\cdot  \notag\\
&\left \langle \bar{\bm{\lambda}}^{k}-\bm{\lambda}^{*}-\frac{\gamma^{k}}{J}\sum_{i\in \mathcal{H}}\nabla g_{i}(\bar{\bm{\lambda}}^{k}), \frac{1}{J}\sum_{i\in \mathcal{H}}(\nabla g_{i}(\bar{\bm{\lambda}}^{k})-\nabla g_{i}(\bm{\lambda}_{i}^{k}))\right \rangle \notag\\
\le& \frac{1}{1-v_{1}}\|\bar{\bm{\lambda}}^{k}-\bm{\lambda}^{*}-\frac{\gamma^{k}}{J}\sum_{i\in \mathcal{H}}\nabla g_{i}(\bar{\bm{\lambda}}^{k})\|^{2}\notag\\
&+\frac{(\gamma^{k})^{2}}{1-v_{1}}\|\frac{1}{J}\sum_{i\in \mathcal{H}}(\nabla g_{i}(\bar{\bm{\lambda}}^{k})-\nabla g_{i}(\bm{\lambda}_{i}^{k}))\|^{2}\notag\\
&+ \frac{v_{2}^{-1}\gamma^{k}}{1-v_{1}}\|\frac{1}{J}\sum_{i\in \mathcal{H}}(\nabla g_{i}(\bar{\bm{\lambda}}^{k})-\nabla g_{i}(\bm{\lambda}_{i}^{k}))\|^{2}\notag\\
&+\frac{v_{2}\gamma^{k}}{1-v_{1}}\|\bar{\bm{\lambda}}^{k}-\bm{\lambda}^{*}-\frac{\gamma^{k}}{J}\sum_{i\in \mathcal{H}}\nabla g_{i}(\bar{\bm{\lambda}}^{k})\|^{2}\notag\\
\le&\frac{1+v_{2}\gamma^{k}}{1-v_{1}}\|\bar{\bm{\lambda}}^{k}-\bm{\lambda}^{*}-\frac{\gamma^{k}}{J}\sum_{i\in \mathcal{H}}\nabla g_{i}(\bar{\bm{\lambda}}^{k})\|^{2}\notag\\
&+\frac{\gamma^{k}(\gamma^{k}+v_{2}^{-1})H}{(1-v_{1})J^{2}}\sum_{i\in \mathcal{H}}\|\nabla g_{i}(\bar{\bm{\lambda}}^{k})-\nabla g_{i}(\bm{\lambda}_{i}^{k})\|^{2},\notag
\end{align}
where $v_{2}>0$ is any positive constant. To drive the first inequality, we use $2\bm{a}^{\top}\bm{b}\le v^{-1}\|\bm{a}\|^{2}+v\|\bm{b}\|^{2}$ for any $v>0$. The last inequality holds because $(a_{1}+\cdots+a_{H})^{2}\le H(a_{1}^{2}+\cdots+a_{H}^{2})$. Next, we analyze the first and second terms at the right-hand side of \eqref{eq_proof_theorem2-b-T3-3} in turn.

According to the fact that $\sum_{i\in \mathcal{H}}\nabla g_{i}(\bm{\lambda}^{*})=\bm{0}$, we have
\begin{align}
\label{eq_proof_theorem2-b-T3-4}
&\frac{1+v_{2}\gamma^{k}}{1-v_{1}}\|\bar{\bm{\lambda}}^{k}-\bm{\lambda}^{*}-\frac{\gamma^{k}}{J}\sum_{i\in \mathcal{H}}\nabla g_{i}(\bar{\bm{\lambda}}^{k})\|^{2}\\
=& \frac{1+v_{2}\gamma^{k}}{1-v_{1}}\|\bar{\bm{\lambda}}^{k}-\bm{\lambda}^{*}-\frac{\gamma^{k}}{J}\sum_{i\in \mathcal{H}}\nabla g_{i}(\bar{\bm{\lambda}}^{k})-\frac{\gamma^{k}}{J}\sum_{i\in \mathcal{H}}\nabla g_{i}(\bm{\lambda}^{*})\|^{2}\notag 
\end{align}
\begin{align}
=&\frac{(1+v_{2}\gamma^{k})(\gamma^{k})^{2}H^{2}}{(1-v_{1})J^{2}}\|\frac{1}{H}\sum_{i\in \mathcal{H}}\nabla g_{i}(\bar{\bm{\lambda}}^{k})-\frac{1}{H}\sum_{i\in \mathcal{H}}\nabla g_{i}(\bm{\lambda}^{*})\|^{2}\notag\\
&+\frac{1+v_{2}\gamma^{k}}{1-v_{1}}\|\bar{\bm{\lambda}}^{k}-\bm{\lambda}^{*}\|^{2}- \frac{2\gamma^{k}(1+v_{2}\gamma^{k})H}{(1-v_{1})J}\cdot \notag\\
&\left \langle  \bar{\bm{\lambda}}^{k}-\bm{\lambda}^{*},\frac{1}{H}\sum_{i\in \mathcal{H}}\nabla g_{i}(\bar{\bm{\lambda}}^{k})-\frac{1}{H}\sum_{i\in \mathcal{H}}\nabla g_{i}(\bm{\lambda}^{*}) \right \rangle  .\notag
\end{align}
For $\left \langle  \bar{\bm{\lambda}}^{k}-\bm{\lambda}^{*},\frac{1}{H}\sum_{i\in \mathcal{H}}\nabla g_{i}(\bar{\bm{\lambda}}^{k})-\frac{1}{H}\sum_{i\in \mathcal{H}}\nabla g_{i}(\bm{\lambda}^{*}) \right \rangle $,  the last term at the right-hand side of \eqref{eq_proof_theorem2-b-T3-4}, we have the following bound. According to Lemma \ref{l2}, $g_{i}(\cdot)$ is $\frac{1}{HL_{f}}$-strongly convex and $\frac{1}{Hu_{f}}$-smooth. By Lemma 3 in \cite{b-Kun-Yuan-2016}, since $\frac{1}{H}\sum_{i\in \mathcal{H}}g_{i}(\cdot)$ is $\frac{1}{HL_{f}}$-strongly convex and $\frac{1}{Hu_{f}}$-smooth, we have

\begin{align}
\label{eq_proof_theorem2-b-T3-5}
&\left \langle  \bar{\bm{\lambda}}^{k}-\bm{\lambda}^{*},\frac{1}{H}\sum_{i\in \mathcal{H}}\nabla g_{i}(\bar{\bm{\lambda}}^{k})-\frac{1}{H}\sum_{i\in \mathcal{H}}\nabla g_{i}(\bm{\lambda}^{*}) \right \rangle\\
\ge & \alpha \|\frac{1}{H}\sum_{i\in \mathcal{H}}\nabla g_{i}(\bar{\bm{\lambda}}^{k})-\frac{1}{H}\sum_{i\in \mathcal{H}}\nabla g_{i}(\bm{\lambda}^{*})\|^{2}+\beta\|\bar{\bm{\lambda}}^{k}-\bm{\lambda}^{*}\|^{2},\notag
\end{align}
where $\alpha=\frac{Hu_{f}L_{f}}{u_{f}+L_{f}}$ and $\beta=\frac{1}{H(u_{f}+L_{f})}$. Substituting \eqref{eq_proof_theorem2-b-T3-5} into \eqref{eq_proof_theorem2-b-T3-4} and rearranging the terms, we have
\begin{align}
\label{eq_proof_theorem2-b-T3-6}
&\frac{1+v_{2}\gamma^{k}}{1-v_{1}}\|\bar{\bm{\lambda}}^{k}-\bm{\lambda}^{*}-\frac{\gamma^{k}}{J}\sum_{i\in \mathcal{H}}\nabla g_{i}(\bar{\bm{\lambda}}^{k})\|^{2}\\
\le& \frac{(1+v_{2}\gamma^{k})(1-2\gamma^{k}\beta\cdot\frac{H}{J})}{1-v_{1}}\|\bar{\bm{\lambda}}^{k}-\bm{\lambda}^{*}\|^{2}\notag\\
&+\frac{(1+v_{2}\gamma^{k})((\gamma^{k})^{2}\cdot\frac{H}{J}-2\gamma^{k}\alpha)}{1-v_{1}}\cdot\notag\\
&\frac{H}{J}\|\frac{1}{H}\sum_{i\in \mathcal{H}}\nabla g_{i}(\bar{\bm{\lambda}}^{k})-\frac{1}{H}\sum_{i\in \mathcal{H}}\nabla g_{i}(\bm{\lambda}^{*})\|^{2} \notag\\
\le &\frac{(1+v_{2}\gamma^{k})(1-\gamma^{k}\beta\cdot\frac{H}{J})}{1-v_{1}}\|\bar{\bm{\lambda}}^{k}-\bm{\lambda}^{*}\|^{2},\notag
\end{align}
where the last inequality holds with a proper step size $\gamma^{k}$ satisfying $(\gamma^{k})^{2}\cdot\frac{H}{J}-2\gamma^{k}\alpha\le 0$.

Since $g_{i}(\cdot)$ is $\frac{1}{Hu_{f}}$-smooth, we obtain
\begin{align}
\label{eq_proof_theorem2-b-T3-7}
&\frac{\gamma^{k}(\gamma^{k}+v_{2}^{-1})}{(1-v_{1})H}\cdot\frac{H^{2}}{J^{2}}\sum_{i\in \mathcal{H}}\|\nabla g_{i}(\bar{\bm{\lambda}}^{k})-\nabla g_{i}(\bm{\lambda}_{i}^{k})\|^{2}\\
\le& \frac{\gamma^{k}(\gamma^{k}+v_{2}^{-1})}{(1-v_{1})H^{3}u_{f}^{2}}\cdot\frac{H^{2}}{J^{2}}\sum_{i\in \mathcal{H}}\|\bm{\lambda}_{i}^{k}-\bar{\bm{\lambda}}^{k}\|^{2}\notag\\
= &\frac{\gamma^{k}(\gamma^{k}+v_{2}^{-1})}{(1-v_{1})H^{3}u_{f}^{2}}\cdot\frac{H^{2}}{J^{2}}\|\Lambda^{k}-\frac{1}{H}\bm{1}\bm{1}^{\top}\Lambda^{k}\|_{F}^{2}.\notag
\end{align}
Substituting \eqref{eq_proof_theorem2-b-T3-6} and \eqref{eq_proof_theorem2-b-T3-7} into \eqref{eq_proof_theorem2-b-T3-3} and rearranging the terms, we obtain
\begin{align}
\label{eq_proof_theorem2-b-T3}
T_{3}\le & \frac{(1+v_{2}\gamma^{k})(1-\gamma^{k}\beta\cdot\frac{H}{J})}{1-v_{1}}\|\bar{\bm{\lambda}}^{k}-\bm{\lambda}^{*}\|^{2} \\
&+\frac{\gamma^{k}(\gamma^{k}+v_{2}^{-1})}{(1-v_{1})H^{3}u_{f}^{2}}\cdot\frac{H^{2}}{J^{2}}\|\Lambda^{k}-\frac{1}{H}\bm{1}\bm{1}^{\top}\Lambda^{k}\|_{F}^{2}.\notag
\end{align}

Substituting \eqref{eq_proof_theorem2-b-T1}, \eqref{eq_proof_theorem2-b-T2} and \eqref{eq_proof_theorem2-b-T3} into \eqref{eq_proof_theorem2-b-1} and rearranging the terms, we have
\begin{align}
\label{eq_proof_theorem2-b-2}
&\|\bar{\bm{\lambda}}^{k+1}-\bm{\lambda}^{*}\|^{2}\\
\le & \frac{(1+v_{2}\gamma^{k})(1-\gamma^{k}\beta\cdot\frac{H}{J})}{1-v_{1}}\|\bar{\bm{\lambda}}^{k}-\bm{\lambda}^{*}\|^{2} \notag\\
+ & \frac{\gamma^{k}(\gamma^{k}+v_{2}^{-1})}{(1-v_{1})H^{3}u_{f}^{2}}\cdot\frac{H^{2}}{J^{2}}\|\Lambda^{k}-\frac{1}{H}\bm{1}\bm{1}^{\top}\Lambda^{k}\|_{F}^{2}\notag\\
+ & \frac{8\rho^{2}}{v_{1}}\|\Lambda^{k+\frac{1}{2}}-\frac{1}{H}\bm{1}\bm{1}^{\top}\Lambda^{k+\frac{1}{2}}
\|^{2}_{F}+\frac{2\chi^{2}}{v_{1}H}\|\Lambda^{k+\frac{1}{2}}-\frac{1}{H}\bm{1}\bm{1}^{\top}\Lambda^{k+\frac{1}{2}}
\|_{F}^{2}.\notag
\end{align}
Setting $v_{1}=\frac{\gamma^{k}\beta\cdot\frac{H}{J}}{4}\in(0,1)$ and $v_{2}=\frac{\beta\cdot\frac{H}{J}}{2(1-\gamma^{k}\beta\cdot\frac{H}{J})}>0$, from \eqref{eq_proof_theorem2-b-2} we have
\begin{align}
\label{eq_proof_theorem2-b-3}
&\|\bar{\bm{\lambda}}^{k+1}-\bm{\lambda}^{*}\|^{2}\\
\le & (1-\frac{\gamma^{k}\beta\cdot\frac{H}{J}}{4})\|\bar{\bm{\lambda}}^{k}-\bm{\lambda}^{*}\|^{2}+\frac{4\gamma^{k}}{\beta H^{2}u_{f}^{2}J}\|\Lambda^{k}-\frac{1}{H}\bm{1}\bm{1}^{\top}\Lambda^{k}\|_{F}^{2}\notag\\
&+ \frac{8(4\rho^{2}H+\chi^{2})J}{\gamma^{k} \beta H^{2}}\|\Lambda^{k+\frac{1}{2}}-\frac{1}{H}\bm{1}\bm{1}^{\top}\Lambda^{k+\frac{1}{2}}
\|_{F}^{2}.\notag
\end{align}

Substituting \eqref{lemma3-2} in Lemma \ref{l3} into \eqref{eq_proof_theorem2-b-3} and rearranging the terms, we obtain
\begin{align}
\label{eq_proof_theorem2-b-4}
&\|\bar{\bm{\lambda}}^{k+1}-\bm{\lambda}^{*}\|^{2}\\
\le & (1-\frac{\gamma^{k}\beta\cdot\frac{H}{J}}{4})\|\bar{\bm{\lambda}}^{k}-\bm{\lambda}^{*}\|^{2}+\frac{4\gamma^{k}}{\beta H^{2}u_{f}^{2}J}\|\Lambda^{k}-\frac{1}{H}\bm{1}\bm{1}^{\top}\Lambda^{k}\|_{F}^{2}\notag\\
&+ \frac{24(4\rho^{2}H+\chi^{2})J}{\gamma^{k} \beta H^{2}}\|\Lambda^{k}-\frac{1}{H}\bm{1}\bm{1}^{\top}\Lambda^{k}
\|_{F}^{2}\notag\\
&+ \frac{48\gamma^{k}\delta^{2}(4\rho^{2}H+\chi^{2})H}{\beta J}.\notag
\end{align}
Based on Lemma \ref{l4}, we can rewrite \eqref{eq_proof_theorem2-b-4} as
\begin{align}
\label{eq_proof_theorem2-b-5}
&\|\bar{\bm{\lambda}}^{k+1}-\bm{\lambda}^{*}\|^{2}\\
\le & (1-\frac{\gamma^{k}\beta\cdot\frac{H}{J}}{4})\|\bar{\bm{\lambda}}^{k}-\bm{\lambda}^{*}\|^{2}+\frac{72(\gamma^{k})^{3}\delta^{2}H}{\beta\epsilon^{3}u_{f}^{2}J^{3}}\notag\\
&+\frac{432\gamma^{k}(4\rho^{2}H+\chi^{2})\delta^{2}H}{\beta\epsilon^{3} J} + \frac{48\gamma^{k}\delta^{2}(4\rho^{2}H+\chi^{2})H}{\beta J}.\notag
\end{align}
Set a proper decreasing step size $\gamma^{k}=\frac{4}{\beta\cdot\frac{H}{J}(k_{0}+k)}$, where $k_{0}> 1$ is a any positive integer. Thus, $1-\frac{\gamma^{k} \beta\cdot\frac{H}{J}}{4}=1-\frac{1}{k_{0}+k}$ and \eqref{eq_proof_theorem2-b-5} can be rewritten as
\begin{align}
\label{eq_proof_theorem2-b-6}
&\|\bar{\bm{\lambda}}^{k+1}-\bm{\lambda}^{*}\|^{2}\\
\le & (1-\frac{1}{k_{0}+k})\|\bar{\bm{\lambda}}^{k}-\bm{\lambda}^{*}\|^{2}+ \frac{1}{(k_{0}+k)^{3}}\cdot \frac{4608\delta^{2}}{\beta^{4}\epsilon^{3}u_{f}^{2}H^{2}}\notag\\
&+\frac{1}{k_{0}+k}\cdot\frac{192\delta^{2}}{\beta^{2}}\cdot(1+\frac{9}{\epsilon^{3}})\cdot(4\rho^{2}H+\chi^{2}).\notag
\end{align}
Then we rewrite \eqref{eq_proof_theorem2-b-6} recursively and obtain
\begin{align}
\label{eq_proof_theorem2-b-7}
&\|\bar{\bm{\lambda}}^{k+1}-\bm{\lambda}^{*}\|^{2}\\
& \le \underbrace{\prod_{k'=0}^{k}(1-\frac{1}{k_{0}+k-k'})}_{T_{4}}\|\bar{\bm{\lambda}}^{0}-\bm{\lambda}^{*}\|^{2}\notag\\
&+[\prod_{k'=0}^{k-1}(1-\frac{1}{k_{0}+k-k'})\frac{1}{(k_{0}+0)^{3}}+\cdots+\notag\\
&\underbrace{\prod_{k'=0}^{0}\frac{(1-\frac{1}{k_{0}+k-k'})}{(k_{0}+k-1)^{3}}+\frac{1}{(k_{0}+k)^{3}}}_{T_{5}}]\cdot \frac{4608\delta^{2}}{\beta^{4}\epsilon^{3}u_{f}^{2}H^{2}}\notag\\
&+[\prod_{k'=0}^{k-1}(1-\frac{1}{k_{0}+k-k'})\frac{1}{k_{0}+0}+\cdots+\notag\\
&\underbrace{\prod_{k'=0}^{0}\frac{(1-\frac{1}{k_{0}+k-k'})}{k_{0}+k-1}+\frac{1}{k_{0}+k}]}_{T_{6}}\cdot\frac{192\delta^{2}(4\rho^{2}H+\chi^{2})(1+\frac{9}{\epsilon^{3}})}{\beta^{2}}.\notag
\end{align}
Next, we analyze $T_{4}$, $T_{5}$ and $T_{6}$ in turn.

For $T_{4}$, we have
\begin{align}
\label{eq_proof_theorem2-b-T4}
T_{4}=&\prod_{k'=0}^{k}(1-\frac{1}{k_{0}+k-k'})\\
=&\frac{k_{0}+k-1}{k_{0}+k}\cdot\frac{k_{0}+k-2}{k_{0}+k-1}\cdot\cdots\cdot\frac{k_{0}-1}{k_{0}}\notag\\
=&\frac{k_{0}-1}{k_{0}+k}.\notag
\end{align}
For $T_{5}$, we have
\begin{align}
\label{eq_proof_theorem2-b-T5}
T_{5}=&\frac{k_{0}}{k_{0}+k}\cdot\frac{1}{(k_{0}+0)^{3}}+\frac{k_{0}+1}{k_{0}+k}\cdot\frac{1}{(k_{0}+1)^{3}}+\cdots\\
&+\frac{k_{0}+k-1}{k_{0}+k}\cdot\frac{1}{(k_{0}+k-1)^{3}}+\frac{1}{(k_{0}+k)^{3}}\notag\\
=& \frac{1}{k_{0}+k}[\frac{1}{(k_{0})^{2}}+\frac{1}{(k_{0}+1)^{2}}+\cdots\notag\\
&+\frac{1}{(k_{0}+k-1)^{2}}+\frac{1}{(k_{0}+k)^{2}}]\notag\\
\le& \frac{1}{k_{0}+k}\cdot\frac{1}{k_{0}-1}.\notag
\end{align}
To drive the last inequality, we use $\sum_{k'=k_{0}}^{k}\frac{1}{(k')^{2}}\le \frac{1}{k_{0}-1}$.
For $T_{6}$, we have
\begin{align}
\label{eq_proof_theorem2-b-T6}
T_{6}=&[\frac{k_{0}}{k_{0}+k}\cdot\frac{1}{k_{0}+0}+\frac{k_{0}+1}{k_{0}+k}\cdot\frac{1}{k_{0}+1}+\cdots\\
&+\frac{k_{0}+k-1}{k_{0}+k}\cdot\frac{1}{k_{0}+k-1}+\frac{1}{k_{0}+k}]\notag\\
=& \frac{k+1}{k_{0}+k}.\notag
\end{align}

Substituting \eqref{eq_proof_theorem2-b-T4}, \eqref{eq_proof_theorem2-b-T5} and \eqref{eq_proof_theorem2-b-T6}  into \eqref{eq_proof_theorem2-b-7}, we have
\begin{align}
\label{eq_proof_theorem2-b-8}
&\|\bar{\bm{\lambda}}^{k+1}-\bm{\lambda}^{*}\|^{2}\\
\le&\frac{k_{0}-1}{k_{0}+k}\|\bar{\bm{\lambda}}^{0}-\bm{\lambda}^{*}\|^{2}+\frac{1}{(k_{0}+k)(k_{0}-1)}\cdot\frac{4608\delta^{2}}{\beta^{4}\epsilon^{3}u_{f}^{2}H^{2}}\notag\\
&+\frac{k+1}{k_{0}+k}\cdot \frac{192\delta^{2}}{\beta^{2}}\cdot(1+\frac{9}{\epsilon^{3}})\cdot(4\rho^{2}H+\chi^{2}).\notag
\end{align}
Applying the triangle inequality into \eqref{eq_proof_theorem2-b-8}, we obtain
\begin{align}
\label{eq_proof_theorem2-b-9}
&\|\bar{\bm{\lambda}}^{k+1}-\bm{\lambda}^{*}\|\\
\le&\sqrt{\frac{k_{0}-1}{k_{0}+k}}\|\bar{\bm{\lambda}}^{0}-\bm{\lambda}^{*}\|+\sqrt{\frac{1}{(k_{0}+k)(k_{0}-1)}\cdot\frac{4608\delta^{2}}{\beta^{4}\epsilon^{3}u_{f}^{2}H^{2}}}\notag\\
&+\sqrt{\frac{k+1}{k_{0}+k}\cdot \frac{192\delta^{2}}{\beta^{2}}\cdot(1+\frac{9}{\epsilon^{3}})\cdot(4\rho^{2}H+\chi^{2})}.\notag
\end{align}
By Lemma \ref{l4} and using the step size $\gamma^{k}=\frac{4}{\beta\cdot\frac{H}{J}(k_{0}+k)}$, we obtain
\begin{align}
\label{eq_proof_theorem2-b-10}
\sum_{i\in \mathcal{H}}\|\bm{\lambda}_{i}^{k+1}-\bar{\bm{\lambda}}^{k+1}\|^{2}=&\|\Lambda^{k+1}-\frac{1}{|\mathcal{H}|}\bm{1}\bm{1}^{\top}\Lambda^{k+1}\|^{2}_{F}\\
\le &\frac{18(\gamma^{k+1})^{2}\delta^{2}H^{3}}{\epsilon^{3} \cdot J^{2}}\notag\\
=& \frac{288\delta^{2}H}{(k_{0}+k+1)^{2}\epsilon^{3}\beta^{2}}.\notag
\end{align}
By $(a_{1}+\cdots+a_{H})^{2}\le H(a_{1}^{2}+\cdots+a_{H}^{2})$ and \eqref{eq_proof_theorem2-b-10}, we have
\begin{align}
\label{eq_proof_theorem2-b-11}
\sum_{i\in \mathcal{H}}\|\bm{\lambda}_{i}^{k+1}-\bar{\bm{\lambda}}^{k+1}\|=&\sqrt{H}\cdot\sqrt{\sum_{i\in \mathcal{H}}\|\bm{\lambda}_{i}^{k+1}-\bar{\bm{\lambda}}^{k+1}\|^{2}}\\
\le &\sqrt{H}\cdot \sqrt{\frac{288\delta^{2}H}{(k_{0}+k+1)^{2}\epsilon^{3}\beta^{2}}}.\notag
\end{align}

Combining \eqref{eq_proof_theorem2-b-9} and \eqref{eq_proof_theorem2-b-11} yields
\begin{align}
\label{eq_proof_theorem2-b-12}
&\sum_{i\in \mathcal{H}}\|\bm{\lambda}_{i}^{k+1}-\bm{\lambda}^{*}\|\\
\le& \sum_{i\in \mathcal{H}}\|\bm{\lambda}_{i}^{k+1}-\bar{\bm{\lambda}}^{k+1}\|+H\cdot\|\bar{\bm{\lambda}}^{k+1}-\bm{\lambda}^{*}\|\notag\\
\le&H\cdot\sqrt{\frac{k_{0}-1}{k_{0}+k}}\|\bar{\bm{\lambda}}^{0}-\bm{\lambda}^{*}\|+\sqrt{\frac{1}{(k_{0}+k)(k_{0}-1)}\cdot\frac{4608\delta^{2}}{\beta^{4}\epsilon^{3}u_{f}^{2}}}\notag\\
&+H\cdot\sqrt{\frac{k+1}{k_{0}+k}\cdot \frac{192\delta^{2}}{\beta^{2}}\cdot(1+\frac{9}{\epsilon^{3}})\cdot(4\rho^{2}H+\chi^{2})}\notag\\
&+\sqrt{H}\cdot \sqrt{\frac{288\delta^{2}H}{(k_{0}+k+1)^{2}\epsilon^{3}\beta^{2}}}.\notag
\end{align}
Taking $k \to +\infty$, we obtain
\begin{align}
\label{eq_proof_eq22}
&\limsup_{k\to +\infty}\sum_{i\in \mathcal{H}}\|\bm{\lambda}_{i}^{k+1}-\bm{\lambda}^{*}\|\\
\le&\sqrt{\frac{192\delta^{2}H^{2}}{\beta^{2}}\cdot(1+\frac{9}{\epsilon^{3}})\cdot(4\rho^{2}H+\chi^{2})}.\notag
\end{align}
And the dual variables in our proposed Byzantine-resilient algorithm asymptotically converge to a neighborhood of the optimum at a rate of $O(\frac{1}{\sqrt{k}})$.

\subsection{Part $b$ of Theorem \ref{t2}}
Based on $\|\Lambda^{k+1}-\frac{1}{H}\bm{1}\bm{1}^{\top}\Lambda^{k+1}
\|^{2}_{F}\le \frac{18(\gamma^{k+1})^{2}\delta^{2}H^{3}}{\epsilon^{3} \cdot J^{2}}$ in Lemma \ref{l4} and the fact $\|\Lambda^{k+1}-\frac{1}{H}\bm{1}\bm{1}^{\top}\Lambda^{k+1}
\|^{2}_{F}=\sum_{i\in \mathcal{H}}\|\bm{\lambda}_{i}^{k+1}-\bar{\bm{\lambda}}^{k+1}\|^{2}$, with a proper decreasing step size $\gamma^{k}=O(\frac{1}{k})$, taking $k \to +\infty$, we obtain
\begin{align}
\label{eq-proof-theorem2-c-1}
\lim_{k\to +\infty}\sum_{i\in \mathcal{H}}\|\bm{\lambda}_{i}^{k+1}-\bar{\bm{\lambda}}^{k+1}\|= 0.
\end{align}
We conclude by summarizing the conditions on the step size $\gamma^{k}$ in Theorem \ref{t2}. It must satisfy $(\gamma^{k})^{2}\cdot\frac{H}{J}-2\gamma^{k}\alpha\le 0$, $\frac{18(\gamma^{k})^{2}}{\epsilon u_{f}^{2}J^{2}}\le \frac{(2-\epsilon)\epsilon^{2}}{3(3-\epsilon)}$, $1\le \frac{(\gamma^{k})^{2}}{(\gamma^{k+1)^{2}}}\le \frac{2}{1+(1-\epsilon^{2})}$, as well as $\gamma^{k}\le \frac{u_{f}J}{2\sqrt{3}}$. The specific step size $\gamma^{k}=\frac{4}{\beta\cdot\frac{H}{J}(k_{0}+k)}$ with $k_{0}\ge \max\{\frac{2}{\alpha\beta},\sqrt{\frac{216(3-\epsilon)}{(2-\epsilon)\epsilon u_{f}^{2}H^{2}\beta^{2}}},\frac{1}{\sqrt{\frac{2}{1+(1-\epsilon^{2})}}-1},\frac{8\sqrt{3}}{u_{f}H \beta}\}$ satisfies these conditions.

\subsection{Part $c$ of Theorem \ref{t2}}
The Lagrangian function of \eqref{eq_oracle-problem} is
\begin{align}
\label{eq_Lagrangian-oracle-problem}
\begin{split}
\mathcal{L}\left ( \bm{\Theta} ;  \bm{\lambda} \right ):=
\frac{1}{H}\sum_{i\in \mathcal{H}}f_{i}\left ( \bm{\theta}_{i} \right )+ \left \langle \bm{\lambda }, \frac{1}{H}\sum_{i\in \mathcal{H}}\bm{\theta}_{i}-\bm{s} \right \rangle.
\end{split}
\end{align}

Since $\bm{\Theta}^{*}$ is the optimal solution of the primal problem \eqref{eq_oracle-problem}, we have $ \frac{1}{H}\sum_{i\in \mathcal{H}}\bm{\theta}_{i}^{*} = \bm{s}$ and $\bm{\theta}^{*}\in C$. According to \eqref{eq_Lagrangian-oracle-problem}, for any dual variable $\bm{\lambda}$ we have
\begin{align}
\label{eq_proof_theorem2-a-1}
\mathcal{L}\left ( \bm{\Theta}^{*} ;  \bm{\lambda} \right )=&\frac{1}{H}\sum_{i\in \mathcal{H}}f_{i}\left ( \bm{\theta}_{i}^{*} \right )+ \left \langle \bm{\lambda }, \frac{1}{H}\sum_{i\in \mathcal{H}}\bm{\theta}_{i}^{*}-\bm{s} \right \rangle\\
=& f(\bm{\Theta}^{*}). \notag
\end{align}
By Assumption \ref{a2}, the duality gap is zero. According to \eqref{eq_dual-problem_oracle}, for any $\bm{\lambda}$ we obtain $g(\bm{\lambda}^{*})=-f(\bm{\Theta}^{*})=-\mathcal{L}\left ( \bm{\Theta}^{*} ;  \bm{\lambda} \right )$. Now we introduce a vector $^{\dagger}\bm{\Theta}^{k+1}:=[_{}^{\dagger}\bm{\theta}_{1}^{k+1}; \cdots;  _{}^{\dagger}\bm{\theta}_{H}^{k+1}]$, where $_{}^{\dagger}\bm{\theta}_{i}^{k+1}=\arg\min_{\bm{\theta}_{i}\in C_{i}} \{\bm{\theta}_{i}^{\top}\bar{\bm{\lambda}}^{k+1}+f_{i}(\bm{\theta}_{i})\}$ and $\bar{\bm{\lambda}}^{k+1}:=\frac{1}{H}\sum_{i\in \mathcal{H}}\bm{\lambda}_{i}^{k+1}$. Therefore, we have
\begin{align}
\label{eq_proof_theorem2-a-2}
&g(\bar{\bm{\lambda}}^{k+1})-g(\bm{\lambda}^{*})\\
=&-\inf_{\bm{\Theta}\in C}\mathcal{L}\left ( \bm{\Theta} ;  \bar{\bm{\lambda}}^{k+1} \right )+\mathcal{L}\left ( \bm{\Theta}^{*} ;  \bar{\bm{\lambda}}^{k+1} \right )\notag\\
=&-\mathcal{L}\left ( ^{\dagger}\bm{\Theta}^{k+1} ;  \bar{\bm{\lambda}}^{k+1} \right )+\mathcal{L}\left ( \bm{\Theta}^{*} ;  \bar{\bm{\lambda}}^{k+1} \right ).\notag
\end{align}

Assumption \ref{a1} shows that the local cost function $f_{i}(\cdot)$ is $u_{f}$-strongly convex. Further using the definition of $\mathcal{L}\left ( \bm{\Theta} ;  \bm{\lambda} \right )$ in \eqref{eq_Lagrangian-oracle-problem}, we know that $\mathcal{L}\left ( \bm{\Theta} ;  \bm{\lambda} \right )$ is $u_{f}$-strongly convex with respect to $\bm{\Theta}$. Therefore, we have
\begin{align}
\label{eq_proof_theorem2-a-3}
&\mathcal{L}\left ( \bm{\Theta}^{*} ;  \bar{\bm{\lambda}}^{k+1} \right )-\mathcal{L}\left ( ^{\dagger}\bm{\Theta}^{k+1} ;  \bar{\bm{\lambda}}^{k+1} \right )\\
\ge &\nabla \mathcal{L}^{\top}\left ( ^{\dagger}\bm{\Theta}^{k+1} ;  \bar{\bm{\lambda}}^{k+1} \right )(\bm{\Theta}^{*}-\ ^{\dagger}\bm{\Theta}^{k+1})+\frac{u_{f}}{2}\| ^{\dagger}\bm{\Theta}^{k+1}-\bm{\Theta}^{*}\|^{2}.\notag
\end{align}
Combining \eqref{eq_proof_theorem2-a-2} and \eqref{eq_proof_theorem2-a-3}, we obtain
\begin{align}
\label{eq_proof_theorem2-a-4}
&g(\bar{\bm{\lambda}}^{k+1})-g(\bm{\lambda}^{*})\\
\ge& \nabla \mathcal{L}^{\top}\left ( ^{\dagger}\bm{\Theta}^{k+1} ;  \bar{\bm{\lambda}}^{k+1} \right )(\bm{\Theta}^{*}-\ ^{\dagger}\bm{\Theta}^{k+1})+\frac{u_{f}}{2}\| ^{\dagger}\bm{\Theta}^{k+1}-\bm{\Theta}^{*}\|^{2}\notag\\
\ge&\frac{u_{f}}{2}\| ^{\dagger}\bm{\Theta}^{k+1}-\bm{\Theta}^{*}\|^{2}.\notag
\end{align}
To drive the last inequality, we use optimality condition of $_{}^{\dagger}\bm{\theta}_{i}^{k+1}=\arg\min_{\bm{\theta}_{i}\in C_{i}} \{\bm{\theta}_{i}^{\top}\bar{\bm{\lambda}}^{k+1}+f_{i}(\bm{\theta}_{i})\}$ \cite[Proposition 2.1.2]{b-Dimitri-P.-Bertsekas-1999}.

According to Lemma \ref{l2}, $g_{i}(\bm{\lambda})$ is smooth with constant $\frac{1}{Hu_{f}}$. Therefore, function $g(\bm{\lambda })=\sum_{i\in\mathcal{H}}g_{i}(\bm{\lambda })$ is smooth with constant $\frac{1}{u_{f}}$. This fact leads to
\begin{align}
\label{eq_proof_theorem2-a-5}
&g(\bar{\bm{\lambda}}^{k+1})-g(\bm{\lambda}^{*})\\
\le& \nabla g^{\top}(\bm{\lambda}^{*})(\bar{\bm{\lambda}}^{k+1}-\bm{\lambda}^{*})+\frac{1}{2u_{f}}\|\bar{\bm{\lambda}}^{k+1}-\bm{\lambda}^{*}\|^{2}\notag\\
=& \frac{1}{2u_{f}}\|\bar{\bm{\lambda}}^{k+1}-\bm{\lambda}^{*}\|^{2}.\notag
\end{align}
To drive the last equality, we use the fact that $\nabla g(\bm{\lambda}^{*})=\bm{0}$.
Combining \eqref{eq_proof_theorem2-a-4} and \eqref{eq_proof_theorem2-a-5}, we have
\begin{align}
\label{eq_proof_theorem2-a-6}
\|^{\dagger}\bm{\Theta}^{k+1}-\bm{\Theta}^{*}\|^{2} \le \frac{1}{(u_{f})^{2}}\|\bar{\bm{\lambda}}^{k+1}-\bm{\lambda}^{*}\|^{2}.
\end{align}
Combining \eqref{eq_proof_theorem2-b-9} and \eqref{eq_proof_theorem2-a-6}, we obtain
\begin{align}
\label{eq_proof_theorem2-a-7}
&\| ^{\dagger}\bm{\Theta}^{k+1}-\bm{\Theta}^{*}\| \\
\le&\frac{1}{u_{f}}\cdot[\sqrt{\frac{k_{0}-1}{k_{0}+k}}\|\bar{\bm{\lambda}}^{0}-\bm{\lambda}^{*}\|+\sqrt{\frac{1}{(k_{0}+k)(k_{0}-1)}\cdot\frac{288\delta^{2}}{\beta^{4}\epsilon^{3}u_{f}^{2}H^{2}}}\notag\\
&+\sqrt{\frac{k+1}{k_{0}+k}\cdot \frac{192\delta^{2}}{\beta}\cdot(1+\frac{3}{\epsilon^{3}})\cdot(4\rho^{2}H+\chi^{2})}].\notag
\end{align}
Taking $k\to +\infty$, we obtain
\begin{align}
\label{eq_proof_theorem2-a-8}
&\limsup_{k\to +\infty}\| ^{\dagger}\bm{\Theta}^{k+1}-\bm{\Theta}^{*}\|\\
\le& \frac{1}{u_{f}}\cdot \sqrt{ \frac{192\delta^{2}}{\beta}\cdot(1+\frac{3}{\epsilon^{3}})\cdot(4\rho^{2}H+\chi^{2})}. \notag
\end{align}

According to Assumption \ref{a1}, $f_{i}(\bm{\theta}_{i})$ is $u_{f}$-strongly convex. By the conjugate correspondence theorem in \cite{b-Amir-Beck-2017}, the conjugate function $F_{i}^{*}(\bm{\lambda})=\max_{\bm{\theta}_{i}\in C_{i}}\{\bm{\lambda}^{\top}\bm{\theta}_{i}-f_{i}(\bm{\theta}_{i})\}$ is $\frac{1}{u_{f}}$-smooth. In consequence, the gradient $\nabla F^{*}_{i}(-\bm{\lambda})=\arg\min_{\bm{\theta}_{i}\in C_{i}}\{\bm{\lambda }^{\top}\bm{\theta}_{i}+f_{i}(\bm{\theta}_{i})\}$ is $\frac{1}{u_{f}}$-Lipschitz continuous. According to the definition of Lipschitz continuity, we have
\begin{align}
\label{eq_proof_theorem2-a-lipschitz-continuous}
\|\nabla F^{*}_{i}(\bm{\lambda}_{i}^{k+1})-\nabla F^{*}_{i}(\bar{\bm{\lambda}}^{k+1})\| \le \frac{1}{u_{f}}\|\bm{\lambda}_{i}^{k+1}-\bar{\bm{\lambda}}^{k+1}\|.
\end{align}
Based on $_{}^{\dagger}\bm{\theta}_{i}^{k+1}=\arg\min_{\bm{\theta}_{i}\in C_{i}} \{\bm{\theta}_{i}^{\top}\bar{\bm{\lambda}}^{k+1}+f_{i}(\bm{\theta}_{i})\}$ and $\bm{\theta}_{i}^{k+1}=\arg\min_{\bm{\theta}_{i}\in C_{i}} \{\bm{\theta}_{i}^{\top}\bm{\lambda}_{i}^{k+1}+f_{i}(\bm{\theta}_{i})\}$, we obtain $\bm{\theta}_{i}^{k+1}=\nabla F^{*}_{i}(\bm{\lambda}_{i}^{k+1})$ and $_{}^{\dagger}\bm{\theta}_{i}^{k+1}=\nabla F^{*}_{i}(\bar{\bm{\lambda}}^{k+1})$. Substituting them into \eqref{eq_proof_theorem2-a-lipschitz-continuous}, we have
\begin{align}
\label{eq_proof_theorem2-a-10}
\|\bm{\theta}_{i}^{k+1}-\ _{}^{\dagger}\bm{\theta}_{i}^{k+1}\| \le \frac{1}{u_{f}}\|\bm{\lambda}_{i}^{k+1}-\bar{\bm{\lambda}}^{k+1}\|.
\end{align}
Combining \eqref{eq_proof_theorem2-a-10}, $^{\dagger}\bm{\Theta}^{k+1}:=[_{}^{\dagger}\bm{\theta}_{1}^{k+1}; \cdots;  _{}^{\dagger}\bm{\theta}_{H}^{k+1}]$ and $\bm{\Theta}^{k+1}:=[\bm{\theta}_{1}^{k+1}; \cdots; \bm{\theta}_{H}^{k+1}]$, we obtain
\begin{align}
\label{eq_proof_theorem2-a-11}
\|\bm{\Theta}^{k+1}-\ ^{\dagger}\bm{\Theta}^{k+1}\| \le \frac{1}{u_{f}}\sum_{i\in \mathcal{H}}\|\bm{\lambda}_{i}^{k+1}-\bar{\bm{\lambda}}^{k+1}\|.
\end{align}
Combining \eqref{eq_proof_theorem2-b-11} and  \eqref{eq_proof_theorem2-a-11}, we have
\begin{align}
\label{eq_proof_theorem2-a-12}
\|\bm{\Theta}^{k+1}-\ ^{\dagger}\bm{\Theta}^{k+1}\| \le \frac{\sqrt{H}}{u_{f}}\cdot \sqrt{\frac{288\delta^{2}H}{(k_{0}+k+1)^{2}\epsilon^{3}\beta^{2}}}.
\end{align}
Taking $k\to +\infty$, we obtain
\begin{align}
\label{eq_proof_theorem2-a-13}
\lim_{k\to +\infty}\|\bm{\Theta}^{k+1}-\ ^{\dagger}\bm{\Theta}^{k+1}\| =0.
\end{align}

Combining \eqref{eq_proof_theorem2-a-8} and \eqref{eq_proof_theorem2-a-13} yields
\begin{align}
\label{eq_proof_theorem2-a-9}
&\limsup_{k\to +\infty}\|\bm{\Theta}^{k+1}-\bm{\Theta}^{*}\|\\
\le &\limsup_{k\to +\infty}\|\bm{\Theta}^{k+1}-\ ^{\dagger}\bm{\Theta}^{k+1}\|+\limsup_{k\to +\infty}\|^{\dagger}\bm{\Theta}^{k+1}-\bm{\Theta}^{*}\|\notag\\
\le& \frac{1}{u_{f}}\cdot \sqrt{ \frac{192\delta^{2}}{\beta}\cdot(1+\frac{3}{\epsilon^{3}})\cdot(4\rho^{2}H+\chi^{2})}. \notag
\end{align}
And the primal variables in our proposed Byzantine-resilient algorithm asymptotically converge to a neighborhood of the optimum at a rate of $O(\frac{1}{\sqrt{k}})$.

\subsection{Supporting Lemmas}
\begin{lemma}
\label{l1}
Under Assumption \ref{a1}, for any $\bm{\lambda} \in \mathbb{R}^D$, the maximum distance between the honest agents' local dual gradients and their average, denoted by $\max_{i\in \mathcal{H}}\|\nabla g_{i}(\bm{\lambda}) - \frac{1}{H} \sum_{i\in \mathcal{H}} \nabla g_{i}(\bm{\lambda}) \|^{2}$, is bounded by some positive constant $\delta^2$.

\textbf{Proof}.
Recalling the definition of the local dual gradient $\nabla g_{i}(\bm{\lambda})=-\frac{1}{H}\arg\min_{\bm{\theta}_{i}\in C_{i}} \{\bm{\lambda}^{\top}\bm{\theta}_{i}+f_{i}(\bm{\theta}_{i})\}+\frac{1}{H}\bm{s}$, we have
\begin{align}
&\max_{i\in \mathcal{H}}\|\nabla g_{i}(\bm{\lambda}) - \frac{1}{H} \sum_{i\in \mathcal{H}} \nabla g_{i}(\bm{\lambda}) \|^{2}\\
=&\max_{i\in \mathcal{H}}\|-\frac{1}{H}\arg\underset{\bm{\theta}_{i}\in C_{i}}{\min}\{\bm{\lambda}^{\top}\bm{\theta}_{i}+f_{i}(\bm{\theta}_{i})\}+\frac{1}{H}\bm{s}\notag\\
&+\frac{1}{H^{2}}\sum_{i\in \mathcal{H}}\arg\underset{\bm{\theta}_{i}\in C_{i}}{\min}\{\bm{\lambda}^{\top}\bm{\theta}_{i}+f_{i}(\bm{\theta}_{i})\}-\frac{1}{H}\bm{s}\|^{2}\notag
\end{align}
\begin{align}
=&\max_{i\in \mathcal{H}}\|\frac{1}{H}(\frac{1}{H}\sum_{i\in \mathcal{H}}\arg\underset{\bm{\theta}_{i}\in C_{i}}{\min}\{\bm{\lambda}^{\top}\bm{\theta}_{i}+f_{i}(\bm{\theta}_{i})\}\notag\\
&-\arg\underset{\bm{\theta}_{i}\in C_{i}}{\min}\{\bm{\lambda}^{\top}\bm{\theta}_{i}+f_{i}(\bm{\theta}_{i})\})\|^{2}.\notag
\end{align}
According to Assumption \ref{a1}, the local constraint sets $C_{i}$ are bounded by hypothesis, and we know that $\max_{i\in \mathcal{H}}\|\nabla g_{i}(\bm{\lambda}) - \frac{1}{H} \sum_{i\in \mathcal{H}} \nabla g_{i}(\bm{\lambda}) \|^{2}$ is also bounded by some positive constant, which we denoted as $\delta^{2}$.
\end{lemma}

\begin{lemma}
\label{l2}
Under Assumption \ref{a1}, for any honest agent $i \in \mathcal{H}$, the local dual function $g_{i}(\bm{\lambda})$ is strongly convex with constant $\frac{1}{HL_{f}}$ and smooth with constant $\frac{1}{Hu_{f}}$.

\textbf{Proof}.
According to Assumption \ref{a1}, $f_{i}(\cdot)$ is $u_{f}$-strongly convex and $L_{f}$-smooth. By the conjugate correspondence theorem \cite{b-Amir-Beck-2017}, its conjugate function $F_{i}^{*}(\bm{\lambda})=\max_{\bm{\theta}_{i}\in C_{i}} \{\bm{\lambda}^{\top}\bm{\theta}_{i}-f_{i}(\bm{\theta}_{i})\}$ is $\frac{1}{L_{f}}$-strongly convex and $\frac{1}{u_{f}}$-smooth. By the definition $g_{i}(\bm{\lambda})=\frac{1}{H}F_{i}^{*}(-\bm{\lambda })+\frac{1}{H}\bm{\lambda }^{\top}\bm{s}$ $(\forall i \in \mathcal{H})$, we know $g_{i}(\bm{\lambda})$ is $\frac{1}{HL_{f}}$-strongly convex and $\frac{1}{Hu_{f}}$-smooth.
\end{lemma}

\begin{lemma}
\label{l3}
Define a matrix $\Lambda^{k+\frac{1}{2}}=[\cdots,\bm{\lambda}_{i}^{k+\frac{1}{2}},\cdots]\in \mathbb{R}^{H\times D} $ that collects the dual variables $\bm{\lambda}_{i}^{k+\frac{1}{2}}$ of all honest agents $i\in \mathcal{H}$ generated by Algorithm \ref{alg2}. Under Assumption \ref{a1}, we have
\begin{align}
\label{lemma3-1}
&\|\Lambda^{k+\frac{1}{2}}-\frac{1}{H}\bm{1}\bm{1}^{\top}\Lambda^{k+\frac{1}{2}}
\|^{2}_{F}\\
\le&(\frac{1}{1-v}+\frac{6(\gamma^{k})^{2}}{v\cdot u_{f}^{2}J^{2}})\|\Lambda^{k}-\frac{1}{H}\bm{1}\bm{1}^{\top}\Lambda^{k}
\|^{2}_{F}+\frac{3(\gamma^{k})^{2}\delta^{2}H^{3}}{v\cdot J^{2}},\notag
\end{align}
where $v$ is any positive constant in $(0,1)$. If $v=\frac{1}{2}$ and the step size $\gamma^{k}\le \frac{u_{f}J}{2\sqrt{3}}$, this further yields
\begin{align}
\label{lemma3-2}
&\|\Lambda^{k+\frac{1}{2}}-\frac{1}{H}\bm{1}\bm{1}^{\top}\Lambda^{k+\frac{1}{2}}
\|^{2}_{F}\\
\le&3\|\Lambda^{k}-\frac{1}{H}\bm{1}\bm{1}^{\top}\Lambda^{k}
\|^{2}_{F}+\frac{6(\gamma^{k})^{2}\delta^{2}H^{3}}{ J^{2}}.\notag
\end{align}

\textbf{Proof}.
According to the update of $\bm{\lambda}_{i}^{k+\frac{1}{2}}$ in Algorithm \ref{alg2} and the fact $\|\Lambda^{k+\frac{1}{2}}-\frac{1}{H}\bm{1}\bm{1}^{\top}\Lambda^{k+\frac{1}{2}}
\|^{2}_{F}=\sum_{i\in \mathcal{H}}\|\bm{\lambda}_{i}^{k+\frac{1}{2}}-\bar{\bm{\lambda}}^{k+\frac{1}{2}}\|^{2}$, we have
\begin{align}
\label{eq-proof-lemma3-1-1}
&\|\Lambda^{k+\frac{1}{2}}-\frac{1}{H}\bm{1}\bm{1}^{\top}\Lambda^{k+\frac{1}{2}}
\|^{2}_{F}=\sum_{i\in \mathcal{H}}\|\bm{\lambda}_{i}^{k+\frac{1}{2}}-\bar{\bm{\lambda}}^{k+\frac{1}{2}}\|^{2} \\
=&\sum_{i\in \mathcal{H}}\|\bm{\lambda}_{i}^{k}-\gamma^{k}\cdot \frac{H}{J}\nabla g_{i}(\bm{\lambda}_{i}^{k})-\bar{\bm{\lambda}}^{k}+\gamma^{k}\cdot\frac{1}{J}\sum_{i\in \mathcal{H}}\nabla g_{i}(\bm{\lambda}_{i}^{k})\|^{2}\notag\\
=& \sum_{i\in \mathcal{H}}\|\bm{\lambda}_{i}^{k}-\bar{\bm{\lambda}}^{k}-\gamma^{k}\cdot\frac{H}{J}(\nabla g_{i}(\bm{\lambda}_{i}^{k})-\frac{1}{H}\sum_{i\in \mathcal{H}}\nabla g_{i}(\bm{\lambda}_{i}^{k}))\|^{2}\notag\\
\le&\frac{(\gamma^{k})^{2}}{v}\cdot\frac{H^{2}}{J^{2}}\sum_{i\in \mathcal{H}}\|\nabla g_{i}(\bm{\lambda}_{i}^{k})-\frac{1}{H}\sum_{i\in \mathcal{H}}\nabla g_{i}(\bm{\lambda}_{i}^{k}))\|^{2} \notag\\
&+\frac{1}{1-v}\sum_{i\in \mathcal{H}}\|\bm{\lambda}_{i}^{k}-\bar{\bm{\lambda}}^{k}\|^{2}, \notag
\end{align}
where $v\in (0,1)$ is any positive constant. To drive the last inequality, we use the fact that $\|\bm{a}+\bm{b}\|^{2}\le \frac{1}{v}\|\bm{a}\|^{2}+\frac{1}{1-v}\|\bm{b}\|^{2}$ for any positive constant $v\in (0,1)$.

For $\frac{(\gamma^{k})^{2}}{v}\cdot\frac{H^{2}}{J^{2}}\sum_{i\in \mathcal{H}}\|\nabla g_{i}(\bm{\lambda}_{i}^{k})-\frac{1}{H}\sum_{i\in \mathcal{H}}\nabla g_{i}(\bm{\lambda}_{i}^{k}))\|^{2}$, the first term at the right-hand side of \eqref{eq-proof-lemma3-1-1}, we have
\begin{align}
\label{eq-proof-lemma3-1-2}
&\frac{(\gamma^{k})^{2}}{v}\cdot\frac{H^{2}}{J^{2}}\sum_{i\in \mathcal{H}}\|\nabla g_{i}(\bm{\lambda}_{i}^{k})-\frac{1}{H}\sum_{i\in \mathcal{H}}\nabla g_{i}(\bm{\lambda}_{i}^{k}))\|^{2}
\end{align}
\begin{align}
=&\frac{(\gamma^{k})^{2}}{v}\cdot\frac{H^{2}}{J^{2}}\sum_{i\in \mathcal{H}}\|\nabla g_{i}(\bm{\lambda}_{i}^{k})-\nabla g_{i}(\bar{\bm{\lambda}}^{k})+\nabla g_{i}(\bar{\bm{\lambda}}^{k})\notag\\
&-\frac{1}{H}\sum_{i\in \mathcal{H}}\nabla g_{i}(\bar{\bm{\lambda}}^{k}))+\frac{1}{H}\sum_{i\in \mathcal{H}}\nabla g_{i}(\bar{\bm{\lambda}}^{k}))-\frac{1}{H}\sum_{i\in \mathcal{H}}\nabla g_{i}(\bm{\lambda}_{i}^{k}))\|^{2} \notag\\
\le&\frac{3(\gamma^{k})^{2}}{v}\cdot\frac{H^{2}}{J^{2}}\sum_{i\in \mathcal{H}}\|\nabla g_{i}(\bm{\lambda}_{i}^{k})-\nabla g_{i}(\bar{\bm{\lambda}}^{k})\|^{2}\notag\\
&+\frac{3(\gamma^{k})^{2}}{v}\cdot\frac{H^{2}}{J^{2}}\sum_{i\in \mathcal{H}}\|\nabla g_{i}(\bar{\bm{\lambda}}^{k})-\frac{1}{H}\sum_{i\in \mathcal{H}}\nabla g_{i}(\bar{\bm{\lambda}}^{k}))\|^{2} \notag\\
&+\frac{3(\gamma^{k})^{2}}{v}\cdot\frac{H^{2}}{J^{2}}\sum_{i\in \mathcal{H}}\|\frac{1}{H}\sum_{i\in \mathcal{H}}\nabla g_{i}(\bar{\bm{\lambda}}^{k}))-\frac{1}{H}\sum_{i\in \mathcal{H}}\nabla g_{i}(\bm{\lambda}_{i}^{k}))\|^{2}.\notag
\end{align}
According to Lemmas \ref{l1} and \ref{l2}, we further have
\begin{align}
\label{eq-proof-lemma3-1-3}
&\frac{(\gamma^{k})^{2}}{v}\cdot\frac{H^{2}}{J^{2}}\sum_{i\in \mathcal{H}}\|\nabla g_{i}(\bm{\lambda}_{i}^{k})-\frac{1}{H}\sum_{i\in \mathcal{H}}\nabla g_{i}(\bm{\lambda}_{i}^{k}))\|^{2}\\
\le& \frac{6(\gamma^{k})^{2}}{v\cdot u_{f}^{2}J^{2}}\sum_{i\in\mathcal{H}}\|\bm{\lambda}_{i}^{k}-\bar{\bm{\lambda}}^{k}\|^{2}+\frac{3(\gamma^{k})^{2}\delta^{2}H^{3}}{v\cdot J^{2}}.\notag
\end{align}

Substituting \eqref{eq-proof-lemma3-1-3} into \eqref{eq-proof-lemma3-1-1} and then rearranging the terms, we obtain
\begin{align}
\label{eq-proof-lemma3-1}
&\|\Lambda^{k+\frac{1}{2}}-\frac{1}{H}\bm{1}\bm{1}^{\top}\Lambda^{k+\frac{1}{2}}
\|^{2}_{F}=\sum_{i\in \mathcal{H}}\|\bm{\lambda}_{i}^{k+\frac{1}{2}}-\bar{\bm{\lambda}}^{k+\frac{1}{2}}\|^{2}\\
\le& (\frac{1}{1-v}+\frac{6(\gamma^{k})^{2}}{v \cdot u_{f}^{2}J^{2}})\sum_{i\in \mathcal{H}}\|\bm{\lambda}_{i}^{k}-\bar{\bm{\lambda}}^{k}\|^{2}+\frac{3(\gamma^{k})^{2}\delta^{2}H^{3}}{v\cdot J^{2}}\notag\\
=&(\frac{1}{1-v}+\frac{6(\gamma^{k})^{2}}{v \cdot u_{f}^{2}J^{2}})\|\Lambda^{k}-\frac{1}{H}\bm{1}\bm{1}^{\top}\Lambda^{k}
\|^{2}_{F}+\frac{3(\gamma^{k})^{2}\delta^{2}H^{3}}{v\cdot J^{2}}.\notag
\end{align}
\end{lemma}

\begin{lemma}
\label{l4}
Define a matrix $\Lambda^{k+1}=[\cdots,\bm{\lambda}_{i}^{k+1},\cdots]\in \mathbb{R}^{H\times D} $ that collects the dual variables $\bm{\lambda}_{i}^{k+1}$ of all honest agents $i\in \mathcal{H}$ generated by Algorithm \ref{alg2}. Suppose that the robust aggregation rules $AGG_i$ satisfy \eqref{eq_definition1} in Definition \ref{d1}. Under Assumptions \ref{a1} and \ref{a3}, if the contraction constant $\rho$ satisfies
$\rho<\frac{1-\kappa}{8\sqrt{H}}$, we have
\begin{align}
\label{lemma4}
\|\Lambda^{k+1}-\frac{1}{H}\bm{1}\bm{1}^{\top}\Lambda^{k+1}
\|^{2}_{F}\le \frac{18(\gamma^{k+1})^{2}\delta^{2}H^{3}}{\epsilon^{3} \cdot J^{2}},
\end{align}
where $\epsilon:=1-\kappa-8\rho\sqrt{H}$.

\textbf{Proof}. For any positive constant $w\in (0,1)$, we have
\begin{align}
\label{eq-proof-lemma4-1}
 &\|\Lambda^{k+1}-\frac{1}{H}\bm{1}\bm{1}^{\top}\Lambda^{k+1}
\|^{2}_{F} \\
=& \|\Lambda^{k+1}-\frac{1}{H}\bm{1}\bm{1}^{\top}\Lambda^{k+1}+E\Lambda^{k+\frac{1}{2}}-E\Lambda^{k+\frac{1}{2}}\notag\\
&+\frac{1}{H}\bm{1}\bm{1}^{\top}E\Lambda^{k+\frac{1}{2}}-\frac{1}{H}\bm{1}\bm{1}^{\top}E\Lambda^{k+\frac{1}{2}}\|^{2}_{F} \notag\\
\le& \underbrace{\frac{1}{1-w}\|E\Lambda^{k+\frac{1}{2}}-\frac{1}{H}\bm{1}\bm{1}^{\top}E\Lambda^{k+\frac{1}{2}}\|^{2}_{F}}_{T_{7}}+\underbrace{\frac{2}{w}\|\Lambda^{k+1}-E\Lambda^{k+\frac{1}{2}}\|^{2}_{F}}_{T_{8}}\notag\\
&+\underbrace{\frac{2}{w}\|\frac{1}{H}\bm{1}\bm{1}^{\top}\Lambda^{k+1}-\frac{1}{H}\bm{1}\bm{1}^{\top}E\Lambda^{k+\frac{1}{2}}\|^{2}_{F}}_{T_{9}}. \notag
\end{align}
Next, we analyze $T_{7}$, $T_{8}$ and $T_{9}$ in turn.

\noindent\textbf{Bounding $T_{7}$:}
According to Assumption \ref{a3}, we have
\begin{align}
\label{eq-proof-lemma4-2}
T_{7}=&\frac{1}{1-w}\|E\Lambda^{k+\frac{1}{2}}-\frac{1}{H}\bm{1}\bm{1}^{\top}E\Lambda^{k+\frac{1}{2}}\|^{2}_{F}\\
=&\frac{1}{1-w}\|(I-\frac{1}{H}\bm{1}\bm{1}^{\top})E\Lambda^{k+\frac{1}{2}}\|^{2}_{F} \notag\\
=&\frac{1}{1-w}\|(I-\frac{1}{H}\bm{1}\bm{1}^{\top})E(I-\frac{1}{H}\bm{1}\bm{1}^{\top})\Lambda^{k+\frac{1}{2}}\|^{2}_{F}\notag\\
\le&\frac{1}{1-w}\|(I-\frac{1}{H}\bm{1}\bm{1}^{\top})E\|^{2}\|(I-\frac{1}{H}\bm{1}\bm{1}^{\top})\Lambda^{k+\frac{1}{2}}\|^{2}_{F} \notag\\
=&\frac{\kappa}{1-w}\|\Lambda^{k+\frac{1}{2}}-\frac{1}{H}\bm{1}\bm{1}^{\top}\Lambda^{k+\frac{1}{2}}
\|^{2}_{F},  \notag
\end{align}
where the last inequality holds because of Assumption \ref{a3} and the fact that $\|AB\|_{F}^{2}\le \|A\|^{2}\|B\|_{F}^{2}$.

\noindent\textbf{Bounding $T_{8}$:} According to the update of $\bm{\lambda}_{i}^{k+1}$ in Algorithm \ref{alg2} and \eqref{eq_definition1} in Definition \ref{d1}, we have
\begin{align}
\label{eq-proof-lemma4-3}
T_{8}=&\frac{2}{w}\|\Lambda^{k+1}-E\Lambda^{k+\frac{1}{2}}\|^{2}_{F}\\
=&\frac{2}{w}\sum_{i\in\mathcal{H} }\|\bm{\lambda}_{i}^{k+1}-\bar{\bm{\lambda}}_{i}^{k+\frac{1}{2}}\|^{2} \notag\\
=&\frac{2}{w}\sum_{i\in\mathcal{H} }\|AGG_{i}(\bm{\lambda}_{i}^{k+\frac{1}{2}},\{\check{\bm{\lambda}}_{j}^{k+\frac{1}{2}}\}_{j\in \mathcal{N}_{i}})-\bar{\bm{\lambda}}_{i}^{k+\frac{1}{2}}\|^{2} \notag\\
\le&\frac{2}{w}\sum_{i\in \mathcal{H}}\rho^{2}\max_{j\in \mathcal{N}_{i}\cap\mathcal{H}\cup i}\|\bm{\lambda}_{j}^{k+\frac{1}{2}}-\bar{\bm{\lambda}}_{i}^{k+\frac{1}{2}}\|^{2} \notag\\
=& \frac{2\rho^{2}}{w}\sum_{i\in \mathcal{H}}\max_{j\in \mathcal{N}_{i}\cap\mathcal{H}\cup i}\|\bm{\lambda}_{j}^{k+\frac{1}{2}}-\bar{\bm{\lambda}}^{k+\frac{1}{2}}+\bar{\bm{\lambda}}^{k+\frac{1}{2}}-\bar{\bm{\lambda}}_{i}^{k+\frac{1}{2}}\|^{2} \notag\\
\le& \frac{4\rho^{2}}{w}\sum_{i\in \mathcal{H}}[\max_{j\in \mathcal{N}_{i}\cap\mathcal{H}\cup i}\|\bm{\lambda}_{j}^{k+\frac{1}{2}}-\bar{\bm{\lambda}}^{k+\frac{1}{2}}\|^{2}+\|\bar{\bm{\lambda}}^{k+\frac{1}{2}}-\bar{\bm{\lambda}}_{i}^{k+\frac{1}{2}}\|^{2}]\notag\\
\le&\frac{4\rho^{2}}{w}\sum_{i\in \mathcal{H}}[\max_{i\in \mathcal{H}}\|\bm{\lambda}_{i}^{k+\frac{1}{2}}-\bar{\bm{\lambda}}^{k+\frac{1}{2}}\|^{2}+\max_{i\in \mathcal{H}}\|\bar{\bm{\lambda}}^{k+\frac{1}{2}}-\bar{\bm{\lambda}}_{i}^{k+\frac{1}{2}}\|^{2}]\notag\\
=&\frac{8\rho^{2}H}{w}\max_{i\in \mathcal{H}}\|\bm{\lambda}_{i}^{k+\frac{1}{2}}-\bar{\bm{\lambda}}^{k+\frac{1}{2}}\|^{2}\notag\\
\le&\frac{8\rho^{2}H}{w} \|\Lambda^{k+\frac{1}{2}}-\frac{1}{H}\bm{1}\bm{1}^{\top}\Lambda^{k+\frac{1}{2}}
\|^{2}_{F},\notag
\end{align}
where the last inequality holds as $\max_{i\in \mathcal{H}}\|\bm{\lambda}_{i}^{k+\frac{1}{2}}-\bar{\bm{\lambda}}^{k+\frac{1}{2}}\|^{2}\le \|\Lambda^{k+\frac{1}{2}}-\frac{1}{H}\bm{1}\bm{1}^{\top}\Lambda^{k+\frac{1}{2}}
\|^{2}_{F}$.

\noindent\textbf{Bounding $T_{9}$:}
Likewise, according to the update of $\bm{\lambda}_{i}^{k+1}$ in Algorithm \ref{alg2} and \eqref{eq_definition1} in Definition \ref{d1}, we have
\begin{align}
\label{eq-proof-lemma4-4}
T_{9}=&\frac{2}{w}\|\frac{1}{H}\bm{1}\bm{1}^{\top}\Lambda^{k+1}-\frac{1}{H}\bm{1}\bm{1}^{\top}E\Lambda^{k+\frac{1}{2}}\|^{2}_{F}\\
=&\frac{2}{w}\|\frac{1}{H}\bm{1}\bm{1}^{\top}(\Lambda^{k+1}-E\Lambda^{k+\frac{1}{2}})\|^{2}_{F} \notag\\
\le & \frac{2}{w} \|\frac{1}{H}\bm{1}\bm{1}^{\top}\|^{2}_{F}\|\Lambda^{k+1}-E\Lambda^{k+\frac{1}{2}}\|^{2}_{F}\notag\\
=&\frac{2}{w}\|\Lambda^{k+1}-E\Lambda^{k+\frac{1}{2}}\|^{2}_{F}\notag\\
\le &\frac{8\rho^{2}H}{w}\max_{i\in \mathcal{H}}\|\bm{\lambda}_{i}^{k+\frac{1}{2}}-\bar{\bm{\lambda}}^{k+\frac{1}{2}}\|^{2}\notag\\
\le&\frac{8\rho^{2}H}{w} \|\Lambda^{k+\frac{1}{2}}-\frac{1}{H}\bm{1}\bm{1}^{\top}\Lambda^{k+\frac{1}{2}}
\|^{2}_{F}.\notag
\end{align}
To drive the last equality, we use the fact $\|\frac{1}{H}\bm{1}\bm{1}^{\top}\|^{2}_{F}=1$. From the last equality to the last inequality, we use the same technique in deriving \eqref{eq-proof-lemma4-3}.

Therefore, substituting \eqref{eq-proof-lemma4-2}, \eqref{eq-proof-lemma4-3} and \eqref{eq-proof-lemma4-4} into \eqref{eq-proof-lemma4-1} and rearranging the terms, we obtain
\begin{align}
\label{eq-proof-lemma4-5}
&\|\Lambda^{k+1}-\frac{1}{H}\bm{1}\bm{1}^{\top}\Lambda^{k+1}
\|^{2}_{F}\\
\le&(\frac{\kappa}{1-w}+\frac{16\rho^{2}H}{w})\|\Lambda^{k+\frac{1}{2}}-\frac{1}{H}\bm{1}\bm{1}^{\top}\Lambda^{k+\frac{1}{2}}
\|^{2}_{F}. \notag
\end{align}
Substituting \eqref{lemma3-1} in Lemma \ref{l3} into \eqref{eq-proof-lemma4-5} and rearranging the terms, we obtain
\begin{align}
\label{eq-proof-lemma4-6}
&\|\Lambda^{k+1}-\frac{1}{H}\bm{1}\bm{1}^{\top}\Lambda^{k+1}
\|^{2}_{F}\\
\le&(\frac{\kappa}{1-w}+\frac{16\rho^{2}H}{w})(\frac{1}{1-v}+\frac{6(\gamma^{k})^{2}}{v\cdot u_{f}^{2}J^{2}}) \|\Lambda^{k}-\frac{1}{H}\bm{1}\bm{1}^{\top}\Lambda^{k}
\|^{2}_{F}\notag\\
&+(\frac{\kappa}{1-w}+\frac{16\rho^{2}H}{w})\cdot\frac{3(\gamma^{k})^{2}\delta^{2}H^{3}}{v\cdot J^{2}}\notag.
\end{align}

By setting the constant $w=4\rho\sqrt{H}\le 1-\kappa$, $\frac{\kappa}{1-w}\le \kappa+w$ holds. Therefore, we can rewrite \eqref{eq-proof-lemma4-6} as
\begin{align}
\label{eq-proof-lemma4-7}
&\|\Lambda^{k+1}-\frac{1}{H}\bm{1}\bm{1}^{\top}\Lambda^{k+1}
\|^{2}_{F}\\
\le& (\kappa+8\rho\sqrt{H})(\frac{1}{1-v}+\frac{6(\gamma^{k})^{2}}{v\cdot u_{f}^{2}J^{2} })\|\Lambda^{k}-\frac{1}{H}\bm{1}\bm{1}^{\top}\Lambda^{k}
\|^{2}_{F}\notag\\
&+(\kappa+8\rho\sqrt{H})\cdot\frac{3(\gamma^{k})^{2}\delta^{2}H^{3}}{v\cdot J^{2}}\notag\\
=&(1-\epsilon)(\frac{1}{1-v}+\frac{6(\gamma^{k})^{2}}{v\cdot u_{f}^{2}J^{2}})\|\Lambda^{k}-\frac{1}{H}\bm{1}\bm{1}^{\top}\Lambda^{k}
\|^{2}_{F}\notag\\
&+(1-\epsilon)\cdot\frac{3(\gamma^{k})^{2}\delta^{2}H^{3}}{v\cdot J^{2}},\notag
\end{align}
where $\epsilon:=1-\kappa-8\rho\sqrt{H}$. The parameter $\rho$ should satisfy $\rho<\frac{1-\kappa}{8\sqrt{H}}$ to guarantee $\epsilon>0$.

Set $v=\frac{\epsilon}{3}$ and a proper step size $\gamma^{k}$ satisfying $\frac{6(\gamma^{k})^{2}}{v\cdot u_{f}^{2}J^{2}}\le \frac{(2-\epsilon)\epsilon^{2}}{3(3-\epsilon)}=\frac{(\epsilon-v-v\epsilon)v}{1-v}$. Therefore, we have $\frac{1}{1-v}+\frac{6(\gamma^{k})^{2}}{v\cdot u_{f}^{2}J^{2}}\le 1+\epsilon$. In consequence, \eqref{eq-proof-lemma4-7} can be rewritten as
\begin{align}
\label{eq-proof-lemma4-8}
&\|\Lambda^{k+1}-\frac{1}{H}\bm{1}\bm{1}^{\top}\Lambda^{k+1}
\|^{2}_{F}\\
\le& (1-\epsilon^{2})\|\Lambda^{k}-\frac{1}{H}\bm{1}\bm{1}^{\top}\Lambda^{k}
\|^{2}_{F}+\frac{9(\gamma^{k})^{2}\delta^{2}H^{3}}{\epsilon \cdot J^{2}}.\notag
\end{align}
Under the conditions $\rho<\frac{1-\kappa}{8\sqrt{H}}$ and $\epsilon\in(0,1)$, we write \eqref{eq-proof-lemma4-8} recursively to yield
\begin{align}
\label{eq-proof-lemma4-9}
&\|\Lambda^{k+1}-\frac{1}{H}\bm{1}\bm{1}^{\top}\Lambda^{k+1}
\|^{2}_{F}\\
\le& (1-\epsilon^{2})^{k+1}\|\Lambda^{0}-\frac{1}{H}\bm{1}\bm{1}^{\top}\Lambda^{0}
\|^{2}_{F}\notag\\
&+\sum_{k'=0}^{k}(1-\epsilon^{2})^{k-k'}\cdot(\gamma^{k'})^{2}\cdot\frac{9\delta^{2}H^{3}}{\epsilon \cdot J^{2}}.\notag
\end{align}

With the same initialization $\bm{\lambda}_{i}^{0}$ for all honest agents $i\in \mathcal{H}$, we can rewrite \eqref{eq-proof-lemma4-9} as
\begin{align}
\label{eq-proof-lemma4-10}
&\|\Lambda^{k+1}-\frac{1}{H}\bm{1}\bm{1}^{\top}\Lambda^{k+1}
\|^{2}_{F}\\
\le& \sum_{k'=0}^{k}(1-\epsilon^{2})^{k-k'}\cdot(\gamma^{k'})^{2}\cdot\frac{9\delta^{2}H^{3}}{\epsilon \cdot J^{2}}.\notag
\end{align}
To bound $\sum_{k'=0}^{k}(1-\epsilon^{2})^{k-k'}\cdot (\gamma^{k'})^{2}$ in \eqref{eq-proof-lemma4-10}, we define $y^{k}$ as
\begin{align}
y^{k}=\sum_{k'=0}^{k-1}(1-\epsilon^{2})^{k-1-k'}\cdot (\gamma^{k'})^{2},
\end{align}
which satisfies the relation $y^{k+1}=(1-\epsilon^{2})y^{k}+(\gamma^{k})^{2}$. Substituting $\psi_{1}=1-\epsilon^{2}\in (0,1)$, $\psi_{2}=1\ge 0$ and $y^{0}=0\le (\gamma^{0})^{2}$ to Lemma \ref{l4-appendix1}, for integer $k\geq0$ and the step size $\gamma^{k}$ satisfying $1\le \frac{(\gamma^{k})^{2}}{(\gamma^{k+1)^{2}}}\le \frac{2}{1+(1-\epsilon^{2})}$, we have
\begin{align}
\label{eq-proof-lemma4-11}
y^{k+1}&=\sum_{k'=0}^{k}(1-\epsilon^{2})^{k-k'}\cdot (\gamma^{k'})^{2}\\
&\le \frac{2}{1-(1-\epsilon^{2})}(\gamma^{k+1})^{2}= \frac{2}{\epsilon^{2}}(\gamma^{k+1})^{2}.\notag
\end{align}
With \eqref{eq-proof-lemma4-11}, we can rewrite \eqref{eq-proof-lemma4-10} as
\begin{align}
\label{eq-proof-lemma3-11}
\|\Lambda^{k+1}-\frac{1}{H}\bm{1}\bm{1}^{\top}\Lambda^{k+1}
\|^{2}_{F}\le \frac{18(\gamma^{k+1})^{2}\delta^{2}H^{3}}{\epsilon^{3} \cdot J^{2}},
\end{align}
which completes the proof.
\end{lemma}

\begin{lemma}
\label{l4-appendix1}
Suppose that for any integer $k\ge 0$, a sequence $\{\gamma^{k}\}$ satisfies
\begin{align}
\label{l4-appendix1-1}
1\le \frac{(\gamma^{k})^{2}}{(\gamma^{k+1})^{2}}\le \frac{2}{1+\psi_{1}}
\end{align}
for some $\psi_{1}\in(0,1)$, and another sequence $\{y^{k}\}$ satisfies
\begin{align}
\label{l4-appendix1-2}
y^{k+1}\le \psi_{1}y^{k}+\psi_{2}(\gamma^{k})^{2} \quad \mbox{and} \quad y^{0}\le \psi_{2}(\gamma^{0})^{2}
\end{align}
for some $\psi_{1}\in(0,1)$ and $\psi_{2}\ge 0$.
Then, $y^{k}$ is upper-bounded by
\begin{align}
y^{k}\le \frac{2\psi_{2}}{1-\psi_{1}}(\gamma^{k})^{2}.
\label{l4-appendix1-3}
\end{align}

\textbf{Proof}. With the conditions $1\le \frac{2}{1+\psi_{1}}$ and $y^{0}\le \psi_{2}(\gamma^{0})^{2}$, we have $y^{0}\le \frac{2\psi_{2}}{1-\psi_{1}}(\gamma^{0})^{2}$. Therefore, when $k=0$, the proposition $y^{k}\le \frac{2\psi_{2}}{1-\psi_{1}}(\gamma^{k})^{2}$ holds.

Now we prove the conclusion by mathematical induction. Suppose that when $k=k'$, the proposition $y^{k'}\le \frac{2\psi_{2}}{1-\psi_{1}}(\gamma^{k'})^{2}$ holds. We analyze when $k=k'+1$, whether $y^{k'+1}\le \frac{2\psi_{2}}{1-\psi_{1}}(\gamma^{k'+1})^{2}$ holds. Combining $y^{k'+1}\le \psi_{1}y^{k'}+\psi_{2}(\gamma^{k'})^{2}$ and $y^{k'}\le \frac{2\psi_{2}}{1-\psi_{1}}(\gamma^{k'})^{2}$, we obtain $y^{k'+1}\le \frac{(1+\psi_{1})\psi_{2}}{1-\psi_{1}}(\gamma^{k'})^{2}$. Since the step size $\gamma^{k'}$ satisfies $\frac{(\gamma^{k'})^{2}}{(\gamma^{k'+1})^{2}}\le \frac{2}{1+\psi_{1}}$, we have $(\gamma^{k'})^{2}\le \frac{2}{1+\psi_{1}}\cdot(\gamma^{k'+1})^{2}$ and conclude that $y^{k'+1}\le \frac{(1+\psi_{1})\psi_{2}}{1-\psi_{1}}(\gamma^{k'})^{2}\le \frac{2\psi_{2}}{1-\psi_{1}}(\gamma^{k'+1})^{2}$. Hence, when $k=k'+1$, $y^{k'+1}\le \frac{2\psi_{2}}{1-\psi_{1}}(\gamma^{k'+1})^{2}$ holds. This completes the proof.
\end{lemma}

\section{Proof of Theorem \ref{t1}}\label{appendix1}

\subsection{Part $a$ of Theorem \ref{t1}}
According to the update of $\bm{\lambda}_{i}^{k+1}$ in Algorithm \ref{alg1}, we have
\begin{align}
\label{eq_proof_theorem1-b-1}
\|\bar{\bm{\lambda}}^{k+1}-\widetilde{\bm{\lambda}}^{*}\|^{2}=\|\frac{1}{J}\sum_{i\in \mathcal{J}}\sum_{j\in \mathcal{N}_{i}\cup \{i\}}\widetilde{e}_{ij}\bm{\lambda}_{j}^{k+\frac{1}{2}}-\widetilde{\bm{\lambda}}^{*}\|^{2}.
\end{align}
According to Assumption \ref{a3}, $\widetilde{E}$ is doubly stochastic. Therefore, we have $\sum_{i\in \mathcal{J}}\sum_{j\in \mathcal{N}_{i}\cup \{i\}}\widetilde{e}_{ij}\bm{\lambda}_{j}^{k+\frac{1}{2}}=\sum_{i\in \mathcal{J}}\bm{\lambda}_{i}^{k+\frac{1}{2}}$. Combining the update of $\bm{\lambda}_{i}^{k+\frac{1}{2}}$ in Algorithm \ref{alg1} and the fact $\sum_{i\in \mathcal{J}}\sum_{j\in \mathcal{N}_{i}\cup \{i\}}\widetilde{e}_{ij}\bm{\lambda}_{j}^{k+\frac{1}{2}}=\sum_{i\in \mathcal{J}}\bm{\lambda}_{i}^{k+\frac{1}{2}}$, we can rewrite \eqref{eq_proof_theorem1-b-1} as
\begin{align}
\label{eq_proof_theorem1-b-2}
&\|\bar{\bm{\lambda}}^{k+1}-\widetilde{\bm{\lambda}}^{*}\|^{2}\\
=& \|\frac{1}{J}\sum_{i\in \mathcal{J}}[\bm{\lambda}_{i}^{k}-\gamma^{k}\nabla \widetilde{g}_{i}(\bm{\lambda}_{i}^{k})]-\widetilde{\bm{\lambda}}^{*}\|^{2}\notag\\
=&\|\bar{\bm{\lambda}}^{k}-\widetilde{\bm{\lambda}}^{*}-\frac{\gamma^{k}}{J}\sum_{i\in \mathcal{J}}\nabla \widetilde{g}_{i}(\bar{\bm{\lambda}}^{k})\notag\\
&+\frac{\gamma^{k}}{J}\sum_{i\in \mathcal{J}}\nabla \widetilde{g}_{i}(\bar{\bm{\lambda}}^{k})-\frac{\gamma^{k}}{J}\sum_{i\in \mathcal{J}}\nabla \widetilde{g}_{i}(\bm{\lambda}_{i}^{k})\|^{2} \notag\\
=&\|\bar{\bm{\lambda}}^{k}-\widetilde{\bm{\lambda}}^{*}-\frac{\gamma^{k}}{J}\sum_{i\in \mathcal{J}}\nabla \widetilde{g}_{i}(\bar{\bm{\lambda}}^{k})\|^{2}\notag\\
&+\frac{(\gamma^{k})^{2}}{J^{2}}\|\sum_{i\in \mathcal{J}}(\nabla \widetilde{g}_{i}(\bar{\bm{\lambda}}^{k})-\nabla \widetilde{g}_{i}(\bm{\lambda}_{i}^{k}))\|^{2}+2\gamma^{k}\cdot\notag\\
&\left \langle \bar{\bm{\lambda}}^{k}-\widetilde{\bm{\lambda}}^{*}-\frac{\gamma^{k}}{J}\sum_{i\in \mathcal{J}}\nabla \widetilde{g}_{i}(\bar{\bm{\lambda}}^{k}), \frac{1}{J}\sum_{i\in \mathcal{J}}(\nabla \widetilde{g}_{i}(\bar{\bm{\lambda}}^{k})-\nabla \widetilde{g}_{i}(\bm{\lambda}_{i}^{k}))\right \rangle\notag\\
\le&\|\bar{\bm{\lambda}}^{k}-\widetilde{\bm{\lambda}}^{*}-\frac{\gamma^{k}}{J}\sum_{i\in \mathcal{J}}\nabla \widetilde{g}_{i}(\bar{\bm{\lambda}}^{k})\|^{2}\notag\\
&+\frac{(\gamma^{k})^{2}}{J^{2}}\|\sum_{i\in \mathcal{J}}(\nabla \widetilde{g}_{i}(\bar{\bm{\lambda}}^{k})-\nabla \widetilde{g}_{i}(\bm{\lambda}_{i}^{k}))\|^{2}\notag\\
&+ \frac{v^{-1}\gamma^{k}}{J^{2}}\|\sum_{i\in \mathcal{J}}(\nabla \widetilde{g}_{i}(\bar{\bm{\lambda}}^{k})-\nabla \widetilde{g}_{i}(\bm{\lambda}_{i}^{k}))\|^{2}\notag\\
&+v\gamma^{k}\|\bar{\bm{\lambda}}^{k}-\widetilde{\bm{\lambda}}^{*}-\frac{\gamma^{k}}{J}\sum_{i\in \mathcal{J}}\nabla \widetilde{g}_{i}(\bar{\bm{\lambda}}^{k})\|^{2}\notag\\
\le&\underbrace{(1+v\gamma^{k})\|\bar{\bm{\lambda}}^{k}-\widetilde{\bm{\lambda}}^{*}-\frac{\gamma^{k}}{J}\sum_{i\in \mathcal{J}}\nabla \widetilde{g}_{i}(\bar{\bm{\lambda}}^{k})\|^{2}}_{T_{10}}\notag\\
&+\underbrace{\frac{\gamma^{k}(\gamma^{k}+v^{-1})}{J}\sum_{i\in \mathcal{J}}\|\nabla \widetilde{g}_{i}(\bar{\bm{\lambda}}^{k})-\nabla \widetilde{g}_{i}(\bm{\lambda}_{i}^{k})\|^{2}}_{T_{11}},\notag
\end{align}
where $v>0$ is any positive constant. To drive the first inequality, we use $2\bm{a}^{\top}\bm{b}\le v^{-1}\|\bm{a}\|^{2}+v\|\bm{b}\|^{2}$ for any $v>0$. The last inequality holds because $(a_{1}+\cdots+a_{J})^{2}\le J(a_{1}^{2}+\cdots+a_{J}^{2})$. Next, we analyze $T_{10}$ and $T_{11}$ in turn.\\
\noindent\textbf{Bounding $T_{10}$:}
According to $\sum_{i\in \mathcal{J}}\nabla \widetilde{g}_{i}(\widetilde{\bm{\lambda}}^{*})=\bm{0}$, we have
\begin{align}
\label{eq_proof_theorem1-b-T1-1}
& \hspace{-1em} T_{10}= (1+v\gamma^{k})\|\bar{\bm{\lambda}}^{k}-\widetilde{\bm{\lambda}}^{*}-\frac{\gamma^{k}}{J}(\sum_{i\in \mathcal{J}}\nabla \widetilde{g}_{i}(\bar{\bm{\lambda}}^{k})-\sum_{i\in \mathcal{J}}\nabla \widetilde{g}_{i}(\widetilde{\bm{\lambda}}^{*}))\|^{2}\notag
\end{align}
\begin{align}
& \hspace{-1em} = (1+v\gamma^{k})\|\bar{\bm{\lambda}}^{k}-\widetilde{\bm{\lambda}}^{*}\|^{2}\\
& \hspace{-1em} +(1+v\gamma^{k})(\gamma^{k})^{2}\|\frac{1}{J}\sum_{i\in \mathcal{J}}\nabla \widetilde{g}_{i}(\bar{\bm{\lambda}}^{k})-\frac{1}{J}\sum_{i\in \mathcal{J}}\nabla \widetilde{g}_{i}(\widetilde{\bm{\lambda}}^{*})\|^{2}\notag\\
& \hspace{-1em} -2\gamma^{k}(1+v\gamma^{k})\cdot \left \langle  \bar{\bm{\lambda}}^{k}-\widetilde{\bm{\lambda}}^{*},\frac{1}{J}\sum_{i\in \mathcal{J}}\nabla \widetilde{g}_{i}(\bar{\bm{\lambda}}^{k})-\frac{1}{J}\sum_{i\in \mathcal{J}}\nabla \widetilde{g}_{i}(\widetilde{\bm{\lambda}}^{*}) \right \rangle.  \notag
\end{align}

Now we analyze the last term at the right-hand side of \eqref{eq_proof_theorem1-b-T1-1}, $\left \langle  \bar{\bm{\lambda}}^{k}-\widetilde{\bm{\lambda}}^{*},\frac{1}{J}\sum_{i\in \mathcal{J}}\nabla \widetilde{g}_{i}(\bar{\bm{\lambda}}^{k})-\frac{1}{J}\sum_{i\in \mathcal{J}}\nabla \widetilde{g}_{i}(\widetilde{\bm{\lambda}}^{*}) \right \rangle$. By Lemma \ref{l2-appendix1}, $\widetilde{g}_{i}(\cdot)$ is strongly convex with constant $\frac{1}{JL_{f}}$ and smooth with constant $\frac{1}{Ju_{f}}$. According to Lemma 3 in \cite{b-Kun-Yuan-2016}, since $\frac{1}{J}\sum_{i\in \mathcal{J}}\widetilde{g}_{i}(\cdot)$ is $\frac{1}{JL_{f}}$-strongly convex and  $\frac{1}{Ju_{f}}$-smooth, we have
\begin{align}
\label{eq_proof_theorem1-b-T1-2}
&\left \langle  \bar{\bm{\lambda}}^{k}-\widetilde{\bm{\lambda}}^{*},\frac{1}{J}\sum_{i\in \mathcal{J}}\nabla \widetilde{g}_{i}(\bar{\bm{\lambda}}^{k})-\frac{1}{J}\sum_{i\in \mathcal{J}}\nabla \widetilde{g}_{i}(\widetilde{\bm{\lambda}}^{*}) \right \rangle \\
&\ge \widetilde{\alpha}\|\frac{1}{J}\sum_{i\in \mathcal{J}}\nabla \widetilde{g}_{i}(\bar{\bm{\lambda}}^{k})-\frac{1}{J}\sum_{i\in \mathcal{J}}\nabla \widetilde{g}_{i}(\widetilde{\bm{\lambda}}^{*})\|^{2}+\widetilde{\beta}\|\bar{\bm{\lambda}}^{k}-\widetilde{\bm{\lambda}}^{*}\|^{2},\notag
\end{align}
where $\widetilde{\alpha}=\frac{J u_{f}L_{f}}{u_{f}+L_{f}}$ and $\widetilde{\beta}=\frac{1}{J(u_{f}+L_{f})}$.
Substituting \eqref{eq_proof_theorem1-b-T1-2} into \eqref{eq_proof_theorem1-b-T1-1} and rearranging the terms, we have
\begin{align}
\label{eq_proof_theorem1-b-T1}
T_{10}\le& (1+v\gamma^{k})(1-2\gamma^{k}\widetilde{\beta})\|\bar{\bm{\lambda}}^{k}-\widetilde{\bm{\lambda}}^{*}\|^{2}\\
&+(1+v\gamma^{k})((\gamma^{k})^{2}-2\gamma^{k}\widetilde{\alpha})\cdot \notag\\
&\hspace{1em} \|\frac{1}{J}\sum_{i\in \mathcal{J}}\nabla \widetilde{g}_{i}(\bar{\bm{\lambda}}^{k})-\frac{1}{J}\sum_{i\in \mathcal{J}}\nabla \widetilde{g}_{i}(\widetilde{\bm{\lambda}}^{*})\|^{2}\notag\\
\le& (1+v\gamma^{k})(1-2\gamma^{k}\widetilde{\beta})\|\bar{\bm{\lambda}}^{k}-\widetilde{\bm{\lambda}}^{*}\|^{2},\notag
\end{align}
where the last inequality holds with a proper step size $\gamma^{k}$ satisfying $(\gamma^{k})^{2}-2\gamma^{k}\widetilde{\alpha}\le 0$.

\noindent\textbf{Bounding $T_{11}$:} Since $\widetilde{g}_{i}(\cdot)$ is $\frac{1}{Ju_{f}}$-smooth, we obtain
\begin{align}
\label{eq_proof_theorem1-b-T2}
T_{11}\le \frac{\gamma^{k}(\gamma^{k}+v^{-1})}{J^{3}u_{f}^{2}}\sum_{i\in \mathcal{J}}\|\bm{\lambda}_{i}^{k}-\bar{\bm{\lambda}}^{k}\|^{2}.
\end{align}
Substituting \eqref{eq_proof_theorem1-b-T1} and \eqref{eq_proof_theorem1-b-T2} into \eqref{eq_proof_theorem1-b-2} and rearranging the terms, we have
\begin{align}
\label{eq_proof_theorem1-b-3}
&\|\bar{\bm{\lambda}}^{k+1}-\widetilde{\bm{\lambda}}^{*}\|^{2}\\
\le & (1+v\gamma^{k})(1-\gamma^{k}\widetilde{\beta})\|\bar{\bm{\lambda}}^{k}-\widetilde{\bm{\lambda}}^{*}\|^{2}\notag\\
&+\frac{\gamma^{k}(\gamma^{k}+v^{-1})}{J^{3}u_{f}^{2}}\sum_{i\in \mathcal{J}}\|\bm{\lambda}_{i}^{k}-\bar{\bm{\lambda}}^{k}\|^{2}\notag\\
=&(1+v\gamma^{k})(1-\gamma^{k}\widetilde{\beta})\|\bar{\bm{\lambda}}^{k}-\widetilde{\bm{\lambda}}^{*}\|^{2}\notag\\
&+\frac{\gamma^{k}(\gamma^{k}+v^{-1})}{J^{3}u_{f}^{2}}\|\widetilde{\Lambda}^{k}-\frac{1}{J}\widetilde{\bm{1}}\widetilde{\bm{1}}^{\top}\widetilde{\Lambda}^{k}\|_{F}^{2}, \notag
\end{align}
where the matrix $\widetilde{\Lambda}:=[\cdots,\bm{\lambda}_{i},\cdots]\in \mathbb{R}^{J\times D} $ collects $\bm{\lambda}_{i}$ of all agents $i\in \mathcal{J}$. To drive the last equality, we use the fact that $\sum_{i\in \mathcal{J}}\|\bm{\lambda}_{i}^{k}-\bar{\bm{\lambda}}^{k}\|^{2}=\|\widetilde{\Lambda}^{k}-\frac{1}{J}\widetilde{\bm{1}}\widetilde{\bm{1}}^{\top}\widetilde{\Lambda}^{k}\|_{F}^{2}$. Based on Lemma \ref{l3-appendix1}, we can rewrite \eqref{eq_proof_theorem1-b-3} as
\begin{align}
\label{eq_proof_theorem1-b-4}
&\|\bar{\bm{\lambda}}^{k+1}-\widetilde{\bm{\lambda}}^{*}\|^{2}\\
\le &(1+v\gamma^{k})(1-\gamma^{k}\widetilde{\beta})\|\bar{\bm{\lambda}}^{k}-\widetilde{\bm{\lambda}}^{*}\|^{2}\notag\\
&+\frac{18(\gamma^{k})^{3}\widetilde{\delta}^{2}(\gamma^{k}+v^{-1})}{\sigma^{3}u_{f}^{2}J^{2}}. \notag
\end{align}

Setting $v=\frac{\widetilde{\beta}}{2(1-\gamma^{k}\widetilde{\beta})}$, we can rewrite \eqref{eq_proof_theorem1-b-4} as
\begin{align}
\label{eq_proof_theorem1-b-5}
\hspace{-1em}\|\bar{\bm{\lambda}}^{k+1}-\widetilde{\bm{\lambda}}^{*}\|^{2}\le (1-\frac{\gamma^{k}\widetilde{\beta}}{2})\|\bar{\bm{\lambda}}^{k}-\widetilde{\bm{\lambda}}^{*}\|^{2}+\frac{36(\gamma^{k})^{3}\widetilde{\delta}^{2}}{\sigma^{3}u_{f}^{2}J^{2}\widetilde{\beta}}.
\end{align}
We further set a proper decreasing step size $\gamma^{k}=\frac{2}{\widetilde{\beta}(k_{0}+k)}$, where $k_{0}> 1$ is any positive integer. Thus, $1-\frac{\gamma^{k}\widetilde{\beta}}{2}=1-\frac{1}{k_{0}+k}$, and \eqref{eq_proof_theorem1-b-5} can be rewritten as
\begin{align}
\label{eq_proof_theorem1-b-6}
&\|\bar{\bm{\lambda}}^{k+1}-\widetilde{\bm{\lambda}}^{*}\|^{2}\\
\le&(1-\frac{1}{k_{0}+k})\|\bar{\bm{\lambda}}^{k}-\widetilde{\bm{\lambda}}^{*}\|^{2}+\frac{1}{(k_{0}+k)^{3}}\cdot\frac{288\widetilde{\delta}^{2}}{\widetilde{\beta}^{4}\sigma^{3}u_{f}^{2}J^{2}}.\notag
\end{align}
Then we rewrite \eqref{eq_proof_theorem1-b-6} recursively and obtain
\begin{align}
\label{eq_proof_theorem1-b-7}
&\|\bar{\bm{\lambda}}^{k+1}-\widetilde{\bm{\lambda}}^{*}\|^{2}\\
&\le \underbrace{\prod_{k'=0}^{k}(1-\frac{1}{k_{0}+k-k'})}_{T_{12}}\|\bar{\bm{\lambda}}^{0}-\widetilde{\bm{\lambda}}^{*}\|^{2}\notag\\
&+\big[\prod_{k'=0}^{k-1}(1-\frac{1}{k_{0}+k-k'})\frac{1}{(k_{0}+0)^{3}}+\cdots+\notag\\
&+\underbrace{\prod_{k'=0}^{0}\frac{(1-\frac{1}{k_{0}+k-k'})}{(k_{0}+k-1)^{3}}+\frac{1}{(k_{0}+k)^{3}}}_{T_{13}}\big]\cdot\frac{288\widetilde{\delta}^{2}}{\widetilde{\beta}^{4}\sigma^{3}u_{f}^{2}J^{2}}.\notag
\end{align}
Next we analyze $T_{12}$ and $T_{13}$ in turn.

For $T_{12}$, we have
\begin{align}
\label{eq_proof_theorem1-b-T3}
T_{12}=&\prod_{k'=0}^{k}(1-\frac{1}{k_{0}+k-k'})\\
=&\frac{k_{0}+k-1}{k_{0}+k}\cdot\frac{k_{0}+k-2}{k_{0}+k-1}\cdot\cdots\cdot\frac{k_{0}-1}{k_{0}}\notag\\
=&\frac{k_{0}-1}{k_{0}+k}.\notag
\end{align}

For $T_{13}$, we have
\begin{align}
\label{eq_proof_theorem1-b-T4}
T_{13}=&\frac{k_{0}}{k_{0}+k}\cdot\frac{1}{(k_{0}+0)^{3}}+\frac{k_{0}+1}{k_{0}+k}\cdot\frac{1}{(k_{0}+1)^{3}}+\cdots\\
&+\frac{k_{0}+k-1}{k_{0}+k}\cdot\frac{1}{(k_{0}+k-1)^{3}}+\frac{1}{(k_{0}+k)^{3}}\notag\\
=& \frac{1}{k_{0}+k}[\frac{1}{(k_{0})^{2}}+\frac{1}{(k_{0}+1)^{2}}+\cdots\notag\\
&+\frac{1}{(k_{0}+k-1)^{2}}+\frac{1}{(k_{0}+k)^{2}}]\notag\\
\le& \frac{1}{k_{0}+k}\cdot\frac{1}{k_{0}-1}.\notag
\end{align}
To drive the last inequality, we use $\sum_{k'=k_{0}}^{k}\frac{1}{{k'}^{2}}\le \frac{1}{k_{0}-1}$.

Substituting \eqref{eq_proof_theorem1-b-T3} and \eqref{eq_proof_theorem1-b-T4} into \eqref{eq_proof_theorem1-b-7}, we have
\begin{align}
\label{eq_proof_theorem1-b-8}
&\|\bar{\bm{\lambda}}^{k+1}-\widetilde{\bm{\lambda}}^{*}\|^{2}\\
\le &\frac{k_{0}-1}{k_{0}+k}\|\bar{\bm{\lambda}}^{0}-\widetilde{\bm{\lambda}}^{*}\|^{2}+\frac{1}{(k_{0}+k)(k_{0}-1)}\cdot \frac{288\widetilde{\delta}^{2}}{\widetilde{\beta}^{4}\sigma^{3}u_{f}^{2}J^{2}}.\notag
\end{align}
Applying the triangle inequality into \eqref{eq_proof_theorem1-b-8}, we obtain
\begin{align}
\label{eq_proof_theorem1-b-9}
&\|\bar{\bm{\lambda}}^{k+1}-\widetilde{\bm{\lambda}}^{*}\|\\
\le &\sqrt{\frac{k_{0}-1}{k_{0}+k}}\|\bar{\bm{\lambda}}^{0}-\widetilde{\bm{\lambda}}^{*}\|+\sqrt{\frac{1}{(k_{0}+k)(k_{0}-1)}\cdot \frac{288\widetilde{\delta}^{2}}{\widetilde{\beta}^{4}\sigma^{3}u_{f}^{2}J^{2}}}.\notag
\end{align}

According to Lemma \ref{l3-appendix1} and using the step size  $\gamma^{k}=\frac{2}{\widetilde{\beta}(k_{0}+k)}$, we have
\begin{align}
\label{eq_proof_theorem1-b-10}
\sum_{i\in \mathcal{J}}\|\bm{\lambda}_{i}^{k+1}-\bar{\bm{\lambda}}^{k+1}\|^{2}=&\|\widetilde{\Lambda}^{k+1}-\frac{1}{J}\widetilde{\bm{1}}\widetilde{\bm{1}}^{\top}\widetilde{\Lambda}^{k+1}\|^{2}_{F}\\
\le& \frac{72\widetilde{\delta}^{2}J}{\sigma^{3}\widetilde{\beta}^{2}(k_{0}+k+1)^{2}}.\notag
\end{align}
By $(a_{1}+\cdots+a_{J})^{2}\le J(a_{1}^{2}+\cdots+a_{J}^{2})$ and \eqref{eq_proof_theorem1-b-10}, we have
\begin{align}
\label{eq_proof_theorem1-b-11}
\sum_{i\in \mathcal{J}}\|\bm{\lambda}_{i}^{k+1}-\bar{\bm{\lambda}}^{k+1}\|\le& \sqrt{J}\sqrt{\sum_{i\in \mathcal{J}}\|\bm{\lambda}_{i}^{k+1}-\bar{\bm{\lambda}}^{k+1}\|^{2}}\\
\le& \sqrt{J}\sqrt{\frac{72\widetilde{\delta}^{2}J}{\sigma^{3}\widetilde{\beta}^{2}(k_{0}+k+1)^{2}}}.\notag
\end{align}

Combining \eqref{eq_proof_theorem1-b-9} and \eqref{eq_proof_theorem1-b-11} yields
\begin{align}
\label{eq_proof_theorem1-b-12}
&\sum_{i\in \mathcal{J}}\|\bm{\lambda}_{i}^{k+1}-\widetilde{\bm{\lambda}}^{*}\|\\
\le& \sum_{i\in \mathcal{J}}\|\bm{\lambda}_{i}^{k+1}-\bar{\bm{\lambda}}^{k+1}\|+J\cdot \|\bar{\bm{\lambda}}^{k+1}-\widetilde{\bm{\lambda}}^{*}\|\notag\\
\le& J\cdot\sqrt{\frac{k_{0}-1}{k_{0}+k}}\|\bar{\bm{\lambda}}^{0}-\widetilde{\bm{\lambda}}^{*}\|+\sqrt{\frac{1}{(k_{0}+k)(k_{0}-1)} \cdot\frac{288\widetilde{\delta}^{2}}{\widetilde{\beta}^{4}\sigma^{3}u_{f}^{2}}}\notag\\
&+\sqrt{J}\sqrt{\frac{72\widetilde{\delta}^{2}J}{\sigma^{3}\widetilde{\beta}^{2}(k_{0}+k+1)^{2}}}.\notag
\end{align}
Taking $k \to +\infty$, we obtain
\begin{align}
\label{eq_proof_theorem1-b-13}
\lim_{k\to +\infty}\sum_{i\in \mathcal{J}}\|\bm{\lambda}_{i}^{k+1}-\widetilde{\bm{\lambda}}^{*}\|
=0.
\end{align}
And the asymptotic convergence speed of the dual variables in the attack-free algorithm is $O(\frac{1}{\sqrt{k}})$.

\subsection{Part $b$ of Theorem \ref{t1}}
Since $\widetilde{\bm{\Theta}}^{*}$ is the optimal solution of the primal problem \eqref{eq_primal-problem}, we have $\frac{1}{J}\sum_{i\in \mathcal{J}}\widetilde{\bm{\theta}}_{i}^{*} = \bm{s}$ and $\widetilde{\bm{\Theta}}^{*}\in \widetilde{C}$. According to \eqref{eq_Lagrangian-primal-problem}, for any dual variable $\widetilde{\bm{\lambda}}$ we have
\begin{align}
\label{eq_proof_theorem1-a-1}
\widetilde{\mathcal{L}}\left ( \widetilde{\bm{\Theta}}^{*} ;  \widetilde{\bm{\lambda}} \right )=&\frac{1}{J}\sum_{i\in \mathcal{J}}f_{i}\left ( \widetilde{\bm{\theta}}_{i}^{*} \right )+ \left \langle \widetilde{\bm{\lambda} }, \frac{1}{J}\sum_{i\in \mathcal{J}}\widetilde{\bm{\theta}}_{i}^{*}-\bm{s} \right \rangle\\
=& \widetilde{f}(\widetilde{\bm{\Theta}}^{*}).\notag
\end{align}
By Assumption \ref{a2}, the duality gap is zero. Therefore, for any $\widetilde{\bm{\lambda}}$ we obtain $\widetilde{g}(\widetilde{\bm{\lambda}}^{*})=-\widetilde{f}(\widetilde{\bm{\Theta}}^{*})=-\widetilde{\mathcal{L}}\left ( \widetilde{\bm{\Theta}}^{*} ;  \widetilde{\bm{\lambda}} \right )$. Now we introduce a vector $_{}^{\dagger}\widetilde{\bm{\Theta}}^{k+1}:=[ _{}^{\dagger}\bm{\theta}_{1}^{k+1}; \cdots;   _{}^{\dagger}\bm{\theta}_{J}^{k+1}]$, where $_{}^{\dagger}\bm{\theta}_{i}^{k+1}=\arg\min_{\bm{\theta}_{i}\in C_{i}} \{\bm{\theta}_{i}^{\top}\bar{\bm{\lambda}}^{k+1}+f_{i}(\bm{\theta}_{i})\}$ and $\bar{\bm{\lambda}}^{k+1}:=\frac{1}{J}\sum_{i\in \mathcal{J}}\bm{\lambda}_{i}^{k+1}$. Therefore, we have
\begin{align}
\label{eq_proof_theorem1-a-2}
&\widetilde{g}(\bar{\bm{\lambda}}^{k+1})-\widetilde{g}(\widetilde{\bm{\lambda}}^{*})\\
=&-\inf_{\widetilde{\bm{\Theta}}\in \widetilde{C}}\widetilde{\mathcal{L}}\left ( \widetilde{\bm{\Theta}} ;  \bar{\bm{\lambda}}^{k+1} \right )+\widetilde{\mathcal{L}}\left ( \widetilde{\bm{\Theta}}^{*} ;  \bar{\bm{\lambda}}^{k+1} \right )\notag\\
=&-\widetilde{\mathcal{L}}\left ( ^{\dagger}\widetilde{\bm{\Theta}}^{k+1} ;  \bar{\bm{\lambda}}^{k+1} \right )+\widetilde{\mathcal{L}}\left ( \widetilde{\bm{\Theta}}^{*} ;  \bar{\bm{\lambda}}^{k+1} \right ).\notag
\end{align}

Assumption \ref{a1} shows that the local cost function $f_{i}(\cdot)$ is $u_{f}$-strongly convex. Further using the definition of $\widetilde{\mathcal{L}}\left ( \widetilde{\bm{\Theta}} ;  \widetilde{\bm{\lambda}} \right )$ in \eqref{eq_Lagrangian-primal-problem}, we know that $\widetilde{\mathcal{L}}\left ( \widetilde{\bm{\Theta}} ;  \widetilde{\bm{\lambda}} \right )$ is $u_{f}$-strongly convex with respect to $\widetilde{\bm{\Theta}}$. Therefore, we have
\begin{align}
\label{eq_proof_theorem1-a-3}
&\widetilde{\mathcal{L}}\left ( \widetilde{\bm{\Theta}}^{*} ;  \bar{\bm{\lambda}}^{k+1} \right )-\widetilde{\mathcal{L}}\left ( ^{\dagger}\widetilde{\bm{\Theta}}^{k+1} ;  \bar{\bm{\lambda}}^{k+1} \right )\\
\ge& \nabla \widetilde{\mathcal{L}}\left ( ^{\dagger}\widetilde{\bm{\Theta}}^{k+1} ;  \bar{\bm{\lambda}}^{k+1} \right )^{\top} (\widetilde{\bm{\Theta}}^{*}-\ ^{\dagger}\widetilde{\bm{\Theta}}^{k+1})+\frac{u_{f}}{2}\|^{\dagger}\widetilde{\bm{\Theta}}^{k+1}-\widetilde{\bm{\Theta}}^{*}\|^{2}.\notag
\end{align}
Combining \eqref{eq_proof_theorem1-a-2} and \eqref{eq_proof_theorem1-a-3}, we obtain
\begin{align}
\label{eq_proof_theorem1-a-4}
&\widetilde{g}(\bar{\bm{\lambda}}^{k+1})-\widetilde{g}(\widetilde{\bm{\lambda}}^{*})\\
\ge& \nabla \widetilde{\mathcal{L}}\left ( ^{\dagger}\widetilde{\bm{\Theta}}^{k+1} ;  \bar{\bm{\lambda}}^{k+1} \right )^{\top}(\widetilde{\bm{\Theta}}^{*}-\ ^{\dagger}\widetilde{\bm{\Theta}}^{k+1})+\frac{u_{f}}{2}\|^{\dagger}\widetilde{\bm{\Theta}}^{k+1}-\widetilde{\bm{\Theta}}^{*}\|^{2}\notag\\
\ge&\frac{u_{f}}{2}\|^{\dagger}\widetilde{\bm{\Theta}}^{k+1}-\widetilde{\bm{\Theta}}^{*}\|^{2}. \notag
\end{align}
To drive the last inequality, we use the optimality condition of $_{}^{\dagger}\bm{\theta}_{i}^{k+1}=\arg\min_{\bm{\theta}_{i}\in C_{i}} \{\bm{\theta}_{i}^{\top}\bar{\bm{\lambda}}^{k+1}+f_{i}(\bm{\theta}_{i})\}$ \cite[Proposition 2.1.2]{b-Dimitri-P.-Bertsekas-1999}.

According to Lemma \ref{l2-appendix1}, $\widetilde{g}_{i}(\widetilde{\bm{\lambda}})$ is smooth with constant $\frac{1}{Ju_{f}}$. Therefore, $\widetilde{g}(\widetilde{\bm{\lambda }}):=\sum_{i\in\mathcal{J}}\widetilde{g}_{i}(\widetilde{\bm{\lambda }})$ is smooth with constant $\frac{1}{u_{f}}$. This fact leads to
\begin{align}
\label{eq_proof_theorem1-a-5}
&\widetilde{g}(\bar{\bm{\lambda}}^{k+1})-\widetilde{g}(\widetilde{\bm{\lambda}}^{*})\\
\le& \nabla \widetilde{g}^{\top}(\widetilde{\bm{\lambda}}^{*})(\bar{\bm{\lambda}}^{k+1}-\widetilde{\bm{\lambda}}^{*})+\frac{1}{2u_{f}}\|\bar{\bm{\lambda}}^{k+1}-\widetilde{\bm{\lambda}}^{*}\|^{2}\notag\\
=& \frac{1}{2u_{f}}\|\bar{\bm{\lambda}}^{k+1}-\widetilde{\bm{\lambda}}^{*}\|^{2}.\notag
\end{align}
To drive the last equality, we use the fact that $\nabla \widetilde{g}(\widetilde{\bm{\lambda}}^{*})=\bm{0}$.
Combining \eqref{eq_proof_theorem1-a-4} and \eqref{eq_proof_theorem1-a-5}, we have
\begin{align}
\label{eq_proof_theorem1-a-6}
\|^{\dagger}\widetilde{\bm{\Theta}}^{k+1}-\widetilde{\bm{\Theta}}^{*}\|^{2} \le \frac{1}{u_{f}^{2}}\|\bar{\bm{\lambda}}^{k+1}-\widetilde{\bm{\lambda}}^{*}\|^{2}.
\end{align}
Combining \eqref{eq_proof_theorem1-b-9} and \eqref{eq_proof_theorem1-a-6}, we obtain
\begin{align}
\label{eq_proof_theorem1-a-7}
&\|^{\dagger}\widetilde{\bm{\Theta}}^{k+1}-\widetilde{\bm{\Theta}}^{*}\|\\ \le& \frac{1}{u_{f}^{2}}\cdot \big[\sqrt{\frac{k_{0}-1}{k_{0}+k}}\|\bar{\bm{\lambda}}^{0}-\widetilde{\bm{\lambda}}^{*}\|\notag\\
&+\sqrt{\frac{1}{(k_{0}+k)(k_{0}-1)}\cdot \frac{288\widetilde{\delta}^{2}(1-\sigma )}{\widetilde{\beta}^{4}\sigma^{3}u_{f}^{2}J^{2}}} \big].\notag
\end{align}
Taking $k\to +\infty$, we obtain
\begin{align}
\label{eq_proof_theorem1-a-8}
\lim_{k\to +\infty}\|^{\dagger}\widetilde{\bm{\Theta}}^{k+1}-\widetilde{\bm{\Theta}}^{*}\|=0.
\end{align}

According to Assumption \ref{a1}, $f_{i}(\bm{\theta}_{i})$ is $u_{f}$-strongly convex. By the conjugate correspondence theorem in \cite{b-Amir-Beck-2017}, we know that the conjugate function $\widetilde{F}_{i}^{*}(\widetilde{\bm{\lambda}})=\max_{\bm{\theta}_{i}\in C_{i}}\{\widetilde{\bm{\lambda}}^{\top}\bm{\theta}_{i}-f_{i}(\bm{\theta}_{i})\}$ is $\frac{1}{u_{f}}$-smooth. Therefore, the gradient $\nabla \widetilde{F}^{*}_{i}(-\widetilde{\bm{\lambda}})=\arg\min_{\bm{\theta}_{i}\in C_{i}}\{\widetilde{\bm{\lambda}}^{\top}\bm{\theta}_{i}+f_{i}(\bm{\theta}_{i})\}$ is $\frac{1}{u_{f}}$-Lipschitz continuous. According to the definition of Lipschitz continuity, we have
\begin{align}
\label{eq_proof_theorem1-a-lipschitz-continuous}
\|\nabla \widetilde{F}^{*}_{i}(\bm{\lambda}_{i}^{k+1})-\nabla \widetilde{F}^{*}_{i}(\bar{\bm{\lambda}}^{k+1})\| \le \frac{1}{u_{f}}\|\bm{\lambda}_{i}^{k+1}-\bar{\bm{\lambda}}^{k+1}\|.
\end{align}
Based on $_{}^{\dagger}\bm{\theta}_{i}^{k+1}=\arg\min_{\bm{\theta}_{i}\in C_{i}} \{\bm{\theta}_{i}^{\top}\bar{\bm{\lambda}}^{k+1}+f_{i}(\bm{\theta}_{i})\}$ and $\bm{\theta}_{i}^{k+1}=\arg\min_{\bm{\theta}_{i}\in C_{i}} \{\bm{\theta}_{i}^{\top}\bm{\lambda}_{i}^{k+1}+f_{i}(\bm{\theta}_{i})\}$, we obtain $\bm{\theta}_{i}^{k+1}=\nabla \widetilde{F}^{*}_{i}(\bm{\lambda}_{i}^{k+1})$ and $_{}^{\dagger}\bm{\theta}_{i}^{k+1}=\nabla \widetilde{F}^{*}_{i}(\bar{\bm{\lambda}}^{k+1})$. Substituting $\bm{\theta}_{i}^{k+1}=\nabla \widetilde{F}^{*}_{i}(\bm{\lambda}_{i}^{k+1})$ and $_{}^{\dagger}\bm{\theta}_{i}^{k+1}=\nabla \widetilde{F}^{*}_{i}(\bar{\bm{\lambda}}^{k+1})$ into \eqref{eq_proof_theorem1-a-lipschitz-continuous}, we have
\begin{align}
\label{eq_proof_theorem1-a-10}
\|\bm{\theta}_{i}^{k+1}-\ _{}^{\dagger}\bm{\theta}_{i}^{k+1}\| \le \frac{1}{u_{f}}\|\bm{\lambda}_{i}^{k+1}-\bar{\bm{\lambda}}^{k+1}\|.
\end{align}
Combining \eqref{eq_proof_theorem1-a-10}, $^{\dagger}\widetilde{\bm{\Theta}}^{k+1}:=[_{}^{\dagger}\bm{\theta}_{1}^{k+1}; \cdots;  _{}^{\dagger}\bm{\theta}_{J}^{k+1}]$ and $\widetilde{\bm{\Theta}}^{k+1}:=[\bm{\theta}_{1}^{k+1}; \cdots; \bm{\theta}_{J}^{k+1}]$, we obtain

\begin{align}
\label{eq_proof_theorem1-a-11}
\|\widetilde{\bm{\Theta}}^{k+1}-\ ^{\dagger}\widetilde{\bm{\Theta}}^{k+1}\| \le \frac{1}{u_{f}}\sum_{i\in \mathcal{J}}\|\bm{\lambda}_{i}^{k+1}-\bar{\bm{\lambda}}^{k+1}\|.
\end{align}
Combining \eqref{eq_proof_theorem1-b-11} and \eqref{eq_proof_theorem1-a-11}, we have
\begin{align}
\label{eq_proof_theorem1-a-12}
\|\widetilde{\bm{\Theta}}^{k+1}-\ ^{\dagger}\widetilde{\bm{\Theta}}^{k+1}\| \le \frac{\sqrt{J}}{u_{f}}\cdot\sqrt{\frac{72\widetilde{\delta}^{2}J}{\sigma^{3}\widetilde{\beta}^{2}(k_{0}+k+1)^{2}}}.
\end{align}
Taking $k\to +\infty$, we obtain
\begin{align}
\label{eq_proof_theorem1-a-13}
\lim_{k\to +\infty}\|\widetilde{\bm{\Theta}}^{k+1}-\ ^{\dagger}\widetilde{\bm{\Theta}}^{k+1}\| =0.
\end{align}

Combining \eqref{eq_proof_theorem1-a-8} and \eqref{eq_proof_theorem1-a-13} yields
\begin{align}
\label{eq_proof_theorem1-a-9}
&\lim_{k\to +\infty}\|\widetilde{\bm{\Theta}}^{k+1}-\widetilde{\bm{\Theta}}^{*}\|\\
\leq &\lim_{k\to +\infty}\|\widetilde{\bm{\Theta}}^{k+1}-\ ^{\dagger}\widetilde{\bm{\Theta}}^{k+1}\|+\lim_{k\to +\infty}\|^{\dagger}\widetilde{\bm{\Theta}}^{k+1}-\widetilde{\bm{\Theta}}^{*}\|=0\notag.
\end{align}
And the asymptotic convergence speed of the primal variables in the attack-free algorithm is also $O(\frac{1}{\sqrt{k}})$.

We conclude by summarizing the conditions on the step size $\gamma^{k}$ in Theorem \ref{t1}. It must satisfy $(\gamma^{k})^{2}-2\gamma^{k}\widetilde{\alpha}\le 0$, $\frac{6(\gamma^{k})^{2}}{ u_{f}^{2}J^{2}}\le\frac{(2-\sigma)\sigma^{2}}{3(3-\sigma)}$ and $1\le \frac{(\gamma^{k})^{2}}{(\gamma^{k+1})^2}\le \frac{2}{1+(1-\sigma^{2})}$. The step size $\gamma^{k}=\frac{2}{\widetilde{\beta}(k_{0}+k)}$ with $k_{0}\ge \max\{\frac{1}{\widetilde{\alpha}\widetilde{\beta}},\sqrt{\frac{72(3-\sigma)}{(2-\sigma)\sigma^{2}u_{f}^{2}J^{2}\widetilde{\beta}^{2}}},\frac{1}{\sqrt{\frac{2}{1+(1-\sigma^{2})}}-1}\}$ satisfies these conditions.

\subsection{Supporting Lemmas}
\begin{lemma}
\label{l1-appendix1}
Under Assumption \ref{a1}, for any $\widetilde{\bm{\lambda}} \in \mathbb{R}^D$, the maximum distance between the local dual gradients and their average, denoted by  $\max_{i\in \mathcal{J}}\|\nabla \widetilde{g}_{i}(\widetilde{\bm{\lambda}}) - \frac{1}{J} \sum_{i\in \mathcal{J}} \nabla \widetilde{g}_{i}(\widetilde{\bm{\lambda}}) \|^{2}$, is bounded by some positive constant $\widetilde{\delta}^{2}$.

\textbf{Proof}. Recalling the definition of the local dual gradient $\nabla \widetilde{g}_{i}(\widetilde{\bm{\lambda}})=-\frac{1}{J}\arg\min_{\bm{\theta}_{i}\in C_{i}} \{\widetilde{\bm{\lambda}}^{\top}\bm{\theta}_{i}+f_{i}(\bm{\theta}_{i})\}+\frac{1}{J}\bm{s}$, we have
\begin{align}
&\max_{i\in \mathcal{J}}\|\nabla \widetilde{g}_{i}(\widetilde{\bm{\lambda}}) - \frac{1}{J} \sum_{i\in \mathcal{J}} \nabla \widetilde{g}_{i}(\widetilde{\bm{\lambda}}) \|^{2}\\
=&\max_{i\in \mathcal{J}}\|-\frac{1}{J}\arg\underset{\bm{\theta}_{i}\in C_{i}}{\min}\{\widetilde{\bm{\lambda}}^{\top}\bm{\theta}_{i}+f_{i}(\bm{\theta}_{i})\}+\frac{1}{J}\bm{s}\notag\\
&\hspace{2.5em}+\frac{1}{J^{2}}\sum_{i\in \mathcal{J}}\arg\underset{\bm{\theta}_{i}\in C_{i}}{\min}\{\widetilde{\bm{\lambda}}^{\top}\bm{\theta}_{i}+f_{i}(\bm{\theta}_{i})\}-\frac{1}{J}\bm{s}\|^{2}\notag\\
=&\max_{i\in \mathcal{J}}\|\frac{1}{J}(\frac{1}{J}\sum_{i\in \mathcal{J}}\arg\underset{\bm{\theta}_{i}\in C_{i}}{\min}\{\widetilde{\bm{\lambda}}^{\top}\bm{\theta}_{i}+f_{i}(\bm{\theta}_{i})\}\notag\\
&\hspace{2.5em}-\arg\underset{\bm{\theta}_{i}\in C_{i}}{\min}\{\widetilde{\bm{\lambda}}^{\top}\bm{\theta}_{i}+f_{i}(\bm{\theta}_{i})\})\|^{2}.\notag
\end{align}
According to Assumption \ref{a1}, the local constraint sets $C_{i}$ are bounded by hypothesis, and we know that $\max_{i\in \mathcal{J}}\|\nabla g_{i}(\widetilde{\bm{\lambda}}) - \frac{1}{J} \sum_{i\in \mathcal{J}} \nabla g_{i}(\widetilde{\bm{\lambda}}) \|^{2}$ is also bounded by some positive constant, which we denote as $\widetilde{\delta}^{2}$.
\end{lemma}

\begin{lemma}
\label{l2-appendix1} Under Assumption \ref{a1}, for any agent $i$, the local dual function $\widetilde{g}_{i}(\widetilde{\bm{\lambda}})$ is strongly convex with constant $\frac{1}{JL_{f}}$ and smooth with parameter $\frac{1}{Ju_{f}}$.

\textbf{Proof}.
According to Assumption \ref{a1}, $f_{i}(\cdot)$ is $u_{f}$-strongly convex and $L_{f}$-smooth. By the conjugate correspondence theorem \cite{b-Amir-Beck-2017}, its conjugate function $\widetilde{F}_{i}^{*}(\widetilde{\bm{\lambda}})=\max_{\bm{\theta}_{i}\in C_{i}}\{\widetilde{\bm{\lambda}}^{\top}\bm{\theta}_{i}-f_{i}(\bm{\theta}_{i})\}$ is $\frac{1}{L_{f}}$-strongly convex and $\frac{1}{u_{f}}$-smooth. By the definition of $\widetilde{g}_{i}(\widetilde{\bm{\lambda}})=\frac{1}{J}\widetilde{F}_{i}^{*}(-\widetilde{\bm{\lambda }})+\frac{1}{J}\widetilde{\bm{\lambda }}^{\top}\bm{s}$, we know $\widetilde{g}_{i}(\widetilde{\bm{\lambda}})$ is $\frac{1}{JL_{f}}$-strongly convex and $\frac{1}{Ju_{f}}$-smooth.
\end{lemma}

\begin{lemma}
\label{l3-appendix1}
Define a matrix $\widetilde{\Lambda}^{k+1}:=[\cdots,\bm{\lambda}_{i}^{k+1},\cdots]\in \mathbb{R}^{J\times D} $ that collects the dual variables $\bm{\lambda}_{i}^{k+1}$ of all agents $i\in \mathcal{J}$ generated by Algorithm \ref{alg1}. Under Assumptions \ref{a1} and \ref{a3}, we have
\begin{align}
\label{lemma3-appendix1}
\|\widetilde{\Lambda}^{k+1}-\frac{1}{J}\widetilde{\bm{1}}\widetilde{\bm{1}}^{\top}\widetilde{\Lambda}^{k+1}
\|^{2}_{F}\le \frac{18(\gamma^{k+1})^{2}\widetilde{\delta}^{2}J}{\sigma^{3}},
\end{align}
where $\sigma=1- \widetilde{\kappa} \in (0,1)$.

\textbf{Proof}.
Define $\nabla \widetilde{g}(\widetilde{\Lambda})=[\cdots,\nabla \widetilde{g}_{i}(\bm{\lambda}_{i}),\cdots]\in \mathbb{R}^{J \times D}$ to collect the dual gradients $\nabla\widetilde{g}_{i}(\bm{\lambda}_{i})$ of all agents $i\in \mathcal{J}$. With these notations, we can rewrite the updates of  $\bm{\lambda}_{i}^{k+1}$ and $\bm{\lambda}_{i}^{k+\frac{1}{2}}$ in Algorithm \ref{alg1} in compact forms of
\begin{align}
\label{lemma3-appendix1-proof-1}
\widetilde{\Lambda}^{k+\frac{1}{2}}=\widetilde{\Lambda}^{k}-\gamma^{k}\nabla \widetilde{g}(\widetilde{\Lambda}^{k}),
\end{align}
\begin{align}
\label{lemma3-appendix1-proof-2}
\widetilde{\Lambda}^{k+1}=\widetilde{E}\widetilde{\Lambda}^{k+\frac{1}{2}}.
\end{align}
Combining \eqref{lemma3-appendix1-proof-1} and \eqref{lemma3-appendix1-proof-2}, and also using the fact that $\widetilde{E}$ is doubly stochastic by Assumption \ref{a3}, we have
\begin{align}
\label{lemma3-appendix1-proof-3}
&\|\widetilde{\Lambda}^{k+1}-\frac{1}{J}\widetilde{\bm{1}}\widetilde{\bm{1}}^{\top}\widetilde{\Lambda}^{k+1}
\|^{2}_{F}\\
=&\|\widetilde{E}(\widetilde{\Lambda}^{k}-\gamma^{k}\nabla \widetilde{g}(\widetilde{\Lambda}^{k}))-\frac{1}{J}\widetilde{\bm{1}}\widetilde{\bm{1}}^{\top}\widetilde{E}(\widetilde{\Lambda}^{k}-\gamma^{k}\nabla \widetilde{g}(\widetilde{\Lambda}^{k}))\|_{F}^{2}\notag\\
=&\|\widetilde{E}\widetilde{\Lambda}^{k}-\frac{1}{J}\widetilde{\bm{1}}\widetilde{\bm{1}}^{\top}\widetilde{\Lambda}^{k}-\widetilde{E}\gamma^{k}\nabla \widetilde{g}(\widetilde{\Lambda}^{k})+\frac{1}{J}\widetilde{\bm{1}}\widetilde{\bm{1}}^{\top}\gamma^{k}\nabla \widetilde{g}(\widetilde{\Lambda}^{k})\|_{F}^{2}\notag\\
\le&\frac{1}{1-v}\|\widetilde{E}\widetilde{\Lambda}^{k}-\frac{1}{J}\widetilde{\bm{1}}\widetilde{\bm{1}}^{\top}\widetilde{\Lambda}^{k}\|_{F}^{2}\notag\\
&+\frac{1}{v}\|\widetilde{E}\gamma^{k}\nabla \widetilde{g}(\widetilde{\Lambda}^{k})-\frac{1}{J}\widetilde{\bm{1}}\widetilde{\bm{1}}^{\top}\gamma^{k}\nabla \widetilde{g}(\widetilde{\Lambda}^{k})\|_{F}^{2} \notag\\
=&\frac{1}{1-v}\|(\widetilde{E}-\frac{1}{J}\widetilde{\bm{1}}\widetilde{\bm{1}}^{\top})(\widetilde{\Lambda}^{k}-\frac{1}{J}\widetilde{\bm{1}}\widetilde{\bm{1}}^{\top}\widetilde{\Lambda}^{k})\|_{F}^{2}\notag\\
&+\frac{(\gamma^{k})^{2}}{v}\|(\widetilde{E}-\frac{1}{J}\widetilde{\bm{1}}\widetilde{\bm{1}}^{\top})(\nabla \widetilde{g}(\widetilde{\Lambda}^{k})-\frac{1}{J}\widetilde{\bm{1}}\widetilde{\bm{1}}^{\top}\nabla \widetilde{g}(\widetilde{\Lambda}^{k}))\|_{F}^{2} \notag\\
\le& \frac{1}{1-v}\|\widetilde{E}-\frac{1}{J}\widetilde{\bm{1}}\widetilde{\bm{1}}^{\top}\|^{2}\|\widetilde{\Lambda}^{k}-\frac{1}{J}\widetilde{\bm{1}}\widetilde{\bm{1}}^{\top}\widetilde{\Lambda}^{k}\|_{F}^{2}\notag\\
&+\frac{(\gamma^{k})^{2}}{v}\|\widetilde{E}-\frac{1}{J}\widetilde{\bm{1}}\widetilde{\bm{1}}^{\top}\|^{2}\|\nabla \widetilde{g}(\widetilde{\Lambda}^{k})-\frac{1}{J}\widetilde{\bm{1}}\widetilde{\bm{1}}^{\top}\nabla \widetilde{g}(\widetilde{\Lambda}^{k})\|_{F}^{2},\notag
\end{align}
where $v\in (0,1)$ is any positive constant. To drive the last inequality, we use the fact that $\|AB\|_{F}^{2}\le \|A\|^{2}\|B\|_{F}^{2}$.
By Assumption \ref{a3}, $\widetilde{\kappa} := \|\widetilde{E}-\frac{1}{J}\widetilde{\bm{1}}\widetilde{\bm{1}}^{\top}\|^{2} < 1$. Thus, we have
\begin{align}
\label{lemma3-appendix1-proof-4}
&\|\widetilde{\Lambda}^{k+1}-\frac{1}{J}\widetilde{\bm{1}}\widetilde{\bm{1}}^{\top}\widetilde{\Lambda}^{k+1}
\|^{2}_{F}\\
\le& \frac{\widetilde{\kappa}}{1-v}\|\widetilde{\Lambda}^{k}-\frac{1}{J}\widetilde{\bm{1}}\widetilde{\bm{1}}^{\top}\widetilde{\Lambda}^{k}\|_{F}^{2}\notag\\
&+\frac{(\gamma^{k})^{2}\widetilde{\kappa}}{v}\|\nabla \widetilde{g}(\widetilde{\Lambda}^{k})-\frac{1}{J}\widetilde{\bm{1}}\widetilde{\bm{1}}^{\top}\nabla \widetilde{g}(\widetilde{\Lambda}^{k})\|_{F}^{2}.\notag
\end{align}

We bound the term $\frac{(\gamma^{k})^{2}\widetilde{\kappa}}{v}\|\nabla \widetilde{g}(\widetilde{\Lambda}^{k})-\frac{1}{J}\widetilde{\bm{1}}\widetilde{\bm{1}}^{\top}\nabla \widetilde{g}(\widetilde{\Lambda}^{k})\|_{F}^{2}$ at the right-hand side of \eqref{lemma3-appendix1-proof-4} as
\begin{align}
\label{lemma3-appendix1-proof-5}
&\frac{(\gamma^{k})^{2}\widetilde{\kappa}}{v}\|\nabla \widetilde{g}(\widetilde{\Lambda}^{k})-\frac{1}{J}\widetilde{\bm{1}}\widetilde{\bm{1}}^{\top}\nabla \widetilde{g}(\widetilde{\Lambda}^{k})\|_{F}^{2}
\end{align}
\begin{align}
=&\frac{(\gamma^{k})^{2}\widetilde{\kappa}}{v}\sum_{i\in \mathcal{J}}\|\nabla \widetilde{g}_{i}(\bm{\lambda}_{i}^{k})-\frac{1}{J}\sum_{i\in \mathcal{J}}\nabla \widetilde{g}_{i}(\bm{\lambda}_{i}^{k})\|^{2}\notag\\
=&\frac{(\gamma^{k})^{2}\widetilde{\kappa}}{v}\sum_{i\in \mathcal{J}}\|\nabla \widetilde{g}_{i}(\bm{\lambda}_{i}^{k})-\nabla \widetilde{g}_{i}(\bar{\bm{\lambda}}^{k})+\nabla \widetilde{g}_{i}(\bar{\bm{\lambda}}^{k})\notag\\
&-\frac{1}{J}\sum_{i\in \mathcal{J}}\nabla \widetilde{g}_{i}(\bar{\bm{\lambda}}^{k})+\frac{1}{J}\sum_{i\in \mathcal{J}}\nabla \widetilde{g}_{i}(\bar{\bm{\lambda}}^{k})
-\frac{1}{J}\sum_{i\in \mathcal{J}}\nabla \widetilde{g}_{i}(\bm{\lambda}_{i}^{k})\|^{2} \notag\\
\le&\frac{3(\gamma^{k})^{2}\widetilde{\kappa}}{v}\sum_{i\in \mathcal{J}}\|\nabla \widetilde{g}_{i}(\bm{\lambda}_{i}^{k})-\nabla \widetilde{g}_{i}(\bar{\bm{\lambda}}^{k})\|^{2}\notag\\
&+\frac{3(\gamma^{k})^{2}\widetilde{\kappa}}{v}\sum_{i\in \mathcal{J}}\|\nabla \widetilde{g}_{i}(\bar{\bm{\lambda}}^{k})-\frac{1}{J}\sum_{i\in \mathcal{J}}\nabla \widetilde{g}_{i}(\bar{\bm{\lambda}}^{k})\|^{2}\notag\\
&+\frac{3(\gamma^{k})^{2}\widetilde{\kappa}}{v}\sum_{i\in \mathcal{J}}\|\frac{1}{J}\sum_{i\in \mathcal{J}}\nabla \widetilde{g}_{i}(\bar{\bm{\lambda}}^{k})
-\frac{1}{J}\sum_{i\in \mathcal{J}}\nabla \widetilde{g}_{i}(\bm{\lambda}_{i}^{k})\|^{2}. \notag
\end{align}
According to Lemmas \ref{l1-appendix1} and \ref{l2-appendix1}, we have
\begin{align}
\label{lemma3-appendix1-proof-6}
&\frac{(\gamma^{k})^{2}\widetilde{\kappa}}{v}\|\nabla \widetilde{g}(\widetilde{\Lambda}^{k})-\frac{1}{J}\widetilde{\bm{1}}\widetilde{\bm{1}}^{\top}\nabla \widetilde{g}(\widetilde{\Lambda}^{k})\|_{F}^{2}\\
\le& \frac{6(\gamma^{k})^{2}\widetilde{\kappa}}{v\cdot u_{f}^{2}J^{2}}\sum_{i\in \mathcal{J}}\|\bm{\lambda}_{i}^{k}-\bar{\bm{\lambda}}^{k}\|^{2}+\frac{3(\gamma^{k})^{2}\widetilde{\delta}^{2}\widetilde{\kappa} J}{v}\notag\\
=&\frac{6(\gamma^{k})^{2}\widetilde{\kappa}}{v\cdot u_{f}^{2}J^{2}}\|\widetilde{\Lambda}^{k}-\frac{1}{J}\widetilde{\bm{1}}\widetilde{\bm{1}}^{\top}\widetilde{\Lambda}^{k}\|_{F}^{2}+\frac{3(\gamma^{k})^{2}\widetilde{\delta}^{2}\widetilde{\kappa} J}{v}.\notag
\end{align}

Substituting \eqref{lemma3-appendix1-proof-6} into \eqref{lemma3-appendix1-proof-4} and rearranging the terms, we now obtain
\begin{align}
\label{lemma3-appendix1-proof-7}
&\|\widetilde{\Lambda}^{k+1}-\frac{1}{J}\widetilde{\bm{1}}\widetilde{\bm{1}}^{\top}\widetilde{\Lambda}^{k+1}
\|^{2}_{F}\\
\le& (\frac{\widetilde{\kappa}}{1-v}+\frac{6(\gamma^{k})^{2}\widetilde{\kappa}}{v\cdot u_{f}^{2}J^{2}})\|\widetilde{\Lambda}^{k}-\frac{1}{J}\widetilde{\bm{1}}\widetilde{\bm{1}}^{\top}\widetilde{\Lambda}^{k}\|_{F}^{2}+\frac{3(\gamma^{k})^{2}\widetilde{\delta}^{2}\widetilde{\kappa} J}{v}\notag\\
=&(1-\sigma)(\frac{1}{1-v}+\frac{6(\gamma^{k})^{2}}{v\cdot u_{f}^{2}J^{2}})\|\widetilde{\Lambda}^{k}-\frac{1}{J}\widetilde{\bm{1}}\widetilde{\bm{1}}^{\top}\widetilde{\Lambda}^{k}\|_{F}^{2}+\frac{3(\gamma^{k})^{2}\widetilde{\delta}^{2}J}{v},\notag
\end{align}
where $\sigma=1- \widetilde{\kappa} \in (0,1)$. Setting $v=\frac{\sigma}{3}$ and a proper step size $\gamma^{k}$ satisfying $\frac{6(\gamma^{k})^{2}}{ u_{f}^{2}J^{2}}\le\frac{(2-\sigma)\sigma^{2}}{3(3-\sigma)}=\frac{(\sigma-v-v\sigma)v}{1-v}$, we have $\frac{1}{1-v}+\frac{6(\gamma^{k})^{2}}{v\cdot u_{f}^{2}J^{2}}\le 1+\sigma$. With this, we can rewrite \eqref{lemma3-appendix1-proof-7} as
\begin{align}
\label{lemma3-appendix1-proof-8}
&\|\widetilde{\Lambda}^{k+1}-\frac{1}{J}\widetilde{\bm{1}}\widetilde{\bm{1}}^{\top}\widetilde{\Lambda}^{k+1}
\|^{2}_{F}\\
\le &(1-\sigma^{2})\|\widetilde{\Lambda}^{k}-\frac{1}{J}\widetilde{\bm{1}}\widetilde{\bm{1}}^{\top}\widetilde{\Lambda}^{k}\|_{F}^{2}+\frac{9(\gamma^{k})^{2}\widetilde{\delta}^{2}J}{\sigma}.\notag
\end{align}
Using \eqref{lemma3-appendix1-proof-8} recursively yields
\begin{align}
\label{lemma3-appendix1-proof-9}
&\|\widetilde{\Lambda}^{k+1}-\frac{1}{J}\widetilde{\bm{1}}\widetilde{\bm{1}}^{\top}\widetilde{\Lambda}^{k+1}
\|^{2}_{F}\\
\le &(1-\sigma^{2})^{k+1}\|\widetilde{\Lambda}^{0}-\frac{1}{J}\widetilde{\bm{1}}\widetilde{\bm{1}}^{\top}\widetilde{\Lambda}^{0}\|_{F}^{2}\notag\\
&+\sum_{k'=0}^{k}(1-\sigma^{2})^{k-k'}\cdot(\gamma^{k'})^{2}\cdot\frac{9\widetilde{\delta}^{2}J}{\sigma}.\notag
\end{align}

With the same initialization $\bm{\lambda}_{i}^{0}$ for all agents $i\in \mathcal{J}$, we can rewrite \eqref{lemma3-appendix1-proof-9} as
\begin{align}
\label{lemma3-appendix1-proof-10}
&\|\widetilde{\Lambda}^{k+1}-\frac{1}{J}\widetilde{\bm{1}}\widetilde{\bm{1}}^{\top}\widetilde{\Lambda}^{k+1}
\|^{2}_{F}\\
\le &\sum_{k'=0}^{k}(1-\sigma^{2})^{k-k'}\cdot(\gamma^{k'})^{2}\cdot\frac{9\widetilde{\delta}^{2}J}{\sigma}.\notag
\end{align}
To bound $\sum_{k'=0}^{k}(1-\sigma^{2})^{k-k'}\cdot(\gamma^{k'})^{2}$ in \eqref{lemma3-appendix1-proof-10}, we define $\widetilde{y}^{k}$ as \begin{align}
\label{lemma3-appendix1-proof-11}
\widetilde{y}^{k}=\sum_{k'=0}^{k-1}(1-\sigma^{2})^{k-1-k'}\cdot (\gamma^{k'})^{2},
\end{align}
which satisfies the relation $\widetilde{y}^{k+1}=(1-\sigma^{2})\widetilde{y}^{k}+(\gamma^{k})^{2}$.

Substituting $\widetilde{y}^{k+1} = y^{k+1}$, $\psi_{1}=1-\sigma^{2}\in (0,1)$, $\psi_{2}=1\ge 0$ and $\widetilde{y}^{0}=0\le (\gamma^{0})^{2}$ into Lemma \ref{l4-appendix1}, for any integer $k \ge 0$ and the step size $\gamma^{k}$ satisfying $1\le \frac{(\gamma^{k})^{2}}{(\gamma^{k+1})^2}\le \frac{2}{1+(1-\sigma^{2})}$, we have
\begin{align}
\label{lemma3-appendix1-proof-12}
\widetilde{y}^{k+1} & =\sum_{k'=0}^{k}(1-\sigma^{2})^{k-k'}\cdot (\gamma^{k'})^{2} \\
& \le \frac{2}{1-(1-\sigma^{2})}(\gamma^{k+1})^{2} = \frac{2}{\sigma^{2}}(\gamma^{k+1})^{2}.\notag
\end{align}
With \eqref{lemma3-appendix1-proof-12}, we can rewrite \eqref{lemma3-appendix1-proof-10} as
\begin{align}
\label{lemma3-appendix1-proof-13}
\|\widetilde{\Lambda}^{k+1}-\frac{1}{J}\widetilde{\bm{1}}\widetilde{\bm{1}}^{\top}\widetilde{\Lambda}^{k+1}
\|^{2}_{F}\le \frac{18(\gamma^{k+1})^{2}\widetilde{\delta}^{2}J}{\sigma^{3}},
\end{align}
which completes the proof.
\end{lemma}
\end{appendices}

\balance

%

\end{document}